\documentclass[12pt]{article}
\usepackage{fullpage,graphicx,psfrag,amsmath,amsfonts,verbatim,cleveref}
\usepackage[small,bf]{caption}
\usepackage{url}
\usepackage{booktabs}

\usepackage{bm,bbm}
\usepackage{amsmath,amsthm,amssymb,amsfonts,cancel}
\usepackage{thmtools,enumitem,subfigure}
\usepackage{algorithm}
\usepackage{algpseudocode}

\usepackage{mathtools,dsfont}
\usepackage{caption}

\newtheorem{theorem}{Theorem}

\newcommand{\problem}{$\mathrm{SRMNS}$}

\newtheorem{lemma}{Lemma}
\newtheorem{assumption}{Assumption}
\newtheorem{claim}{Claim}
\newtheorem{definition}{Definition}
\newtheorem{proposition}{Proposition}

\usepackage{mathtools}
\newcommand{\opt}{\ensuremath{\texttt{OPT}}}
\newcommand{\obj}{\ensuremath{\texttt{OBJ}}}

\newcommand{\calH}{\mathcal{H}}

\newcommand{\halmos}{\qed}

\newcommand{\tildeN}{\ensuremath{\widetilde{N}^f}}
\newcommand{\Alg}{\ensuremath{\texttt{ALG}}}

\mathchardef\mhyphen="2D 

\newcommand{\frall}{\ensuremath{\,\forall\,}}

\DeclareMathOperator*{\argmax}{arg\,max}
\DeclareMathOperator*{\argmin}{arg\,min}

\newtheorem*{definition*}{Definition}
\newtheorem*{lemma*}{Lemma}
\newtheorem*{theorem*}{Theorem}
\newtheorem*{claim*}{Claim}

\usepackage[suppress]{color-edits}
\addauthor{kz}{magenta}
\addauthor{df}{blue}
\usepackage{bbm}
\title{Overbooking with bounded loss}
\author{Daniel Freund \and Jiayu (Kamessi) Zhao}

\begin{document}
\maketitle

\begin{abstract}
We study a classical problem in revenue management: quantity-based single-resource revenue management with no-shows. In this problem, a firm observes a sequence of~$T$ customers requesting a service. Each arrival is drawn independently from a known distribution of~$k$ different types, and the firm needs to decide irrevocably whether to accept or reject requests in an online fashion. The firm has a capacity of resources~$B$, and wants to maximize its profit.  Each accepted service request yields a type-dependent revenue and has a type-dependent probability of requiring a resource once all arrivals have occurred (or, be a \emph{no-show}). If the number of accepted arrivals that require a resource at the end of the horizon is greater than $B$, the firm needs to pay a fixed compensation for each service request that it cannot fulfill. With a clairvoyant, that knows all arrivals ahead of time, as a benchmark, we provide an algorithm with a uniform additive loss bound, i.e., its expected loss is independent of $T$. This improves upon prior works achieving \kzedit{$\Omega(\sqrt{T})$} guarantees.
\end{abstract}

\newpage
\tableofcontents
\newpage

\section{Introduction}
We study the canonical (quantity-based) single-resource revenue management problem with no-shows (\problem). In this problem we consider a firm (e.g., airline, restaurant, hotel, rental car agency) selling multiple products/fares of a single resource (e.g., seats, tables, rooms, cars) over a known, finite, time horizon of~$T$ periods. The firm's resource capacity~$B$ and the product prices are exogenously set, and the firm's objective is to maximize the total profit earned by controlling the availability of the different products over time, i.e., in each period a customer arrives to purchase a particular product, and the firm decides whether or not to \emph{accept} the customer. For each accepted customer it receives a product-dependent reward. At the end of the time horizon, each accepted customer is a \emph{no-show} with some (product-dependent) probability, and thus does not consume any resources. The firm is allowed to \emph{overbook}, i.e., it may accept more customers than it has capacity for, but it pays a fixed compensation (\emph{denied-service cost}) for every customer that was accepted, shows up, and cannot be served due to a lack of capacity. We study the usual scaling for this problem in which the number of customer arrivals~$T$ (demand) and the firm's supply~$B$ (supply) are scaled large \dfedit{while all other parameters remain constant}.   
Our main result is the first algorithm 
that achieves a uniform additive loss guarantee relative to a clairvoyant optimal algorithm that knows the arrivals \emph{a priori} (sometimes referred to as the hindsight optimum); uniform additive loss means that the loss is bounded independent of either $B$ or $T$.\footnote{Some papers \cite{vera2019bayesian,gupta2020interior,bumpensanti2018re} refer to this loss as \emph{regret}, but we avoid that terminology to avoid confusion with other literatures, e.g., the bandit literature \cite{agrawal2016efficient} that involve a different information/arrival structure.} 

\subsection{Motivation}\dfedit{
\problem\; was originally motivated by yield management for airlines for which the single resource consists of seats in a given cabin on a plane \cite{rothstein1971airline}. \cite{talluri2006theory} give an excellent overview of its history in the airline industry that also helps illustrate the mapping between the mathematical model and the real-world setting: up until the 1970s the airline industry was regulated, with commercial airlines  offering only two service options: first-class and coach-class with fares set centrally across carriers by the regulator (Civil Aeronautics Board). After deregulation airlines introduced new discounted fare classes with specific requirements/conditions, such as minimum-stay conditions, required round-trip travel, non-refundability, or the inclusion of a Saturday-night stay to segment demand into different categories. For example, a Saturday-night stay requirement would tend to filter out business travelers while not affecting as many leisure travelers. In essence, these requirements and conditions create different products for the same resource, e.g., a seat in the coach cabin, such that (i) revenue yields vary by product, and (ii) demand for different products is roughly independent --- the latter assumes, e.g., that leisure travelers are more price-sensitive and unwilling to pay the higher fare, whereas business travelers are price-insensitive and thus unwilling to fulfill the requirements to meet the lower fare. Having created such different products, the airlines' real-time decision is whether or not to allow reservations for a particular product. In practice, when an airline closes a given product to reservations, even though seats remain available on a flight, then a customer for such products would be quoted a higher fare based on the products that remain open. As demand is assumed independent across products, such a customer would not make a reservation for the higher fare, and the request is effectively rejected.
Notice that, assuming independence of demands across fares, the decision to accept (reject) demand for a specific product is mathematically equivalent to the decision to have that product be open (closed) to incoming reservation requests.
While originating in the airline industry, overbooking has since}  been adopted much more widely, including for lodging~\cite{bitran1996managing,bitran1995application}, rentals cars~\cite{geraghty1997revenue}, the restaurant industry~\cite{opentable}, and the nonprofit sector~\cite{metters1999yield}. 
Moreover, the single-resource problem specifically is often not only used in isolation, but also sometimes appears as a subproblem when solving problems with more than one resource, i.e., in \emph{network revenue management}~\cite{kunnumkal2012randomized}.  Thus, \problem\; is one of the fundamental \dfedit{admission control} problems in revenue management.

\subsection{Technical challenges and algorithmic techniques}

Without no-shows the quantity-based single-leg revenue management problem can be solved efficiently via dynamic programming: in each period, an optimal decision only requires knowledge of the remaining capacity and the remaining number of periods. More generally, when no-show probabilities are equal across products, an optimal decision requires {knowing} only the total count of admitted customers, not the respective product requested by each admitted customer. Thus, in this case as well there is no curse of dimensionality, and a dynamic program solves the problem efficiently to optimality, as was first shown by \cite{rothstein1971airline}.  
\dfedit{This reflects the requirements of Assumption 4.1 in \cite{talluri2006theory}: (i) no-show probabilities are the same for all customers, (ii) no-shows are independent across customers, (iii) no-shows are independent of the time of the reservation, and (iv) denied-service costs are the same across customers. 
Though these assumptions help obtain a tractable problem, \cite{talluri2006theory}, and later \cite{gallego2019overbooking} (see Chapter 3.3 therein), emphasize that (i) and (iv) are  restrictive assumptions in practice. Our work shows how to relax (i) with bounded loss. However, we keep assumption (iv), i.e., we assume that the denied-service costs are the same across different products; obtaining bounded loss when loosening (iv) in addition to (i) remains an interesting avenue for future research.
}

Though the case with heterogeneous no-show probabilities has been studied extensively, there are few results with provable guarantees for this problem \cite{dai2019network,kunnumkal2012randomized}. This is in contrast to the network revenue management problem without overbooking for which a number of algorithms with uniform loss guarantees have been developed over the past decade, and especially in the last few years \cite{jasin2012re, bumpensanti2018re,vera2019bayesian,sun2020near}. We fill this gap by adapting the compensated coupling technique of \cite{vera2019bayesian} \dfedit{(see proof of Theorem \ref{thm:coupling} and paragraph thereafter)} to derive the first such uniform loss guarantee for a revenue management problem with overbooking, specifically for \problem\; with heterogeneous no-show probabilities.

\dfedit{
\subsubsection*{Nonlinear objective}
The key technical challenge in adapting the compensated coupling technique lies in the fact that the denied-service costs for the expected number of no-shows gives rise to a nonlinear objective. Existing results developing/applying the compensated coupling technique~\cite{vera2019bayesian,vera2019online,freund2019uniform} rely on a property that \cite{freund2019uniform} refer to as $\delta$-insensitivity: intuitively, this means that the optimal solution is Lipschitz continuous in the number of arrivals of each type over the entire horizon; notice that this is a Lipschitz property about the optimal solution itself, not just its objective. For linear objectives, $\delta$-insensitivity follows from standard results in linear perturbation analysis~\cite{mangasarian1987lipschitz}. In contrast, establishing such a Lipschitz property directly for our setting is difficult due to the nonlinear objective.

\subsubsection*{Index policies} Leveraging compensated coupling requires us to take the following detour: instead of benchmarking directly against the complicated clairvoyant optimal solution, we define a suboptimal clairvoyant solution for which we show that it (i) is Lipschitz continuous in the number of arrivals of each type, and (ii) has bounded loss relative to the clairvoyant optimal solution. In order to define this solution we take inspiration from the greedy heuristic for the one-dimensional knapsack problem, which sorts items by their value/weight ratio and includes items in that order until no further items fit. As the expected ``weight'' of a customer is the probability that they will show up, we similarly sort customers by the ratio between their revenue and their probability of showing up (which we refer to as \emph{critical ratio}), and accept customers in that order, i.e., there exists some threshold on the critical ratio such that all types with critical ratio higher than the threshold are accepted whereas all types with critical ratio lower than the threshold are rejected. We refer to such a solution as an \emph{index solution} as the type indices are, without loss of generality, assumed to be ordered by critical ratios; for the index that has critical ratio equal to the index policy's threshold there is no restriction on how many arrivals should be accepted and we refer to that type as having the \emph{threshold index}. The main technical innovations to adapt compensated coupling then require us to establish that optimal index solutions fulfill (i) and (ii).

\subsubsection*{Lipschitz continuity and local optimality} 
For a given type, keeping the number of accepted arrivals constant for all other types, we can identify the optimal number of arrivals to accept through a marginal analysis. In particular, we should accept an additional request of that type if the expected compensation for that arrival (the compensation multiplied by the probability that they show up \emph{and} the already accepted arrivals consume all of the capacity) is smaller-equal to the revenue of that type. The optimal number of arrivals to accept for that type then follows a newsvendor-like condition   --- for the case of a single type, \cite{gallego2019overbooking} show in their Proposition 3.1 that this locally optimal solution is indeed optimal.\footnote{Considering only the offline optimal index solution, an added benefit of the newsvendor-like condition is that it allows us to find the threshold index, as well as the locally optimal number of accepted customers for that index, through binary search.} With heterogeneous types with varying no-show probabilities, optimal solutions can be more complicated: accepting an additional customer of one type may require rejecting some number of customers of another type which in turn allows for more customers of a third type\ldots. In a nutshell, the potential for such loops makes it hard to either find the globally optimal offline solution or prove $\delta$-insensitivity. In a loal sense, however, matters are not as complicated: one of our main technical insights (Lemma \ref{lem:local_sensitivity}) shows that the locally optimal number of accepted customers of one type is Lipschitz continuous in the number of accepted customers of a second type (keeping all other types fixed). 

The reasoning above allows us to derive that index solutions fulfill (i) as, roughly speaking, an increase in the accepted arrivals of one type, only affects the number of arrivals at the threshold index (other types may be similarly affected in their local optimality, but for all other types we accept either all or none of the arrivals). To derive (ii) we first show that 
the optimal solution \emph{almost} resembles an index solution (Lemma \ref{lem:chernoff_difference}). In particular, the optimal solution, for all but one type (which would correspond to the threshold index in the index solution), either accepts all (but a constant number of) arrivals or rejects all (but a constant number of) arrivals. Combined with the local Lipschitz property applied to the threshold index, this guarantees that the difference in the number of accepted arrivals between the clairvoyant optimal and the clairvoyant index solution is bounded by a constant for every type, implying in particular that the clairvoyant index solution has bounded loss relative to the clairvoyant optimal.

\subsubsection*{Instance-dependence}
The reasoning in the above paragraphs involve various instance-dependent constants, i.e., the Lipschitz continuity and the proximity between the clairvoyant index and the clairvoyant optimal solution. In particular, these depend on the magnitude of  the no-show probabilities, and the ratio of value and no-show probability. This makes sense intuitively: if one type has a probability of $\epsilon$ of showing up, then a small decrease in accepted arrivals of another type (with constant probability of showing up) should intuitively cause the locally optimal number of the former type to increase by roughly $1/\epsilon$ --- this intuition is driven again by the greedy Knapsack heuristic that replaces probabilities by expectations. When arrival probabilities are allowed to be that small, our $O(1)$-guarantees break down and, in fact, we show formally in Appendix \ref{app:instance_loss} that any online policy must incur $\Omega(\sqrt{T})$ loss when no-show probabilities are allowed to be arbitrarily small. In that same appendix we also show a shortcoming of index policies specifically: we give an example to illustrate their inability to accurately account for the trade-off between risk and reward when the critical ratios of two types are close to each other (and getting closer as the horizon gets longer). We identify an example where a very simple online policy obtains $O(1)$ loss, whereas even the clairvoyant index policy that knows the arrivals incurs $\Omega(\sqrt{T})$; intuitively, this is caused by the index policy preferring a type with infinitesimally larger critical ratio, even though that type introduces significant costs due to higher variability. For a more detailed discussion of how our guarantees depend on the problem parameters we refer the reader to the paragraphs after Lemma~\ref{lem:local_sensitivity} and Lemma \ref{lem:chernoff_difference}.
\subsection{Related work}
\subsubsection*{Compensated coupling} As we alluded to above, our work is most closely related to the recently introduced \emph{compensated coupling} technique for online stochastic decision-making problems \cite{vera2019bayesian,vera2019online,freund2019uniform}. 
The novelty of our work, relative to these, lies in (i) the structural results for \problem\; that bound the losses incurred in the detour to the index policy, (ii) the extension of the compensated coupling technique to a nonlinear objective, and (iii) the first uniform loss guarantees for a dynamic revenue management problem with overbooking.} Our work also relates closely to the following recent results: \cite{arlotto2018logarithmic,bray2019does} prove logarithmic lower bounds on the optimal loss guarantee when types have values that come from continuous distributions,~e.g.,~$U[0,1]$, rather than coming from a finite set. Similarly, \cite{freund2019uniform} prove lower bounds on uniform loss guarantees that are based on the arrival probabilities of each type (when the latter are not iid). Both of these consider special cases of our setting with overbooking, and thus apply as well.  \cite{vera2019online,chen2021linear} prove uniform loss guarantees for settings where the arrival probabilities are unknown a priori (but iid), and need to be learned --- we believe that this part of the analysis of \cite{vera2019online} would extend to our setting in a straightforward manner. Finally, it is worth highlighting the work of \cite{arlotto2019uniformly}, who introduced novel technical ideas when studying the multisecretary problem with known distribution (equivalent to our problem without no-shows/overbooking), and thereby initiated much of the recent work in this space. 

\subsubsection*{Fluid, diffusion, and uniform loss}
The OR literature often distinguishes between fluid and diffusion optimal solutions when proving asymptotic optimality for finite-horizon stochastic control. A solution is optimal on the fluid scale if its loss grows as $o(T)$ over a horizon of length $T$, and optimal on the diffusion scale if its loss grows as $o(\sqrt{T})$; of course, a uniform loss guarantees implies both. For network revenue management (no overbooking), \cite{cooper2002asymptotic} proved the fluid optimality of a static admission policy that is based on the optimization of a deterministic relaxation;   \cite{reiman2008asymptotically} found that resolving that relaxation once at an appropriately selected time, to update the policy, gives the first diffusion optimal algorithm. This was strengthened by \cite{jasin2012re} who resolve in every period to obtain a uniform loss under a technical non-degeneracy assumptions, and further strengthened by \cite{bumpensanti2018re,vera2019bayesian,freund2019uniform} to hold in the absence of such an assumption, and without resolving in every period. 
For problems with overbooking, we know of no algorithms that are optimal on the diffusion scale, as we discuss next. 


\subsubsection*{Overbooking} Our work adds to the literature on dynamic optimization problems in revenue management, and more specifically to the long line of literature on overbooking. \cite{mcgill1999revenue} describe overbooking as having the ``longest research history of any of the components of
the revenue management problem'' since most quantitative research in revenue management before 1972 focused on admission controls for overbooking. Early models of overbooking include \cite{beckmann1958decision, thompson1961statistical, rothstein1971airline}, which all focus on just a single leg. In contrast to most early works, \cite{rothstein1971airline} studied a dynamic problem, not a static one (in the language of \cite{lautenbacher1999underlying}). The distinction, informally, lies in whether admission levels are statically set initially, or made continuously over time. The solution of \cite{rothstein1971airline} is an optimal dynamic program, but it considers only a single fare. For multiple fares, the \emph{Expected Marginal Seat Revenue} heuristics (EMSRa and EMSRb), that are popular in practice to this day, were developed by \cite{belobaba1987air}. Motivated by Littlewood's rule \cite{littlewood1972forecasting}, these heuristics are fundamentally static, yet they can be applied in a dynamic setting (however, they have no provable guarantees); our use of the \emph{critical ratio} of value and no-show probability is motivated by the same ideas (see Section~\ref{sec:model} and Lemma \ref{lem:local_optimality}). 
Much more recently, \cite{aydin2013single} study the setting with multiple fares on a single-leg, and cancellations, but assume uniform no-show and cancellation probabilities; their policies do not have any performance guarantees. \cite{erdelyi2010dynamic} study the network revenue management problem with only no-shows via a decomposition of the network by flight legs, before using state aggregation for the single-leg problem, i.e., while studying the problem with heterogeneous no-show probabilities, they approximate this problem through one with homogeneous no-show probabilities, without finding provable guarantees for this approach --- in Appendix \ref{app:dpd_loss} we prove that, for the single-leg case, their approximation of heterogeneous no-show probabilities with homogeneous ones causes a loss of $\Omega(T)$; crucially this is not a shortcoming of the specific method by which they approximate heterogeneous no-show probabilities with homogeneous ones, but rather a fundamental limitation of any such approximation (see Appendix \ref{app:dpd_loss} for details). \cite{kunnumkal2012randomized} allows for no-shows in a network setting, and finds an asymptotically optimal policy that is based on the randomized LP of \cite{talluri1999randomized}. While asymptotically optimal, this approach incurs $\Omega(\sqrt{T})$ loss over a time horizon of length~$T$, whereas ours incurs~$O(1)$. 
\dfedit{Loss guarantees of that magnitude may also be obtainable using online convex optimization, e.g., as in \cite{agrawal2014fast}; notably they assume significantly less structure on the problem than we do, and therefore their problem is harder (in general, there is a significant body of work on online stochastic allocation problems without overbooking that aims to achieve sublinear regret with less problem structure, e.g., in the context of the AdWord problem \cite{alaei2012online,alaei2013online})}.
Finally, \cite{dai2019network} study the network revenue management problem with no-shows and cancellations and prove a $O(\sqrt{T})$ loss over a time horizon of length~$T$; they explicitly state as an open question whether $O(1)$ guarantees are achievable when overbooking is allowed, which we answer in the affirmative. We refer the interested reader to~\cite{gallego2019revenue} for an overview of further recent work on overbooking.

\subsubsection*{Overview of remaining sections} In Section \ref{sec:model}, we formally define our model and algorithm. In Section~\ref{sec:main} we present our main result using a compensated coupling proof that requires three new auxiliary results. The main technical difficulties lie in obtaining these auxiliary results, which we prove in Section~\ref{sec:proofs} and Appendix~\ref{app:proofs}. In Section~\ref{sec:numerical}, we complement our analytical results with numerical experiments. 

\section{Model}\label{sec:model}
\subsubsection*{Arrivals} Consider a known finite time horizon of $T$ periods. Every period begins with the arrival of a type $j\in[k]$, where $[k]$ denotes the set $\{1,\ldots,k\}$. Each arrival is of type $j$ with probability~$\lambda_j>0$, where $\sum_j\lambda_j=1$ and arrivals are independent. We denote the vector of arrivals $\vec{A}$ where $A_t$ denotes the arrival type in period $t$, and $\vec{A}[t,T]$ denotes the last $(T-t+1)$ entries of $\vec{A}$. We use $N_j$ for the number of arrivals of type $j$ over the entire time horizon and $N_j\kzedit{[}t_1,t_2\kzedit{]}$ for the number of type $j$ arrivals in periods $t_1,\ldots,t_2$, where $t_1<t_2$. When future arrivals $N\kzedit{[}t,T\kzedit{]}$ need to be estimated in period $t$, we denote such an estimate by $\tildeN\kzedit{[}t\kzedit{]}$.

\subsubsection*{Objective} There is an initial capacity $B$. Upon the arrival of each customer, we need to make an irrevocable decision on whether to accept or reject them. Accepting a customer of type~$j$ generates revenue~$v_j<1$. For an algorithm $\Alg$ we denote by $x^\Alg_j[t]$ the number of accepted customers of type~$j$ until the end of period $t$, i.e., in periods~$1,\ldots, t$. Similarly, let $\vec{x}^\Alg[t]$ denote the $k$-dimensional vector of accepted customers of each type until the end of period~$t$. At the end of the time-horizon each accepted customer of type $j$ (independently) consumes one resource with probability~$p_j \in \mathbb{Q}$, and is a no-show (does not consume any resources) with probability $1-p_j$.
The types are ordered such that $v_1/p_1\geq v_2/p_2\geq\cdots\geq v_k/p_k$, where ties are broken in favor of greater value, i.e., if~$v_j/p_j=v_{j+1}/p_{j+1}$, then $v_j>v_{j+1}$ and $p_j>p_{j+1}$. Since this ratio remains important throughout the paper, we define it as the \emph{critical ratio} $q_j = \frac{v_j}{p_j}$ of customer type $j$, and define $\bar{q}_j = 1-\frac{v_j}{p_j}$. 

In our asymptotic analysis, we assume that all parameters of the customer types are assumed to be fixed, whereas $B$ and $T$ grow large. The combined required resources (arrivals) at the end of the time horizon are distributed as $\sum_j Bin(x_j^\Alg[T],p_j)$, where $Bin(x,p)$ denotes the binomial distribution with $x$ trials and success probability $p$. When the arrivals require more than~$B$ resources, a \emph{type-independent compensation} is paid for each \emph{overbooked} customer. We normalize the compensation to $1$ so the total compensation has cost~$(\sum_jBin(x_j^\Alg[T],p_j)-B)^+$ where $(\alpha)^+:=\max\{\alpha,0\}$. Thus, the objective is to maximize
$$\mathbb{E}\left[\sum_j x_j^\Alg[T] v_j-\left(\sum_jBin(x_j^\Alg[T],p_j)-B\right)^+\right].$$

For any $j\in[k]$, we denote by $X_{j,x_j^\Alg[T]} \sim Bin(x_j^\Alg[T],p_j)$ the binomial random variable of accepted customers of type $j$ who consume a resource. Moreover, we adopt a shorthand notation for expected compensation, with $$V_B(\vec{x}) \coloneqq \mathbb{E}\left[\left(\sum_j X_{j,x_j}-B\right)^+\right].$$

\kzedit{Note that we may assume that the critical ratio $q_j < 1$, i.e, $\bar{q}_j > 0$. Customer types $j$ with $q_j \geq 1$ should always be accepted even if we had $0$ capacity (left) because the expected compensation is less than the guaranteed value from accepting, and this will indeed be the case for our policy too. Thus, we focus without loss of generality on types $j$ with $q_j<1$.}




\subsection{Algorithm and benchmark policies}

We now define four ``policies'': the online index policy, the hybrid clairvoyant index objective, the clairvoyant index objective, and the clairvoyant general objective. The first is our algorithm, the others provide useful benchmarks in our analysis. All four are based on deterministic optimization problems that take an estimate $\tildeN\kzedit{[}t\kzedit{]}$ of future arrivals in period $t$ to solve
\begin{equation}\label{eq:offline}
    \max_{\vec{x}}\sum_j v_j \left(x_j[t-1] + x_j\right) + V_B\left(\vec{x}[t-1] + \vec{x} \right) \text{ subject to } 0\leq x_j\leq \tildeN_j\kzedit{[}t\kzedit{]}, \forall j
\end{equation}
In words, assuming future arrivals $\tildeN\kzedit{[}t\kzedit{]}$ in periods $t,t+1,\ldots,T$ and past accepted arrivals $\vec{x}[t-1]$, optimization \eqref{eq:offline} finds the optimal number $\vec{x}$ to accept of each in periods $t,t+1,\ldots,T$.

\subsubsection*{Index policies}
The first three policies are all index policies, meaning intuitively that if $j'>j$, then they prioritize accepting arrivals of type $j$ over arrivals of type $j'$. Since we are in an online setting, this holds only in a forward-looking manner, i.e., it is conditioned on the already accepted customers. Index solutions in period $t$, conditioned on past decisions $\vec{x}[t-1]$, are defined as follows.
\begin{definition}\label{defn:index_policy}
Suppose in period $t$ the remaining arrivals are estimated to be $\tildeN\kzedit{[}t\kzedit{]}$. A solution~$\vec{x}$ to \eqref{eq:offline} is an \emph{index solution} if it has a \emph{threshold index}~$\widetilde{j}\in[k]$ such that for all~$j<\widetilde{j}$ \emph{all} future arrivals of type $j$ are accepted, i.e., $x_j=\tildeN_j\kzedit{[}t\kzedit{]}$, and for all~$j>\widetilde{j}$ \emph{no} future arrivals of type $j$ are accepted, i.e.,~$x_j=0$. As a convention, we always assume that when more than one index could be a threshold index, we uniquely define $\widetilde{j}$ to be the smallest possible one, implying that $x_{\widetilde{j}}>0$ (unless $x_1=0$).
\end{definition}

Notice that the intuition for an index policy comes from the one-dimensional knapsack LP in which a customer of type $j$ takes up exactly $p_j$ capacity, rather than $Bin(1,p_j)$ capacity. In that case, the critical ratio $v_j/p_j$ corresponds to the \emph{density of value} of type $j$, and the optimal LP solution is to sort by density and pack users until the knapsack is filled. Similarly, our index policy in Definition \ref{defn:index_policy} prioritizes customers with higher critical ratio over ones with lower critical ratios.

\begin{algorithm}[ht]
\caption{Online Index Policy}\label{alg:base}
\begin{algorithmic}[1]
\State Initialize $x_{ j}[0]=0\;\frall j$\;
 \State Draw arrival vector $\vec{A'}$ of length $T$ with each coordinate equal to $j$ with probability (i.i.d.) $\lambda_j$\;
 \For{$t=1,\ldots,T$}
  \State Observe type $j$ of arrival in period $t$\;
  \State $\frall j'\in[k]$ set $\tildeN_{j'}\kzedit{[}t\kzedit{]}$ equal to the number of type $j'$ arrivals in $\vec{A}'[t+1,T]+\mathbbm{1}_{\{j=j'\}}$\;
  \State  Find optimal index solution $\vec{x}'$ for \eqref{eq:offline} with inputs $\tildeN_j\kzedit{[}t\kzedit{]}$ and $\vec{x}[t-1]$\;
  \State \textbf{if}\;{$x_j'\geq N_j^f\kzedit{[}t\kzedit{]}/2$:}{
  Accept arrival of type $j$ in period $t$ and set $x_j[t]=x_j[t-1]+1$
  }
  \State \textbf{else}:\;{Reject arrival of type $j$ in period $t$ and set $x_j[t]=x_j[t-1]$}
  \State Set $x_{j'}[t]=x_{j'}[t-1],\frall j'\neq j$
 \EndFor
\end{algorithmic}
\end{algorithm}

\subsubsection*{Online index policy} The online index policy is formally given in Algorithm \ref{alg:base}. Initially, it samples an arrival vector $\vec{A}'$; this is not the real arrival vector but stems from the same distribution. In period~$t$, it observes an arrival $j'$, and uses $\vec{A}'$ to estimate $\tildeN_j\kzedit{[}t\kzedit{]}$. It then solves the optimization problem~\eqref{eq:offline} with  $\tildeN_j\kzedit{[}t\kzedit{]}$ and its own past decisions $\vec{x}[t-1]$ as input subject to finding an index solution. We denote the threshold index of that solution by~$\widetilde{j}$. In period $t$ it accepts~$j'$ if~$j'<\widetilde{j}$ or  if~$j'=\widetilde{j}$ and~$x_{\widetilde{j}}>\tildeN_j\kzedit{[}t\kzedit{]}/2$, and rejects $j'$ otherwise. 
We denote the expected (over~$\vec{A}'$) objective of the online index policy on an arrival sequence $\vec{A}=j_1,j_2,\ldots j_T$ by~$\obj_{\vec{A}}$.

\subsubsection*{Hybrid clairvoyant index objective} Consider in period $t$, with the online index policies' decision in periods $\kzedit{[}1,t-1\kzedit{]}$ the exact future arrivals $N\kzedit{[}t,T\kzedit{]}$, and solve the optimization problem~\eqref{eq:offline} with these arrivals subject to finding an index solution. The resulting objective in period $t$, for arrival sequence $\vec{A}$, is the hybrid clairvoyant index objective $\calH_{\vec{A}}\kzedit{[}t\kzedit{]}$.

\subsubsection*{Clairvoyant index objective} 
Consider the hybrid clairvoyant index objective in period $1$, i.e., before the online index policy has made any decisions. Then it solves the optimization problem~\eqref{eq:offline} subject to (i) the exact arrival sequence $\vec{A}$, and (ii) the constraint of finding an index solution. In line with the above, we denote it by $\calH_{\vec{A}}\kzedit{[}1\kzedit{]}$.

\subsubsection*{Clairvoyant general objective} 
The clairvoyant general objective is based on solving the optimization problem \eqref{eq:offline} subject to the exact arrivals $N_j$; with arrival vector $\vec{A}$ we denote it by~$\widehat{\opt}_{\vec{A}}$.

\vspace{.1in}

A few remarks are made in order to provide some intuition for the defined policies and benchmarks. First, for $t>1$, $\calH_{\vec{A}}\kzedit{[}t\kzedit{]}$ depends on $\vec{A}$ rather than only on $\vec{N}$ because the decisions of the online index policy in periods $1,\ldots,t-1$ are reflected in the value of $\calH_{\vec{A}}\kzedit{[}t\kzedit{]}$. In contrast, $\calH_{\vec{A}}\kzedit{[}1\kzedit{]}$ and  the clairvoyant general objective $\widehat{\opt}_{\vec{N}}$ do not depend on either the decisions of the online index policy or the order of arrivals. Second, the clairvoyant general solution is an optimal solution to~\eqref{eq:offline} but the clairvoyant does not know the realization of each arrival's no-show probability. Such a clairvoyant would be way too powerful to serve as a useful benchmark, as was observed by~\cite{dai2019network}. Nevertheless, our much weaker clairvoyant, which knows the realization of the future arrivals, but not the realization of the no-show probabilities, still provides an upper bound on an optimal algorithm. Finally, it may seem counter-intuitive for the algorithm to draw a sample~$\vec{A}'$ to estimate~$\tildeN\kzedit{[}t\kzedit{]}$ rather than just using the expected value thereof. However, this ensures that our estimates for~$\tildeN\kzedit{[}t\kzedit{]}$ are integral, which simplifies the analysis in light of~$Bin(x_j,p_j)$ not being well-defined for fractional $x_j$.


\subsection{Further notation}\label{ssec:further_notation}

We often rely on the Poisson Binomial distribution $$PBin(x_1,p_1;...;x_k,p_k),$$ which consists of $x_j$ trials, each with success probability $p_j$, for $j\in[k]$. Note that $\sum_j X_{j,x_j} \sim PBin(x_1,p_1;...;x_k,p_k)$ based on this definition. \kzedit{As a shorthand notation, for $X \sim PBin(x_1,p_1;...;x_k,p_k)$, we denote the compensation given $X$ by $W_B(X) = \left(X-B\right)^+$ and the expected compensation by $\mathbb{E}[W_B(X)] = \mathbb{E}\left[ \left(X-B\right)^+\right]$.} 
 
We denote the Cumulative Distribution Function (CDF) of any discrete random variable $X$ at~$x$ by ~$F_X (x)$, the Probability Mass Function (PMF) by $\mathbb{P}[X=x]$, and the quantile function (i.e., the inverse of CDF) by $F_X^{-1} (q)$, where $q \in [0,1]$. We use the notation $\mathbb{P}[E]$ for the probability of event $E$. And lastly, we use the standard notation $\mathbb{N} = {0,1,2,...}$ for natural numbers, and $e_j$ for the $k$-dimensional canonical unit vector with its $j^{th}$ coordinate equal to 1 and all other coordinates equal to 0.

\section{Main result}\label{sec:main}

Our goal is to prove that the expected additive loss between the online index policy and the clairvoyant objective can be bounded by a constant~$M$, i.e., $\mathbb{E}_{\vec{A}}\left[\widehat{\opt}_{\vec{A}}-\obj_{\vec{A}}\right] \leq~M$.\footnote{Throughout this work, we refer to terms as constant if they do not depend on $B$ or $T$.} 
To do so, we derive intermediate benchmarks and bound their expected loss relative to each other. The most important technical ingredients in our analysis use the following \emph{local} optimality condition:
\kzedit{
\begin{definition}\label{defn:local_opt}
Consider a feasible solution $\vec{x}^\star$ to \eqref{eq:offline} with (estimated) future arrivals $\tildeN\kzedit{[}t\kzedit{]}$ and past accepted arrivals $\vec{x}[t-1]$. We call $\vec{x}^\star$ \emph{locally optimal at} $j$ if 
$$
x_j^\star=\max\argmax_{x_j:0\leq x_j \leq \tildeN_j\kzedit{[}t\kzedit{]}}\left\{
x_jv_j-
\mathbb{E}\left[\left(X_{j,x_j[t-1]+x_j}+\sum_{n\neq j} X_{n,x_n[t-1]+x_n^\star}-B\right)^+\right] \right\}.
$$
\end{definition}
Intuitively, $\vec{x}^\star$ is locally optimal at $j$ if $x_j^\star$ is the (maximum) optimal number of future type $j$ customers to accept given acceptance of $x_n^\star$ future customers of type $n$ for each $n\neq j$.} Because this maximum optimal number of customers is uniquely defined, we observe the irritating property that there may be global optima, i.e., solutions that yield the highest possible profit, that do not satisfy our local optimality condition for every type $j$. Fortunately, for our analysis, there always exists some global optimum that is locally optimal for every type $j$. This claim is formalized and proven in Appendix \ref{app:glocal}, and it allows us to compare our solution to a global optimum that fulfills the local optimality condition for every type. Next, in Lemma \ref{lem:local_optimality} we provide an equivalent local optimality condition that is often easier to apply.

\begin{lemma}\label{lem:local_optimality}
For a solution $\vec{x}^\star$ to \eqref{eq:offline} with (estimated) future arrivals $\tildeN\kzedit{[}t\kzedit{]}$ and past accepted arrivals ~$\vec{x}[t-1]$, let $\vec{x} = \vec{x}[t-1] + \vec{x}^\star$. \kzedit{Then, $\vec{x}^\star$ is locally optimal at $j$ if and only if both of the following conditions are satisfied:
\begin{enumerate}[label=(\roman*)]
    \item $x_j^\star = \tildeN_j\kzedit{[}t\kzedit{]}$ or $\mathbb{P}\left[\sum_{n} X_{n,x_n} \geq B\right] > q_j$;
    \item $x_j^\star = 0$ or $\mathbb{P}\left[\sum_{n \neq j} X_{n,x_n} + X_{j,x_j -1} \geq B\right] \leq q_j$.
\end{enumerate}}
\end{lemma}

\kzedit{As noted in the literature review, }Lemma \ref{lem:local_optimality} bears resemblance to classical results in the field of revenue management like Littlewood's rule \cite{littlewood1972forecasting}, the ESMR heuristics \cite{belobaba1987air}, and most notably Proposition 3.1 in \cite{gallego2019overbooking}, which is equivalent for the case that $k=1$.

In a nutshell, the main analytical advantage of local optimality over global optimality is that we can derive the following (local) Lipschitz bound on locally optimal solutions. We apply this repeatedly in our proofs to bound by how much solutions deviate from each other. 
\begin{lemma}\label{lem:local_sensitivity}
There exists a constant $\delta \geq 1$ such that for any solution $\vec{x}^\star$ to \eqref{eq:offline} with $\tildeN\kzedit{[}t\kzedit{]}$ and ~$\vec{x}[t-1]$, it is true that if $\vec{x}^\star$ is locally optimal at $j$, then $\exists l \in \{0,1,...,\delta\}$ such that $\vec{x}^\prime = \vec{x}^\star + e_i - l e_j$ is locally optimal at $j$.
\end{lemma}

\kzedit{We remark that $\delta$ is large only when (a) $p_j$ is small or (b) the critical ratio of any type is either close to 0 or close to 1. For (a), it makes sense that $\delta$ depends inversely on $p_j$, e.g., when $p_j = 1/10,$ then for every unit of additional capacity one would want to accept (roughly) 10 more arrivals of type~$j$. For (b), $\delta$ grows large at the boundary of the critical ratios due to us resorting to the standard normal distribution for a probability bound (see Section \ref{sec:local_analysis}). The resulting bound we obtain is based on the inverse of the CDF of the standard normal, which has a large derivative (only) close to $0$ and~$1$. In our context, assuming $p_j$ to be an actual constant, e.g., in $(0.2,1]$, this corresponds to types that either have very small value (so we may want to reject them regardless) or types that have very large value (so we may want to accept them regardless). Recall also that if the critical ratio was greater than $1$ in our setting, then one should accept an arrival even if one is out of capacity as the revenue received is greater than the expected compensation to be paid (given the probability of the customer showing up).}

The key motivation for our consideration of index solution and local optimality is that Lemma~\ref{lem:local_sensitivity} provides a strong bound to capture how an optimal index solution changes when an additional arrival of some type is accepted. In contrast, for a globally optimal solution, we do not know how to prove such a bound. Next, armed with these lemmas, we can show the following two bounds.

\begin{theorem}\label{thm:index_loss}
There exists a constant $M_1$ such that for every arrival vector $\vec{A}$
$$
\widehat{\opt}_{\vec{A}}-\calH_{\vec{A}}\kzedit{[}1\kzedit{]}
\leq M_1,
$$
which implies in particular that $\mathbb{E}_{\vec{A}}\left[
\widehat{\opt}_{\vec{A}}\right]\leq M_1 +\mathbb{E}_{\vec{A}}\left[ \calH_{\vec{A}}\kzedit{[}1\kzedit{]}\right]$.
\end{theorem}

\begin{lemma}\label{lem:small_loss_probability}
There exists a constant $M_2$ such that
$$
\sum_{t=2}^T \mathbb{P}\left[\calH_{\vec{A}}\kzedit{[}t-1\kzedit{]}-\calH_{\vec{A}}\kzedit{[}t\kzedit{]}>0\right]\leq M_2.
$$
\end{lemma}

From these bounds we derive, using compensated coupling \cite{vera2019bayesian,freund2019uniform}, the result in Theorem \ref{thm:coupling}.
\begin{theorem}\label{thm:coupling}
There exists a constant $M$ such that $\mathbb{E}_{\vec{A}}\left[\widehat{\opt}_{\vec{A}}-\obj_{\vec{A}}\right] \leq M$.
\end{theorem}
\emph{Proof. }
\begin{align}
\mathbb{E}_{\vec{A}}\left[\widehat{\opt}_{\vec{A}}-\obj_{\vec{A}}\right] \notag\\
\leq M_1 + \mathbb{E}_{\vec{A}}\left[\calH_{\vec{A}}\kzedit{[}1\kzedit{]}-\obj_{\vec{A}}\right] \tag{Theorem \ref{thm:index_loss}}\\
= M_1 + \mathbb{E}_{\vec{A}}\left[\sum_{t=2}^T\calH_{\vec{A}}\kzedit{[}t-1\kzedit{]}-\calH_{\vec{A}}\kzedit{[}t\kzedit{]}\right] \notag\\
\leq M_1 + \sum_{t=2}^T \mathbb{P}\left[\calH_{\vec{A}}\kzedit{[}t-1\kzedit{]}-\calH_{\vec{A}}\kzedit{[}t\kzedit{]}>0\right] \tag{$\star$}\\
\leq M_1 + M_2 \tag{Lemma \ref{lem:small_loss_probability}}.
\end{align}
For the equality above, observe that we have a telescoping sum in which the first term is  $\calH_{\vec{A}}\kzedit{[}1\kzedit{]}$, and the last term is~$\calH_{\vec{A}}\kzedit{[}T\kzedit{]}=\obj_{\vec{A}}$ \kzedit{as (i) the hybrid clairvoyant in period $T$ is bound to make the same decisions as the algorithm in periods $1,\ldots, T-1$ and (ii) the online index policy makes an optimal decision in period $T$ (when there is no uncertainty about future periods)}. The starred inequality holds because any action in one period affects the objective, through lost revenue or incurred compensation, by at most $1$. \halmos

Notice that the proof is based on a sample path-wise coupling of the online index policy and the clairvoyant objective.
Then, in each period, when the hybrid clairvoyant loses in objective due to a decision by the online policy, the technique ``compensates" the clairvoyant for that loss. The expected loss of the online index policy can then be interpreted as the total compensation paid out to the clairvoyant; this, in turn, can be upper-bounded (using that the loss/compensation in any period is bounded by 1) by the sum of the probabilities of having to pay out a compensation.
%


\section{Proofs of auxiliary results}\label{sec:proofs}
This section is dedicated to proving the auxiliary results from Section \ref{sec:main}, \kzedit{and is organized as follows: in Section \ref{sec:local_analysis} we prove Lemma \ref{lem:local_optimality} and Lemma \ref{lem:local_sensitivity} in order to characterize locally optimal solutions and provide a Lipschitz bound for them. Then, in Section~\ref{sec:thm_proof} we prove Theorem~\ref{thm:index_loss} through a series of lemmas: in Lemma \ref{lem:exchange} we show through an exchange argument that for $i<j$ any feasible solution improves when a certain number of type~$j$ customers is replaced by a proportionate number of type-$i$ customers (where the proportion is based on $p_i,p_j$). Lemma \ref{lem:chernoff_difference} and \ref{lem:global_sensitivity} apply this result to show that globally optimal solutions are, in a sense, \emph{not too far} from being index solutions. As the objective of optimization \eqref{eq:offline} is Lipschitz in its solution, Theorem \ref{thm:index_loss} then follows from the globally optimal solution and the optimal index solution not being too different. Finally, in Section \ref{sec:lemma3} we combine the Lipschitz bound on optimal index solutions with the compensated coupling technique to uniformly bound the loss of the online index policy relative to the optimal index solution.
}

\subsection{Proof of Lemma \ref{lem:local_optimality} and Lemma \ref{lem:local_sensitivity}} \label{sec:local_analysis}

\subsubsection*{Proof of Lemma \ref{lem:local_optimality}}

We show both directions by contradiction. First, suppose $\vec{x}^\star$ fulfills Definition~\ref{defn:local_opt}, but~(i) is false. Recalling that the lemma defines $\vec{x} = \vec{x}[t-1] + \vec{x}^\star$, we find that
$$\mathbb{P}\left[\sum_{n} X_{n,x_n} \geq B\right] \leq q_\kzedit{j} \text{ and } x_\kzedit{j}^\star < \tildeN_\kzedit{j}\kzedit{[}t\kzedit{]}.$$
Then $\vec{x}^\star + e_\kzedit{j}$ is feasible because $x_\kzedit{j}^\star + 1 \leq \tildeN_\kzedit{j}\kzedit{[}t\kzedit{]}$; we shall argue that accepting $x_\kzedit{j}^\star+1$ customers of type~$\kzedit{j}$ yields no less profit than accepting $x_\kzedit{j}^\star$, which contradicts that  $\vec{x}^\star$ fulfills Definition~\ref{defn:local_opt}. 

We first observe that the increase in revenue for accepting one more type $\kzedit{j}$ customer is $v_\kzedit{j}$. To quantify the increase in expected compensation from taking solution $\vec{x}^\star + e_\kzedit{j}$ rather than $\vec{x}^\star$, let~$\sum_{n} X_{n,x_n}$ denote the Poisson Binomial random variable capturing the number of customers among $\vec{x} = \vec{x}[t-1]+ \vec{x}^\star$ who consume a resource (see definition in Section~\ref{ssec:further_notation}). Similarly, the random number of customers requiring resources under solution $\vec{x} + e_\kzedit{j}$ is $\sum_{n \neq \kzedit{j}} X_{n,x_n} + X_{\kzedit{j},x_\kzedit{j}+1}$. 

To characterize the compensation paid, suppose first that $\sum_{n} X_{n,x_n} < B$, i.e., the number of accepted customers among $\vec{x}$ who show up is smaller than the capacity~$B$. Then from $X_{\kzedit{j},1}\leq 1$ we also find $\sum_{n \neq \kzedit{j}} X_{n,x_n} + X_{\kzedit{j},x_\kzedit{j}+1} = \sum_{n} X_{n,x_n} + X_{\kzedit{j},1} \leq B$, and no compensation is paid for the additional type $\kzedit{j}$ customer. 

If instead $\sum_{n} X_{n,x_n} \geq B$, then the additional type $\kzedit{j}$ customer decreases the objective by 1 if that customer shows up (with probability $p_\kzedit{j}$). Thus, the increase in expected compensation from an additional type $\kzedit{j}$ customer is
$$V_B\left(\vec{x} + e_\kzedit{j}\right) - V_B\left(\vec{x}\right) = \mathbb{P}\left[\sum_{n} X_{n,x_n} < B\right] \cdot 0 + \mathbb{P}\left[\sum_{n} X_{n,x_n} \geq B\right] \cdot p_\kzedit{j} \cdot 1\leq q_\kzedit{j}\cdot p_\kzedit{j}=v_\kzedit{j}.$$
Notice that the inequality holds when Lemma \ref{lem:local_optimality} (i) is false. 
Thus, $\vec{x}^\star + e_\kzedit{j}$ yields at least the same objective value as $\vec{x}^\star$ does. Since the solution $\vec{x}^\star $ is locally optimal at $\kzedit{j}$ only if $x_\kzedit{j}^\star$ is the largest number that achieves the highest objective value, $\vec{x}^\star$ is not locally optimal at $\kzedit{j}$, which is a contradiction. 

Next, we show that Definition \ref{defn:local_opt} implies Lemma \ref{lem:local_optimality} (ii). We suppose for sake of contradiction that~$\vec{x}^\star$ is locally optimal at $\kzedit{j}$ and (ii) is false. We then show that $\vec{x}^\star - e_\kzedit{j}$ must achieve a higher objective than $\vec{x}^\star$ to derive a contradiction. First, when (ii) is false, we have

$$\mathbb{P}\left[\sum_{n \neq j} X_{n,x_n} + X_{\kzedit{j},x_\kzedit{j}-1} \geq B\right] > q_\kzedit{j} \text{ and } x_\kzedit{j}^\star > 0.$$
Then, $\vec{x}^\star - e_\kzedit{j}$ is feasible because $x_\kzedit{j}^\star - 1 \geq 0$. As above, the increase in expected compensation from taking solution $\vec{x}^\star$ rather than $\vec{x}^\star - e_\kzedit{j}$ is
\begin{align*}
    &V_B\left(\vec{x}\right) - V_B\left(\vec{x} - e_\kzedit{j}\right)\\
    =& \mathbb{P}\left[\sum_{n \neq \kzedit{j}} X_{n,x_n} + X_{\kzedit{j},x_\kzedit{j}-1} < B\right] \cdot 0 + \mathbb{P}\left[\sum_{n \neq \kzedit{j}} X_{n,x_n} + X_{\kzedit{j},x_\kzedit{j}-1} \geq B\right] \cdot p_\kzedit{j}> v_\kzedit{j},
\end{align*}
where the inequality holds when (ii) is false. Since the increase in revenue from accepting $x_\kzedit{j}^\star$ rather than $x_\kzedit{j}^\star-1$ type $\kzedit{j}$ customers is $v_\kzedit{j}$, $\vec{x}^\star - e_\kzedit{j}$ yields a better objective value than $\vec{x}^\star$ does, which contradicts that $\vec{x}^\star$ is locally optimal at $\kzedit{j}$. 



Next we show the other direction, i.e., that $\vec{x}^\star$ is locally optimal at $\kzedit{j}$ if both (i) and (ii) are satisfied. We again argue by contradiction, and show that if $\vec{x}^\star$ is not locally optimal at $\kzedit{j}$, then at least one of (i) and (ii) must be false. Specifically, if $\vec{x}^\star$ is not locally optimal at $\kzedit{j}$, then there exists some $l \neq 0$ such that $\vec{x}^\star + l \cdot e_\kzedit{j}$ is locally optimal at $\kzedit{j}$. 

If $l > 0$, we must have $x_\kzedit{j}^\star < \tildeN_\kzedit{j}\kzedit{[}t\kzedit{]}$ since otherwise $x_\kzedit{j}^\star + l > \tildeN_\kzedit{j}\kzedit{[}t\kzedit{]}$ is infeasible. Since we showed that $\vec{x}^\star + l \cdot e_\kzedit{j}$ being locally optimal at $\kzedit{j}$ requires $\vec{x}^\star + l \cdot e_\kzedit{j}$ to fulfill condition (ii) of the lemma,
$$\mathbb{P}\left[\sum_{n \neq \kzedit{j}} X_{n,x_n} + X_{\kzedit{j},x_\kzedit{j}+l-1} \geq B\right] \leq q_\kzedit{j}$$
holds true. Then, since $l \geq 1$,
$$\mathbb{P}\left[\sum_{n = 1}^k X_{n,x_n} \geq B\right] \leq \mathbb{P}\left[\sum_{n \neq \kzedit{j}} X_{n,x_n} + X_{\kzedit{j},x_\kzedit{j}+l-1} \geq B\right] \leq q_\kzedit{j}.$$
Thus, $\mathbb{P}\left[\sum_{n} X_{n,x_n} \geq B\right] \leq q_\kzedit{j}$ and $x_\kzedit{j}^\star < \tildeN_\kzedit{j}\kzedit{[}t\kzedit{]}$, which contradicts (i).

If $l < 0$, we must have $x_\kzedit{j}^\star > 0$ since otherwise $x_\kzedit{j}^\star + l < 0$ is infeasible. Since we showed that~$\vec{x}^\star + l \cdot e_\kzedit{j}$ being locally optimal at $\kzedit{j}$ requires $\vec{x}^\star + l \cdot e_\kzedit{j}$ to fulfills condition (i) of the lemma,
$$\mathbb{P}\left[\sum_{n \neq \kzedit{j}} X_{n,x_n} + X_{\kzedit{j},x_\kzedit{j}+l} \geq B\right] > q_\kzedit{j}$$
holds true. 
Then, since $l\leq -1$,
$$\mathbb{P}\left[\sum_{n \neq \kzedit{j}} X_{n,x_n} + X_{\kzedit{j},x_\kzedit{j}-1} \geq B\right] \geq \mathbb{P}\left[\sum_{n \neq \kzedit{j}} X_{n,x_n} + X_{\kzedit{j},x_\kzedit{j}+l} \geq B\right] > q_\kzedit{j}.$$
Thus, $\mathbb{P}\left[\sum_{n \neq \kzedit{j}} X_{n,x_n} + X_{\kzedit{j},x_\kzedit{j}-1} \geq B\right] > q_\kzedit{j}$ and $x_\kzedit{j}^\star > 0$, which contradicts (ii).

Thus, if $\vec{x}^\star$ is not locally optimal at $\kzedit{j}$, then at least one of (i) and (ii) must be false. \halmos

\subsubsection*{Proof of Lemma \ref{lem:local_sensitivity}}

Let $\vec{x}$ = $\vec{x}[t-1] + \vec{x}^\star$. Lemma \ref{lem:local_optimality} (i) guarantees that
$$\mathbb{P}\left[\sum_{n} X_{n,x_n} \geq B\right] > q_j, \text{ or } x_j^\star = \tildeN_j\kzedit{[}t\kzedit{]}.$$
If $x_j^\star = \tildeN_j\kzedit{[}t\kzedit{]}$, then no additional type $j$ customers may be accepted; otherwise, for any $l \leq -1,$
$$\mathbb{P}\left[\sum_{n \neq i,j} X_{n,x_n} + X_{i,x_i+1} + X_{j,x_j+|l|} \geq B\right] \geq \mathbb{P}\left[\sum_{n} X_{n,x_n} \geq B\right] > q_j$$
means that the objective value of $\vec{x}+e_i+|l| \cdot e_j$ is smaller than that of $\vec{x}+e_i+( |l|-1) \cdot e_j$. Thus,~$l \geq 0$.

We show for each $j$ that there is a constant $\delta_j$ such that $l \leq \delta_j$ is guaranteed. Then, we set~$\delta = \max_j \delta_j$ as a constant that fulfills the lemma. To determine $\delta_j$ and prove the lemma, we first take Claim \ref{cl:berryesseen} and Claim \ref{cl:exchange} below for granted. Recall that $\bar{q}_j = 1-\frac{v_j}{p_j}$ from Section \ref{sec:model}. We denote for a given constant~$\delta^\star \geq 1$ and solution $\vec{x}$: $$Y = \sum_{n \neq j} X_{n,x_n}+X_{j,x_j+\delta^\star}, Z = \sum_{n} X_{n,x_n}, \text{ and } m = \sum_{n = 1}^k x_n.$$

\begin{claim}\label{cl:berryesseen}
There exist constants $\delta_{j}^\star$ and $m_0 \in \kzedit{\mathbb{Z}^+}$ such that
\begin{equation}\label{diff}
    F_{Y}^{-1} (\bar{q}_j)-F_{Z}^{-1} (\bar{q}_j) \geq 1
\end{equation}
holds whenever $m \geq m_0$.
\end{claim}


\begin{claim}\label{cl:exchange}
If $\vec{x} - (\delta_j + 1)e_j + e_i$ is a feasible solution to the optimization problem \eqref{eq:offline}, then
\begin{equation*}
    \mathbb{P}\left[\sum_{n \neq i,j} X_{n,x_n}+X_{j,x_j-\delta_{j}-1} \geq B-X_{i,x_i+1}\right] \leq \mathbb{P}\left[\sum_{n \neq i,j} X_{n,x_n}+X_{j,x_j-\delta_{j}-1} \geq B-1-X_{i,x_i}\right].
\end{equation*}
\end{claim}
Now, take $\delta_j = 1+\delta_{j}^\star+m_0$ with $\delta_{j}^\star$ and $m_0$ from Claim \ref{cl:berryesseen}. Observe that we only need to check cases where $x_j^\star > \delta_j$, because $x_j^\prime = x_j^\star - l \geq 0$ guarantees $l \leq \delta_j$ when $x_j^\star \leq \delta_j$. Thus, we below assume that $x_j^\star > \delta_j$, when showing $l \leq \delta_j$.

For solution $\vec{x}$ with $x_j^\star > \delta_j$, define $$\bar{Y} = \sum_{n \neq j} X_{n,x_n}+X_{j,x_j-1} \text{ and } \bar{Z} = \sum_{n \neq j} X_{n,x_n}+X_{j,x_j-\delta^\star-1}.$$ In particular, $\bar{Y}$ is the random variable $Y$ from above, constructed for solution $\vec{x} - (\delta^\star+1) e_j$, which is guaranteed to be feasible when $x_j^\star > \delta_j > 1+\delta^\star$. Similarly, $\bar{Z}$ is the random variable $Z$ from above, constructed for solution $\vec{x} - (\delta^\star+1)e_j$. Since $$\sum_n x_n - (\delta^\star+1) \geq x_j^\star - (\delta^\star+1) > m_0,$$ we can apply Claim \ref{cl:berryesseen} to $\bar{Y}$ and $\bar{Z}$.

We will argue that $\forall a \geq 1:\mathbb{P}\left[\sum_{n \neq i,j} X_{n,x_n} + X_{i,x_i+1} + X_{j,x_j-\delta^\star-a} \geq B\right]\geq q_j$ to derive a contradiction from Lemma \ref{lem:local_optimality} (i). Lemma \ref{lem:local_optimality} (ii) guarantees that, with $\vec{x}^\star$ locally optimal at $j$, $x_j^\star>0$,
$$\mathbb{P}\left[\bar{Y} \geq B\right] = \mathbb{P}\left[\sum_{n \neq j} X_{n,x_n} + X_{j,x_j-1} \geq B\right] \leq q_j,$$ 
so we know $F_{\bar{Y}} (B-1) \geq \bar{q}_j$. By inverting the CDF into the quantile function, we equivalently find that $F_{\bar{Y}}^{-1} (\bar{q}_j) \leq B-1.$ Then, Claim \ref{cl:berryesseen} implies that
$$F_{\bar{Z}}^{-1} (\bar{q}_j) \leq F_{\bar{Y}}^{-1} (\bar{q}_j) - 1 \leq B-2.$$
By inverting the quantile function for $\bar{Z}$ back into the CDF, we find that $F_{\bar{Z}} (B-2) \geq \bar{q}_j$. Thus, $$\mathbb{P}\left[\sum_{n \neq j} X_{n,x_n}+X_{j,x_j-\delta^\star-1} \geq B-1\right] = \mathbb{P}\left[\bar{Z} \geq B-1\right] \leq q_j.$$

With Claim \ref{cl:exchange}, we obtain
$$\mathbb{P}\left[\sum_{n \neq i,j} X_{n,x_n}+X_{j,x_j-\delta^\star-1} \geq B-X_{i,x_i+1}\right] \leq \mathbb{P}\left[\sum_{n \neq j} X_{n,x_n}+X_{j,x_j-\delta^\star-1} \geq B-1\right] \leq q_j,$$ 

and thus
\begin{align*}
    \forall a \geq 1:\; &\mathbb{P}\left[\sum_{n \neq i,j} X_{n,x_n} + X_{i,x_i+1} + X_{j,x_j-\delta^\star-a} \geq B\right]\\
    \leq& \mathbb{P}\left[\sum_{n \neq i,j} X_{n,x_n} + X_{i,x_i+1} + X_{j,x_j-\delta^\star-1} \geq B \right]\\
    =& \mathbb{P}\left[\sum_{n \neq i,j} X_{n,x_n} + X_{j,x_j-\delta^\star-1} \geq B-X_{i,x_i+1} \right] \leq q_j.
\end{align*}

It follows from Lemma \ref{lem:local_optimality} (i) that $\vec{x}^\prime = \vec{x}^\star + e_i - l e_j$ is not locally optimal at $j$ when $$x_j^\prime \leq x_j^\star - \delta^\star - 1 < x_j^\star - \delta_j - 1.$$ Thus, there exists $l \leq \delta_j$ such that $\vec{x}^\star + e_i - l e_j$ is locally optimal at $j$. \halmos

We prove Claims \ref{cl:berryesseen} and \ref{cl:exchange} in Appendix \ref{app:cl_berry}.





\subsection{Proof of Theorem \ref{thm:index_loss}}\label{sec:thm_proof}
The proof of Theorem \ref{thm:index_loss} is based on Lemma \ref{lem:local_sensitivity} as well as a second bound. Informally, the second bound implies that, for large enough $T$, the optimal policy starts to resemble an index policy up to constant differences. We deduce that the optimal clairvoyant objective and the optimal clairvoyant index objective are based on solutions to \eqref{eq:offline} that have bounded distance in infinity-norm. Since the objective is Lipschitz, the theorem follows. We formalize this argument through three lemmas. 

\begin{lemma}\label{lem:exchange}
For every $i < j$ there exists a constant $R_{ij}$ such that, for any feasible solutions $\vec{x}^\star+\frac{R_{ij}}{p_i} \cdot~e_i$ and $\vec{x}^\star+\frac{R_{ij}}{p_j} \cdot e_j$ to optimization problem \eqref{eq:offline} with $N\kzedit{[}t,T\kzedit{]}$ and $\vec{x}[t-1]$, the objective of $\vec{x}^\star+\frac{R_{ij}}{p_i} \cdot e_i$ is always greater than that of $\vec{x}^\star+\frac{R_{ij}}{p_j} \cdot e_j$.
\end{lemma}
We prove Lemma $\ref{lem:exchange}$ via an exchange argument. 

\emph{Proof of Lemma \ref{lem:exchange}. } For any constant $R \in \mathbb{Z}^+$ such that $\frac{R}{p_i}$ and $\frac{R}{p_j}$ are both integers, denote $$X_i \sim Bin\left(\frac{R}{p_i},p_i\right), X_j \sim Bin\left(\frac{R}{p_j},p_j\right).$$ Since we leave all customers among $\vec{x} = \vec{x}^\star + \vec{x}[t-1]$ unchanged, we denote the number of unchanged customers who consume a resource by random variable 
$$X \sim PBin\left(x_1,p_1;...;x_k,p_k\right).$$

To prove the lemma we distinguish between the following cases for $i<j$. We prove the Lemma for the first case here, and the second (more complicated) case in Appendix \ref{app:chernoff}.
\begin{enumerate}[label=(\roman*)]
    \item $q_i > q_j$
    \item $q_i = q_j$ and $v_i > v_j$
\end{enumerate}

In case (i), the difference in revenue between any feasible solutions $\vec{x}^\star+\frac{R}{p_i} \cdot e_i$ and $\vec{x}^\star+\frac{R}{p_j} \cdot e_j$ is 
$$\frac{R}{p_i} v_i  - \frac{R}{p_j} v_j = \left(\frac{v_i}{p_i}-\frac{v_j}{p_j}\right)R,$$
which is positive and scales linearly with $R$. Thus, to show that a constant $R$ exists such that the increase in revenue outweighs the increase in compensations, it suffices to verify that the difference in expected compensation is bounded by $o(R)$.

\kzedit{Recall that, by definition, $\mathbb{E}[W_B(X)] = \mathbb{E}\left[ \left(X-B\right)^+\right]$ for random variable~$X$.} We show $\mathbb{E}\left[(X_i-X_j)^+\right]\in O\left(\sqrt{R\log(R)}\right)$ to bound the difference in expected compensation as 
\begin{align*}
    &\kzedit{\mathbb{E}[W_B(X+X_i)] - \mathbb{E}[W_B(X+X_j)]}\\
    \leq& \sum_{x} \mathbb{P}\left[X = x\right] \mathbb{E}\left[(X_i-X_j)^+|X = x\right]\\
    =& \mathbb{E}\left[(X_i-X_j)^+\right]\in O\left(\sqrt{R\log(R)}\right)
\end{align*}

Next, for $\mathbb{E}\left[(X_i-X_j)^+\right]$, we apply a Chernoff bound. Let $\kappa=2 \sqrt{R\log(R)}$. Since

\begin{equation*}
\begin{split}
& \mathbb{E}\left[(X_i-X_j)^+\right]\\
\leq& \mathbb{P}\left[X_i-X_j < \kappa\right] \cdot \kappa + \mathbb{P}\left[X_i-X_j \geq \kappa\right] \cdot \kzedit{\mathbb{E}[\left(X_i-X_j\right)|X_i-X_j \geq \kappa]}\\
\leq& \mathbb{P}\left[X_i-X_j < \kappa\right] \cdot \kappa + \mathbb{P}\left[X_i-X_j \geq \kappa\right] \cdot \frac{R}{p_i},\\
\end{split}
\end{equation*}
it suffices to show $\mathbb{P}\left[X_i-X_j \geq \kappa\right] \in O\left(\sqrt{\frac{\log(R)}{R}}\right)$ to prove the required bound. For~$X_i-X_j \geq \kappa$ to hold, observe that at least one of 
1) $E_1 = \left\{X_i \geq \kappa/2\right\}$; or 2) $E_2 = \left\{X_j \leq \kappa/2 \right\}$ must occur, so 
$$\mathbb{P}\left[X_i-X_j \geq \kappa\right] \leq \max\left(\mathbb{P}\left[E_1\right], \mathbb{P}\left[E_2\right]\right).$$ 
With $\epsilon = \sqrt{\frac{\log(R)}{R}}$ a Chernoff bound gives 
$$\mathbb{P}\left[E_1\right] = \mathbb{P}\left[X_i \geq R + \sqrt{R\log(R)}\right] \leq e^{-\frac{\epsilon^2}{2+\epsilon} R} = {R}^{-\frac{1}{2+\epsilon}} \in o\left(\sqrt{\frac{\log(R)}{R}}\right),$$ 
and a similar Chernoff bound works for $\mathbb{P}\left[E_2\right]$. Thus, we find 
$$\mathbb{E}\left[(X_i-X_j)^+\right] \leq O\left(\sqrt{R\log(R)}\right).$$
Therefore, there exists some $R_{ij}$ such that the increase in revenue from replacing $\frac{R_{ij}}{p_j}$ type $j$ customers with $\frac{R_{ij}}{p_i}$ type $i$ customers  outweighs the loss in additional compensation, i.e.,~$\vec{x}^\star+\frac{R_{ij}}{p_i} \cdot e_i$ yields a better objective than $\vec{x}^\star+\frac{R_{ij}}{p_j} \cdot e_j$ does, and this completes the first part of the proof. \halmos 

\begin{lemma}\label{lem:chernoff_difference}
Consider an optimal solution $\vec{x}^\star$ to the optimization problem \eqref{eq:offline} with future arrivals~$N\kzedit{[}t,T\kzedit{]}$, and past accepted arrivals $\vec{x}[t-1]$. Then for every $i$ and $j$ with $i<j$ there exists a constant $R_{ij}$ such that at least one of the following two is true:
\begin{enumerate}[label=(\roman*)]
    \item $x_i^\star > N\kzedit{[}t,T\kzedit{]}-R_{ij}/p_i$ \item $x_j^\star < R_{ij}/p_j$.
\end{enumerate}
\end{lemma}

\kzedit{Notice that the constant $R_{ij}$ above is large when types have critical ratios that are close to each other. We highlight the intuition here and this limitation of the index policy is further explored in Appendix \ref{app:eg_index}. Intuitively, we find $R_{ij}$ as the smallest $R$ such that $$\left(\frac{v_i}{p_i} - \frac{v_j}{p_j}\right) R \geq 2 \sqrt{R \log(R)},$$ where the left-hand side is the increase in revenue (from replacing $\frac{R}{p_i}$ type~$i$ customers with $\frac{R}{p_j}$ type~$j$ customers) and the right-hand side is an upper bound for the increase in compensation. Thus, $R_{ij}$ roughly scales as $1/\left(\frac{v_i}{p_i} - \frac{v_j}{p_j}\right)^2$. Notice, though, that our upper bound of $2 \sqrt{R \log(R)}$ for the increase in compensation is extremely loose as the expected number of customers who show up among $\frac{R}{p_i}$ type~$i$ customers is in fact the same as that among $\frac{R}{p_j}$ type $j$ customers.} 

\emph{Proof of Lemma \ref{lem:chernoff_difference}. }
We prove the lemma by constructing a contradiction. Take $R_{ij}$ as constructed in Lemma \ref{lem:exchange}. For an optimal solution $\vec{x}^\star$, if neither (i) nor (ii) is true, i.e., 
$$x_i^\star \leq N\kzedit{[}t,T\kzedit{]}-R_{ij}/p_i \text{ and } x_j^\star \geq R_{ij}/p_j,$$ 
then, $\vec{x}^\star - \frac{R_{ij}}{p_j} \cdot e_j + \frac{R_{ij}}{p_i} \cdot e_i$ is a feasible solution, and 
from Lemma \ref{lem:exchange}, we know its objective is greater than that of $\vec{x}^\star$, which contradicts the assumption that $\vec{x}^\star$ is an optimal solution to \eqref{eq:offline}.~\halmos

\kzedit{Now we formally compare a globally optimal solution $x^\star$ with an optimal index solution $x^\prime$. Based on Claim \ref{cl:global_optimality}, we know that there exists  a globally optimal solution $x^\star$ that is locally optimal for every type $j$. For the optimal index solution $x^\prime,$ we only assume that it is locally optimal at its threshold index and it need not be locally optimal at any other index. Since the index policy is allowed to accept any number of requests at its threshold index, we are guaranteed that there exists an index solution that is locally optimal at the threshold index.}

\begin{lemma}\label{lem:global_sensitivity}
Consider an optimal solution $\vec{x}^\star$ that is locally optimal for every type $j$, and an optimal index solution $\vec{x}^\prime$ that is locally optimal at its threshold index $\widetilde{j}$, both for the optimization problem \eqref{eq:offline} with arrivals $\vec{N}$.\footnote{We assume without loss of generality that $N_j>0\forall j$ as types $j$ with $N_j=0$ cannot be accepted by either solution.} With $\delta$ as constructed in Lemma $\ref{lem:local_sensitivity}$, we have
\begin{enumerate}[label=(\roman*)]
    \item $\sum_j \left(x_j^\prime - x_j^\star\right)^+ \leq k \delta \left[1 +\sum_j \left(x_j^\star - x_j^\prime\right)^+\right]$,
    \item $\sum_j \left(x_j^\star - x_j^\prime\right)^+ \leq k \delta \sum_j \left(x_j^\prime - x_j^\star\right)^+$. 
\end{enumerate}
\end{lemma}

\emph{Proof of Lemma \ref{lem:global_sensitivity} (i). } 
{
Recall that we have $x_j'=N_j$  for every $j<\widetilde{j}$, $x_j'=0$ for every $j>\widetilde{j}$, and~$x_{\widetilde{j}}'>0$.}
We start by showing that $\vec{x}^\prime$ is locally optimal at every $j > \widetilde{j}$, and $\vec{x}^\prime - e_{\widetilde{j}}$ is locally optimal at every~$j < \widetilde{j}$. 
From Lemma~\ref{lem:local_optimality}~(i), since $\vec{x}'$ is locally optimal at $\widetilde{j}$,  we find that either~$x_{\widetilde{j}}^\prime = N_{\widetilde{j}},$ in which case $x_{\widetilde{j}+1}^\prime = 0$, and
$$
\mathbb{P}\left[\sum_{n} X_{n,x_n^\prime} \geq B \right] > q_{\widetilde{j}+1} \geq q_j, \forall j > \widetilde{j},\;
\text{ or }\;
\mathbb{P}\left[\sum_{n} X_{n,x_n^\prime} \geq B \right] > q_{\widetilde{j}} \geq q_j, \forall j > \widetilde{j}.$$
Thus, for every $j > \widetilde{j}$ both conditions of Lemma \ref{lem:local_optimality} are fulfilled with 
$$x_j^\prime = 0 \text{ and } \mathbb{P}\left[\sum_{n} X_{n,x_n^\prime} \geq B \right] > q_j,$$
which implies that $\vec{x}^\prime$ is locally optimal at all such $j$.

\noindent Moreover, applying Lemma \ref{lem:local_optimality} (ii) to $\widetilde{j}$, where $x_{\widetilde{j}}^\prime > 0$, we find for all $j<\widetilde{j}$
$$\mathbb{P}\left[\sum_{n \neq \widetilde{j}} X_{n,x_n^\prime} + X_{\widetilde{j},x_{\widetilde{j}}-1} \geq B \right] \leq q_{\widetilde{j}} \leq q_j.$$ 
Thus, for every $j < \widetilde{j}$, we have $x_j^\prime = N_j \text{ and } \mathbb{P}\left[\sum_{n \neq \widetilde{j}} X_{n,x_n^\prime} + X_{\widetilde{j},x_{\widetilde{j}}-1} \geq B \right] \leq q_j,$ so $\vec{x}^\prime - e_{\widetilde{j}}$ fulfills both conditions of Lemma \ref{lem:local_optimality} at every $j < \widetilde{j}$, and is thus locally optimal at each such $j$.  

Next, define $\hat{j}$ to be either $\widetilde{j}$ or $\widetilde{j}-1$ depending on whether $x_{\widetilde{j}}^\prime \geq x_{\widetilde{j}}^\star$ or not, to ensure 
$x_j^\prime \geq x_j^\star$ for all $j \leq \hat{j}$ and $x_j^\prime \leq x_j^\star$ for all $j > \hat{j}$ --- these follow for $j<\widetilde{j}$ from $x_j'=N_j$, for $j>\widetilde{j}$ from $x_j'=0$, and for $j=\widetilde{j}$ from the definition of $\hat{j}$. Then, to prove Lemma \ref{lem:global_sensitivity} (i), it suffices to show for $j\leq\hat{j}$
\begin{equation}\label{eq:delta_bound}x_j^\prime - x_j^\star \leq \delta (1+\sum_{n = \hat{j}+1}^{k} (x_n^\star - x_n^\prime)).
\end{equation}

In particular, when this inequality holds true for $j = \{1,...,\hat{j}\}$, we can sum over $j \leq \hat{j}$ to obtain
$$\sum_{j} \left(x_j^\prime - x_j^\star\right)^+=\sum_{n=1}^{\hat{j}} \left(x_n^\prime - x_n^\star\right)^+ \leq k \delta \left[1 +\sum_{n = \hat{j}+1}^{k} (x_n^\star - x_n^\prime)\right] = k \delta \left[1 +\sum_j \left(x_j^\star - x_j^\prime\right)^+\right].$$

We prove Inequality \eqref{eq:delta_bound} by contradiction. Specifically, suppose $$x_i^\prime - x_i^\star > \delta \left(1+\sum_{n = \hat{j}+1}^{k} (x_n^\star - x_n^\prime)\right)$$ for some $i \leq \hat{j}$. Then  $x_i^\star < N_i$ and we prove that
\begin{equation}\label{eq:loc_opt_1}
\mathbb{P}\left[\sum_{n} X_{n,x_n^\star} \geq B \right]  \leq q_i,
\end{equation}
which contradicts $\vec{x}^\star$ fulfilling the local optimality conditions from  Lemma \ref{lem:local_optimality} at $i$. 

We consider two cases
\begin{enumerate}[label=(\alph*)]
    \item $i < \hat{j} \leq \widetilde{j}$ or $i \leq \hat{j} < \widetilde{j}$
    \item $i = \hat{j} = \widetilde{j}$
\end{enumerate}

For Case (a), let $m = \sum_{n = \hat{j}+1}^{k} (x_n^\star - x_n^\prime)$. We construct a sequence of solutions $$\vec{x}^1,\vec{x}^2,\vec{x}^3,...,\vec{x}^{m+1},$$ where~$\vec{x}^{m+1}$ has the following properties:
\begin{itemize}
    \item $\vec{x}^{m+1}$ is locally optimal at $i$,
    \item $x_i^{m+1} > x_i^\star$,
    \item $x_j^{m+1} = x_j^\prime=N_j\geq x_j^\star$ for all $j \leq \hat{j}$ such that $j \neq i$,
     and
    \item $x_j^{m+1} = x_j^\star$ for all $j > \hat{j}$.
\end{itemize}
With $x_j^{m+1}\geq x_j^\star\forall j$ we have the first inequality, and with $\vec{x}^{m+1}$ locally optimal at $i$ and $x_i^{m+1}>0$ we have (Lemma \ref{lem:local_optimality} (ii)) the second inequality in
$$\mathbb{P}\left[\sum_{n} X_{n,x_n^\star} \geq B \right] \leq \mathbb{P}\left[\sum_{n \neq i} X_{n,x_n^{m+1}} + X_{i,x_i^{m+1}-1} \geq B \right] \leq q_i.$$ 
Thus, such $\vec{x}^{m+1}$ implies Inequality \eqref{eq:loc_opt_1}, i.e., a contradiction to $\vec{x}^\star$ being locally optimal at $i$. We now show how to derive $\vec{x}^{m+1}$ by inductively constructing $\vec{x}^h$ from $\vec{x}^{h-1}$ where for every $h \in \{1,2,...,m+1\}$, $\vec{x}^h$ is locally optimal at $i$, $x_i^{h} > x_i^\star + \delta (m+1-h)$, and $\vec{x}^{m+1}$ fulfills the last two properties above.

Recall that $\vec{x}^\prime - e_{\widetilde{j}}$ is locally optimal at $i<\widetilde{j}$, so by Lemma ~\ref{lem:local_sensitivity} there exists $0 \leq l_1 \leq \delta$ such that
$$\vec{x}^1 = \vec{x}^\prime - e_{\widetilde{j}} + e_{\widetilde{j}}- l_1 e_i = \vec{x}^\prime - l_1 e_i$$
is locally optimal at $i$. Further, since $l_1 \leq \delta$, we know $x_i^1 - x_i^\star > \delta m$. Now, for each $h \in \{1,2,...,m\}$, we repeat this procedure of adding one customer of type $j > \hat{j}$: given a solution $\vec{x}^{h}$ that is locally optimal at~$i$, Lemma \ref{lem:local_sensitivity} allows us to find the $l_{h+1}$ such that $\vec{x}^{h+1} = \vec{x}^{h} + e_j - l_{h+1} e_i$ is locally optimal at~$i$. Moreover, Lemma \ref{lem:local_sensitivity} guarantees that $l_{h+1} \leq \delta$, i.e., $x_i^{h} - x_i^{h+1} \leq \delta$. Thus, the first two properties are maintained throughout. Since we have $x^1_j\geq x_j^\star$ for $j\leq\hat{j}$ when $j\neq i$, and we never remove a customer from such $j$, the third property holds, and since we add customers until $x^h_j=x^\star_j$ for $j>\hat{j}$ the fourth property holds. This completes the proof of Case (a).  








For Case (b)  the sequence of solutions starts at $\vec{x}^\prime$, not  $\vec{x}^1$, because $\vec{x}^\prime$ is already locally optimal at~$i$ when $i = \hat{j} = \widetilde{j}$. Then, we derive the same contradiction to $x_i^\prime - x_i^\star > \delta \sum_{n = \hat{j}+1}^{k} (x_n^\star - x_n^\prime) = \delta m$ by constructing $\vec{x}^{m+1}$ with the properties as above. 

The proof of Lemma \ref{lem:global_sensitivity} (ii) is based on a similar argument and is included in Appendix~\ref{app:global_sensitivity}. \halmos





\emph{Proof of Theorem \ref{thm:index_loss}. }
Consider $p_{min} = \min_{j} p_j < 1$, $R_{max} = \max_{i,j:i<j} R_{ij} \geq 1$ for $R_{ij}$ as constructed in Lemma \ref{lem:chernoff_difference}, and  $\delta \geq 1$ as constructed in Lemma \ref{lem:local_sensitivity}. To prove the theorem we show
$$\widehat{\opt}_{\vec{A}}-\calH_{\vec{A}}\kzedit{[}1\kzedit{]} \leq M_1 := \delta k (k+2) \frac{R_{max}}{p_{min}}, \forall \vec{A},$$
and observe that all terms on the right are independent of $B$ and $T$.

Suppose the clairvoyant general solution, locally optimal at each $j$, given $B$ and $\vec{A}$, is $\vec{x}^\star$ and a clairvoyant index solution, that is locally optimal at its threshold index $\widetilde{j}$, is $\vec{x}^\prime$. Then 

$$\widehat{\opt}_{\vec{A}}-\calH_{\vec{A}}\kzedit{[}1\kzedit{]} = \left(\sum_j v_j x_j^\star - V_B(\vec{x}^\star)\right) - \left(\sum_j v_j x_j^\prime - V_B(\vec{x}^\prime)\right) $$
$$= \left(\sum_j v_j (x_j^\star - x_j^\prime)\right) + \left(V_B(\vec{x}^\prime) - V_B(\vec{x}^\star)\right)\leq  \left(\sum_j(x_j^\star - x_j^\prime)^+\right) + \left(\sum_j(x_j^\prime-x_j^\star)^+\right). $$ 
Notice that the inequality is a seemingly loose bound: it ignores greater revenue for $\vec{x}'$ from types~$j$ with~$x_j'>x_j^\star$, it compensates for every such customer, regardless of whether or not such compensation would be paid, it rounds up $v_j$ to 1 for $j$ with $x_j^\star>x_j'$, and it ignores the compensation for these types. Nonetheless, it is sufficient for our purposes.  We discuss two main cases based on the threshold index $\widetilde{j}$ of $\vec{x}'$, where we recall that $x_j^\prime = N_j \geq x_j^\star$ for $j < \widetilde{j}$, and $x_j^\prime = 0 \leq x_j^\star$ for $j > \widetilde{j}$.


\begin{enumerate}[label=(\roman*)]
    \item $x_j^\star < \frac{R_{max}}{p_{min}}, \forall j \geq \widetilde{j}$. Then we consider
    \begin{enumerate}[label=(\alph*)]
        \item $x_{\widetilde{j}}^\prime \geq x_{\widetilde{j}}^\star$
        \item $x_{\widetilde{j}}^\prime < x_{\widetilde{j}}^\star$
    \end{enumerate}
    \item Find the largest index $j \geq \widetilde{j}$ such that $x_j^\star \geq \frac{R_{max}}{p_{min}}$, denoted by $\hat{j}$. Then we consider
    \begin{enumerate}[label=(\alph*)]
        \item $x_{\hat{j}}^\prime \geq x_{\hat{j}}^\star$
        \item $x_{\hat{j}}^\prime < x_{\hat{j}}^\star$
    \end{enumerate}
\end{enumerate}

We begin with Case (i.a). Since $x_j^\prime = N_j \geq x_j^\star$ for any $j < \widetilde{j}$, 

$$\sum_j (x_j^\star - x_j^\prime) \leq \sum_{j = \widetilde{j}}^k x_j^\star \leq (k-\widetilde{j}+1) \frac{R_{max}}{p_{min}}.$$

Then we find that

\begin{align*}
\sum_{j} (x_j^\prime - x_j^\star)^+\leq & k \delta \left[1 +\sum_{j = \widetilde{j}+1}^{k} (x_j^\star - x_j^\prime)\right]  \tag{Lemma \ref{lem:global_sensitivity} (i)} \\
\leq & k \delta \left[1 + (k-\widetilde{j})\frac{R_{max}}{p_{min}}\right]  \tag{Assumption of Case (i.a)}\\
\leq & \delta k^2 \frac{R_{max}}{p_{min}}\notag.
\end{align*}

 $$\text{Thus, in Case (i.a):  }\;\widehat{\opt}_{\vec{A}}-\calH_{\vec{A}}\kzedit{[}1\kzedit{]} \leq (k-\widetilde{j}+1) \frac{R_{max}}{p_{min}} + \delta k (k+1) \frac{R_{max}}{p_{min}} \leq \delta k (k+2) \frac{R_{max}}{p_{min}}.$$

The arguments for Case (i.b) are similar to those for Case (i.a) and give
$$\sum_{j} (x_j^\prime - x_j^\star)^+ \leq k \delta \left[1 +\sum_{j = \widetilde{j}}^{k} (x_j^\star - x_j^\prime)\right]\leq k \delta \left[1+(k-\widetilde{j}+1) \frac{R_{max}}{p_{min}}\right] \leq \delta k (k+1) \frac{R_{max}}{p_{min}}.$$

 $$\text{Thus, in Case (i):  }\;\widehat{\opt}_{\vec{A}}-\calH_{\vec{A}}\kzedit{[}1\kzedit{]} \leq (k-\widetilde{j}+1) \frac{R_{max}}{p_{min}} + \delta k (k+1) \frac{R_{max}}{p_{min}} \leq \delta k (k+2) \frac{R_{max}}{p_{min}}.$$

For Case (ii.a) we argue as follows. We know that $x_j^\prime = 0 \leq x_j^\star, \forall j > \widetilde{j}$, so for $x_{\hat{j}}^\prime \geq x_{\hat{j}}^\star \geq \frac{R_{max}}{p_{min}}$ to hold we must have $\hat{j}=\widetilde{j}$. Moreover, $x_j^\prime = N_j \geq x_j^\star, \forall j < \hat{j} = \widetilde{j}$. Thus, 
$$
\sum_j(x_j^\star - x_j^\prime)^+ \leq \sum_{j = \hat{j}+1}^k x_j^\star \leq (k-\hat{j}) \frac{R_{max}}{p_{min}}
$$
$$
\text{and } \sum_j(x_j'-x_j^\star)^+ = \sum_{j = 1}^{\hat{j}} (x_j^\prime - x_j^\star) \leq k \delta \left[1 +\sum_{j = \hat{j}+1}^{k} (x_j^\star - x_j^\prime)\right] \leq k \delta \left[1 + (k-\hat{j}) \frac{R_{max}}{p_{min}} \right] \leq \delta k^2 \frac{R_{max}}{p_{min}}
$$
$$\text{imply } \widehat{\opt}_{\vec{A}}-\calH_{\vec{A}}\kzedit{[}1\kzedit{]} \leq (k-\hat{j}) \frac{R_{max}}{p_{min}} + \delta k^2 \frac{R_{max}}{p_{min}} \leq \delta k (k+1) \frac{R_{max}}{p_{min}}.$$

In Case (ii.b), we know from Lemma \ref{lem:chernoff_difference} that, if $x_{\hat{j}}^\star \geq \frac{R_{max}}{p_{min}}$, then $x_j^\star > N_j - \frac{R_{max}}{p_{min}}\; \forall j < \hat{j}$. Thus, $x_j^\prime - x_j^\star \leq \frac{R_{max}}{p_{min}}, \forall j < \hat{j}$, implying
$$\sum_j (x_j'-x_j^\star)^+ \leq \sum_{j = 1}^{\hat{j}-1} (x_j^\prime - x_j^\star) \leq \sum_{j = 1}^{\hat{j}-1}\frac{R_{max}}{p_{min}} \leq (\hat{j}-1) \frac{R_{max}}{p_{min}}.$$
Moreover, with Lemma \ref{lem:global_sensitivity} (ii) implying the second, and Lemma \ref{lem:chernoff_difference} (i) implying the third inequality in 
$$
\sum_j (x_j^\star - x_j^\prime)^+ \leq \sum_{j = \hat{j}}^k (x_j^\star - x_j^\prime)
\leq k \delta \sum_{j = 1}^{\hat{j}-1} (x_j^\prime - x_j^\star)
\leq \delta k^2 \frac{R_{max}}{p_{min}},
$$
$$\text{we conclude }\;\widehat{\opt}_{\vec{A}}-\calH_{\vec{A}}\kzedit{[}1\kzedit{]} \leq \delta k^2 \frac{R_{max}}{p_{min}} + (\hat{j}-1) \frac{R_{max}}{p_{min}} \leq \delta k (k+1) \frac{R_{max}}{p_{min}}.$$
\noindent Thus, in any case, the loss is bounded by $\delta k (k+2) \frac{R_{max}}{p_{min}}$,  independent of $B$, $T$ and ~$\vec{A}$. \halmos

\subsection{Proof of Lemma \ref{lem:small_loss_probability}} \label{sec:lemma3}
The proof of the lemma follows ideas from \cite{vera2019bayesian,freund2019uniform}. Suppose the arrival in period~$t-1$ is of type~$i$, and the online index policy accepts the arrival --- the proof is symmetric, but the notation is more cumbersome, when the arrival is rejected. 
We denote by~$\vec{x}^\star$ the clairvoyant index solution in period~$t-1$ and by $\vec{x}'$ the solution found by Algorithm \ref{alg:base} in period~$t-1$, i.e.,~$\vec{x}^\star$ and~$\vec{x}'$ are both optimal index solutions to \eqref{eq:offline} with past actions $\vec{x}[t-2]$, but~$\vec{x}^\star$ is based on the real arrivals, from~$\vec{A}$, in periods~$t,t+1,\ldots,T$ whereas $\vec{x}'$ is based on the sampled arrivals from $\vec{A'}$. 

Observe first that we must have $x_i'\geq \tildeN_i\kzedit{[}t-1\kzedit{]}/2$ when we are in the case that the online index policy accepts $i$ in period $t-1$ (see Algorithm \ref{alg:base}). Now, if $x_i^\star>0$, then the clairvoyant index policy (in period $t-1$) still accepts at least one arrival of type $i$ in periods $t-1,t,\ldots,T$. But that implies that the clairvoyant index policy objective in period $t-1$ is still achievable in period $t$. Thus, the online index policy does not incur any loss by accepting an arrival of type $i$ in period $t-1$, and we must have $\calH_{\vec{A}}\kzedit{[}t-1\kzedit{]}-\calH_{\vec{A}}\kzedit{[}t\kzedit{]}=0$. We derive that in order to incur a loss when accepting $i$ in period~$t-1$ it must be the case that $x_i^\star=0$ and $x_i'\geq \tildeN_j\kzedit{[}t-1\kzedit{]}/2$, i.e.,
\begin{equation}\label{eq:bound_objective}
\mathbb{P}\left[\calH_{\vec{A}}\kzedit{[}t-1\kzedit{]}-\calH_{\vec{A}}\kzedit{[}t\kzedit{]}>0\right]\leq\mathbb{P}\left[x_i'-x_i^\star> \tildeN_i\kzedit{[}t-1\kzedit{]}/2 -1\right].
\end{equation}

At the same time, since $\vec{x}'$ and $\vec{x}^\star$ are index policies, (i) $x_i^\star=0$ and $N_i\kzedit{[}t-1,T\kzedit{]} \geq 1$ imply that~$x_j^\star=0\frall j\geq i$, and (ii) $x_i'\geq \tildeN_i\kzedit{[}t-1\kzedit{]}/2 > 0$ implies ~$x_j'=\tildeN_j\kzedit{[}t-1\kzedit{]}\frall j<i$. We make the following claim, which we prove in Appendix \ref{app:cl_final_lemma}.
\begin{claim}\label{cl:final_lemma}
With $\vec{x}',\vec{x}^\star$ as described we have $$x_i'-x_i^\star\leq \delta\left[\sum_{j=1}^{i-1} \Big(x_j^\star-x_j')^+\right].$$
\end{claim}
The proof of Claim \ref{cl:final_lemma} is straightforward from Lemma \ref{lem:local_sensitivity} and included in Appendix \ref{app:cl_final_lemma}. Combining the claim with the reasoning above we obtain
\begin{align}
x_i'-x_i^\star \leq \delta\left[\sum_{j=1}^{i-1} \Big(x_j^\star-x_j')^+\right] \tag{Claim \ref{cl:final_lemma}} \\
= \delta \left[\sum_{j=1}^{i-1} \Big(x_j^\star-\tildeN_j\kzedit{[}t-1\kzedit{]}\Big)^+\right] \tag{ii} \\
\leq \delta \left[\sum_{j=1}^{i-1} \Big(N_j\kzedit{[}t-1, T\kzedit{]}-\tildeN_j\kzedit{[}t-1\kzedit{]}\Big)^+\right]. \tag{$x_j^\star\leq N_j\kzedit{[}t-1\kzedit{]}$}
\end{align}

Together with \eqref{eq:bound_objective} this implies that
$$
\mathbb{P}\Big[\calH_{\vec{A}}\kzedit{[}t-1\kzedit{]}-\calH_{\vec{A}}\kzedit{[}t\kzedit{]}>0\Big]
\leq \mathbb{P}\left[\delta \left(\sum_{j=1}^{i-1} \Big(N_j\kzedit{[}t-1,T\kzedit{]}-\tildeN_j\kzedit{[}t-1\kzedit{]}\Big)^+\right)> \tildeN_i\kzedit{[}t-1\kzedit{]}/2 -1\right].
$$
Let $\bar{t}=T-t$. Applying Chernoff bounds to $\vec{N}\kzedit{[}t-1\kzedit{]}$ and $\tildeN\kzedit{[}t-1\kzedit{]}$ with $\epsilon = \frac{\log(\bar{t})}{\lambda_j\sqrt{\bar{t}}}$, we know that 
\begin{eqnarray*}
\text{with $E_1^j=\left\{N_j\kzedit{[}t-1\kzedit{]}-\lambda_j\bar{t} \geq \log(\bar{t})\sqrt{\bar{t}}\right\}$ we have }\; \mathbb{P}\left[E_1^j\right] \leq e^{-\frac{\epsilon^2\lambda_j\bar{t}}{2+\epsilon}} \in o\left(\frac{1}{\bar{t}^2}\right), \forall j \label{eq:chernoff_N}\\
\text{With $E_2^j=\left\{\tildeN_j\kzedit{[}t-1\kzedit{]}-\lambda_j\bar{t}\leq -\log(\bar{t})\sqrt{\bar{t}}\right\}$ we have }\; \mathbb{P}\left[E_2^j\right] \leq e^{-\frac{\epsilon^2\lambda_j\bar{t}}{\epsilon}} \in o\left(\frac{1}{\bar{t}^2}\right),\forall j \label{eq:chernoff_N_tilde}
\end{eqnarray*}

Thus, there exists some constant $t_1 \in \mathbb{Z}^+$ such that $\mathbb{P}\left[E_1^j\right],\mathbb{P}\left[E_2^j\right] < \frac{1}{\bar{t}^2}$ for any $\bar{t} \geq t_1$. A union bound over $E_1^j,E_2^j$ for all $j$ implies that, for $E = \cap_j \left(E_1^j \cup E_2^j\right)^c$, we have 
$$
\mathbb{P}(E) > 1-\frac{2k}{\bar{t}^2}, \forall \bar{t} \geq t_1.
$$
Next, define $t_2\in \mathbb{Z}^+$ as the smallest value such that for $\bar{t}\geq t_2$
$$
\delta k 2\sqrt{\bar{t}}\log(\bar{t})\leq \lambda_i\bar{t}/2 -1-\frac{1}{2}\sqrt{\bar{t}}.
$$
Then, conditioning on event $E$ we find for $\bar{t} \geq t_2$ that
$$
\delta \left[\sum_{j=1}^{i-1} \Big(N_j\kzedit{[}t-1,T\kzedit{]}-\tildeN_j\kzedit{[}t-1\kzedit{]}\Big)^+\right]\leq \delta k 2\sqrt{\bar{t}}\log(\bar{t})\leq \lambda_i\bar{t}/2 -1-\frac{1}{2}\sqrt{\bar{t}}\log(\bar{t})\leq \tildeN_i\kzedit{[}t-1\kzedit{]}/2 -1,
$$
where the first inequality holds conditioned on $E$, and the second with $\bar{t}\geq t_2$. Let $t_3=\max\{t_1,t_2\}$; then for $T-t \geq t_3$ we have $\mathbb{P}\left[\calH_{\vec{A}}\kzedit{[}t-1\kzedit{]}-\calH_{\vec{A}}\kzedit{[}t\kzedit{]}>0\right]<\frac{2k}{(T-t)^2}$. We conclude the proof with
\begin{eqnarray*}
\sum_{t=2}^T\mathbb{P}\left[\calH_{\vec{A}}\kzedit{[}t-1\kzedit{]}-\calH_{\vec{A}}\kzedit{[}t\kzedit{]}>0\right]\\
\leq t_3+\sum_{t=2}^{T-t_3}\mathbb{P}\left[\calH_{\vec{A}}\kzedit{[}t-1\kzedit{]}-\calH_{\vec{A}}\kzedit{[}t\kzedit{]}>0\right]\\
\leq t_3+\sum_{t=2}^{T-t_3} \frac{2k}{(T-t)^2}
\leq t_3+\sum_{t=1}^\infty \frac{2k}{t^2}= t_3+2k\frac{\pi^2}{6}
\end{eqnarray*}
giving a constant bound $M_2:=t_3+2k\frac{\pi^2}{6}$ as required.\footnote{Notice that different $t_1,t_2$ may be needed to capture, through analogous arguments, the case where the online index policy rejects, rather than accepts, an arrival in period $t$.}\hfill\halmos






\section{Numerical results}\label{sec:numerical}


Our numerical results consist of three parts. In the first part, we compare the optimal clairvoyant general solution and the optimal clairvoyant index solution for a few instances, i.e., we focus on the accepted customers rather than the objective. This serves to illustrate both that (i) index solutions are not optimal in general and (ii) as described in Lemma \ref{lem:chernoff_difference}, asymptotically the clairvoyant general and the clairvoyant index policies look ``similar''. Thereafter, in the second part, we focus on the relative performance gaps bounded in Theorem~\ref{thm:index_loss} and Lemma~\ref{lem:small_loss_probability}. Finally, in the third part we show how, for fixed $B$ and $T$, the loss changes when we vary the parameters $v_i$ and $p_i$. All figures from this section are included at the end of this paper.


\subsubsection*{Optimal solutions}\label{sec:opt_sol}


We first observe that the clairvoyant optimal solution is not guaranteed to be an index solution in all settings, although it ``switches'' to an index policy as $B$ scales up. This switching behavior can be intuitively explained by Lemma \ref{lem:chernoff_difference}, which shows that, when $i < j$, accepting $R_{ij}/p_i$ type $i$ customers yields better objective than accepting $R_{ij}/p_j$ types $j$ customers does for some large constant $R_{ij}.$ Thus, when $B$ and $T$ are large, if $R_{ij}/p_i$ type $i$ customers are present, the clairvoyant general policy does not accept more than $R_{ij}/p_j$ types $j$ customers.

We capture this switching behavior in the example below, where we consider a single sample path based on the following parameters: $k = 3, \lambda_1 = 0.3,\lambda_2 = 0.2,\lambda_3 = 0.5.$ Moreover, the values and no-show probabilities are $v_1 = 0.044, v_2 = 0.1, v_3 = 0.06, p_1 = 0.2, p_2 = 0.5, p_3 = 0.3.$ We test for $B \in \{1,...,15\}$ two settings, one in which demand is unconstrained ($T$ is large), so any number of customers of each type can be accepted, and one in which $T=5B$, so the number of customers accepted is constrained by both the capacity and the demand for each type. For each setting, we consider both the clairvoyant general and the clairvoyant index policy.

\begin{figure}[ht]
    \centering
    \includegraphics[scale=0.5]{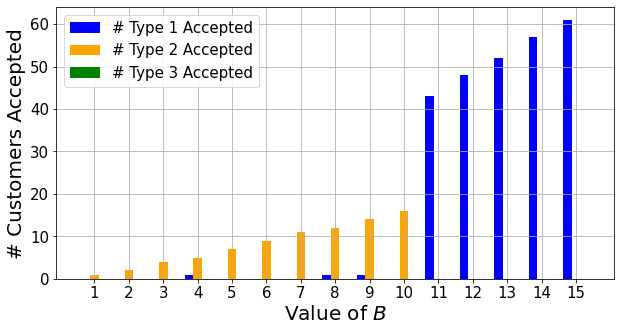}
    \caption{Switching behavior of the clairvoyant general policy under unconstrained demand.}
    \label{fig:1 (a)}
\end{figure}

\begin{figure}
    \centering
    \includegraphics[scale=0.5]{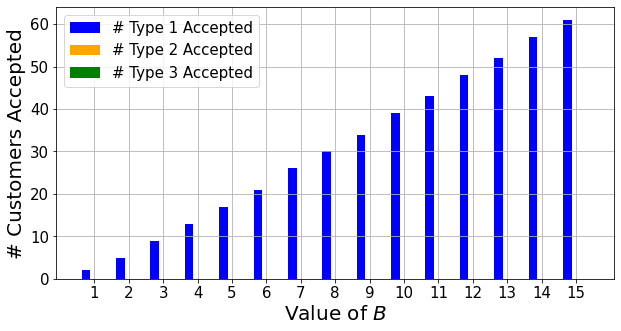}
    \caption{Solution of the  clairvoyant index policy under unconstrained demand.}
    \label{fig:1 (b)}
\end{figure}

In the first setting, we observe a switch from type 2 customers to type 1 customers as $B$ scales up. Figure \ref{fig:1 (a)} shows this switch in the different types of customers accepted by the clairvoyant general policy. On the other hand, Figure \ref{fig:1 (b)} shows that the clairvoyant index policy accepts only type 1 customers, since those have the highest critical ratio, and demand is unconstrained in this example. 

\begin{figure}
    \centering
    \includegraphics[scale=0.5]{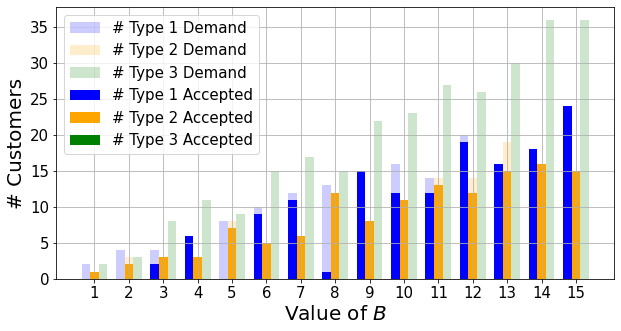}
    \caption{Switching behavior of the clairvoyant general policy with constrained demand.}
    \label{fig:2 (a)}
\end{figure}

\begin{figure}
    \centering
    \includegraphics[scale=0.5]{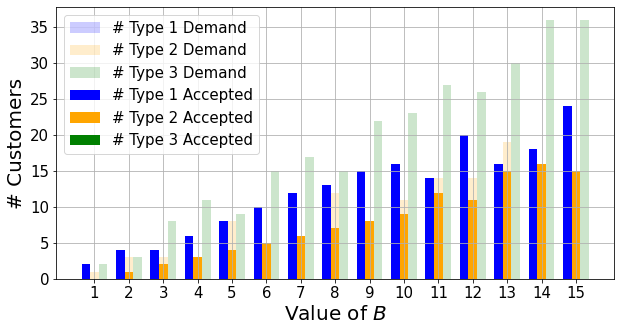}
    \caption{Solution of the clairvoyant index policy with constrained demand.}
    \label{fig:2 (b)}
\end{figure}

The second setting captures a similar behavior when demand is constrained. In Figure \ref{fig:2 (a)} we show both the demand and the accepted number of requests for each type, as we vary $B$. We observe that the clairvoyant general policy in Figure \ref{fig:2 (a)} accepts more customers of type 2 and type 3 customers when $B$ is small, while the clairvoyant index policy in Figure \ref{fig:2 (b)}, by definition, always accepts the types of customers in the order of their indices. As $B$ grows larger, the clairvoyant general policy becomes more similar to an index policy, but sometimes, e.g., with $B=12$, it still does not accept all arrivals of type 1 (in line with Lemma \ref{lem:chernoff_difference}).

\begin{figure}
    \centering
    \includegraphics[scale=0.5]{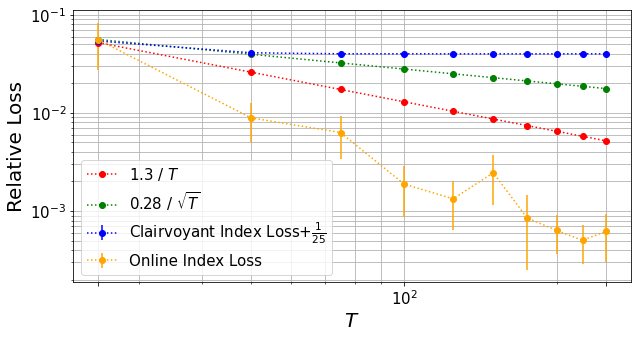}
    \caption{Relative loss of the online and the clairvoyant index policy compared to the clairvoyant general policy as a benchmark (Experiment A).}
    \label{fig:3 (a)}
\end{figure}

\begin{figure}
    \centering
    \includegraphics[scale=0.5]{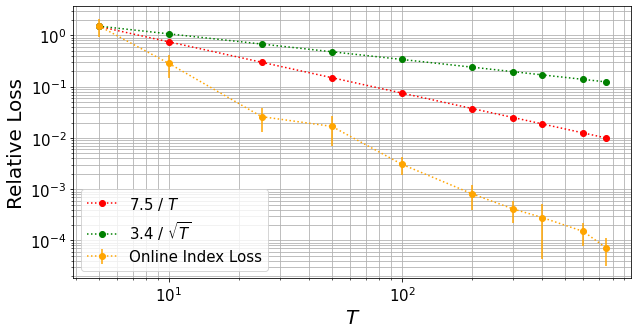}
    \caption{Relative loss of the online index policy compared to the clairvoyant index policy as a benchmark (Experiment A).}
    \label{fig:3 (b)}
\end{figure}

\subsubsection*{Performance loss of clairvoyant index policies}\label{sec:simu_index}

We next turn our attention to the performance loss of the clairvoyant index and the online index policy relative to the clairvoyant general policy. 
As our results show, both  have uniformly bounded loss. We consider two experiments below, corresponding to the case with ``switching" behavior (Experiment~A) and the case without ``switching" behavior (Experiment B). Throughout this part, we show results as a log-log plot, i.e., 1 minus the index policies' objective divided by the objective of the clairvoyant general policy. On these plots, a slope of $-1$, i.e., relative loss scales as $1/B$, indicates a uniform loss guarantee. We include curves proportional to~$\Theta(\frac{1}{\sqrt{B}})$ and $\Theta(\frac{1}{B})$ to better visualize the different scalings.  Further, given that finding the clairvoyant general solution is computationally expensive (even when knowing the arrivals) for a large time horizon, we only evaluate the loss relative to the general clairvoyant for $T\leq 250$.
However, since we find that the clairvoyant index solution seems to have (near-)zero loss anyway, we also include comparisons between the online index policy and the clairvoyant index policy over longer time-horizons.


In Experiment A we have the same parameters as in the first part. Figure \ref{fig:3 (a)} displays the relative loss of the two index policies for $T \in \{25,50,\ldots,250\}$ where $B = T/5$. We observe that the losses of the two policies are uniform as $T$ and $B$ scale up, i.e., the relative loss of the clairvoyant index policy decreases as $1/T$ (or even faster). We remark that the loss of the clairvoyant index solution is so small (indeed, zero after the ``switch", but positive before) that we need to artificially add a constant to it for it to appear on the log-log plot. In Figure \ref{fig:3 (b)} we compare the clairvoyant index policy to the online index policy for a time horizon up to $T=750$.

\begin{figure}
    \centering
    \includegraphics[scale=0.5]{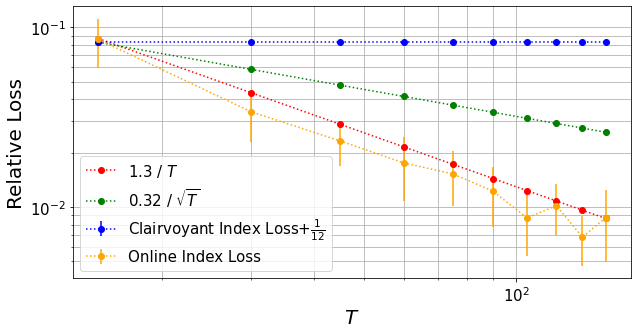}
    \caption{Relative loss of the online and the clairvoyant index policy compared to the clairvoyant general policy as a benchmark (Experiment B)}
    \label{fig:4 (a)}
\end{figure}
\begin{figure}
    \centering
    \includegraphics[scale=0.5]{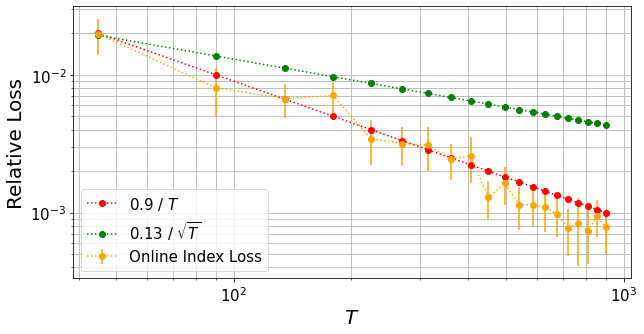}
    \caption{Relative loss of the online index policy compared to the clairvoyant index policy as a benchmark (Experiment B).}
    \label{fig:4 (b)}
\end{figure}

In Experiment B we use a different set of parameters in which we have no switching behavior. We consider $k = 3, \lambda_1 = 0.2,\lambda_2 = 0.3,\lambda_3 = 0.5,$ as well as $v_1 = 0.6, v_2 = 0.4, v_3 = 0.3, p_1 = p_2 = p_3 = 0.8.$ Figure \ref{fig:4 (a)} displays the relative loss. Without the switching behavior the clairvoyant index policy always incurs zero loss, so we add a constant to it for it to appear on the log-log plot. We evaluate both online and clairvoyant index policies for $T \in \{15,10,15,\ldots,150\}$, and $B = T/3$. In Figure \ref{fig:4 (b)} we compare the clairvoyant index policy to the online index policy for a time horizon up to $T=900$.


\subsubsection*{Varying parameters} 
\begin{figure}
    \centering
    \includegraphics[scale=0.5]{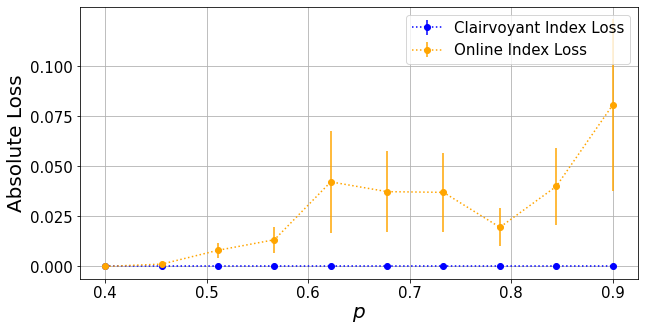}
    \caption{Loss of index policies: $B=10,T=20$, $v_1 = p-0.1, v_2 = p-0.2, v_3 = p-0.3$, where $p_1 = p_2 = p_3 = p$ for $p \in [0.4,0.9]$}
    \label{fig:5}
\end{figure}

\begin{figure}[ht]
    \centering
    \includegraphics[scale=0.5]{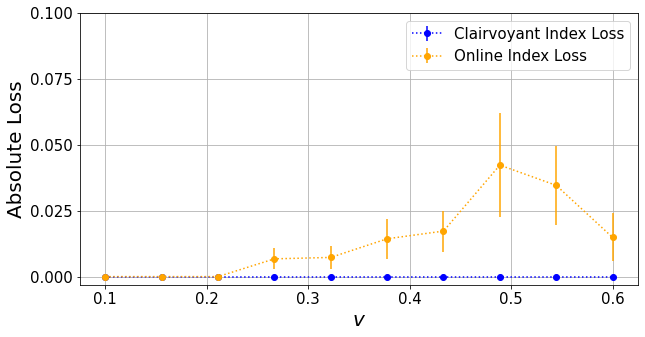}
    \caption{Loss of Index Policies: Revenue per Customer}
    \label{fig:6}
\end{figure}

Finally, we test how the loss bound changes with respect to the other parameters in our model. In Figure \ref{fig:5} and \ref{fig:6} we assume $\lambda_1 = 0.2,\lambda_2 = 0.3,\lambda_3 = 0.5, B = 10, T = 20.$ Moreover, in Figure \ref{fig:5} we let $v_1 = p-0.1, v_2 = p-0.2, v_3 = p-0.3$, where $p_1 = p_2 = p_3 = p$ for $p \in [0.4,0.9]$, and compute the absolute (not relative) loss of the index policies as a function of~$p$. In Figure \ref{fig:6}, we instead assume $v_1 = v_2 = v_3 = v$ and $p_1 = v+0.1,p_2 = v+0.2,p_3 = v+0.3$ for $v \in [0.1,0.6]$, so that we observe how the relative losses of index policies change with respect to the revenue per customer. While the experiments show varying loss for the online clairvoyant index policy as we change these parameters, it is always extremely small relative to the size of the budget and the overall objective for a maximum \emph{relative loss} of about 1\% when compared to the clairvoyant general policy. 

\section{Conclusion}\label{sec:conclusion}
In this work we developed a simple online algorithm for \problem\; with heterogeneous no-show probabilities. In contrast to previous results, that were only optimal on the fluid scale, our algorithm is the first that achieves a uniform loss guarantee. Our key technical innovation is the design of a set of policies, specifically ones based on index solutions, that (i) have bounded loss relative to the clairvoyant over the entire time horizon, and (ii) change in a tractable manner in each period. In doing so, we are able to leverage the novel compensated coupling technique of \cite{vera2019bayesian,freund2019uniform} for a problem with overbooking. 

Our results extend to capture customer requests with (i) product-dependent (heterogeneous) refunds for no-shows (see Appendix \ref{app:refund}), and (ii) different resource requirements (see Appendix \ref{app:multi-unit}), or (iii) arrivals that are not iid. However, there are other ways in which we do not know how to extend our results. The most obvious among those would be to extend our results to a network problem with no-shows, i.e., one with different types of resources. In this case, one may be able to combine a decomposition approach with our technique for individual legs to obtain a uniform loss guarantee (potentially under a non-degeneracy assumption \`a la \cite{jasin2012re}).  Next, it would be interesting to extend our results to capture type-dependent compensation amounts. This complicates the problem, as it is then non-obvious that the clairvoyant general solution asymptotically resembles the clairvoyant index solution --- however, again under some technical assumption, it may be possible to derive such a result.

Finally, the most well-studied extension would be to settings where there are cancellations in addition to no-shows; unfortunately, traditional models of cancellation do not seem to fit under the umbrella of compensated coupling techniques, as the objective is not only based on the action counts, i.e., how often were product requests of each type accepted, but also on the timing of those actions. Thus, this seems to be the least feasible extension of our techniques.

\newpage
\bibliography{template}

\begin{thebibliography}{ABFN13}

\bibitem[ABFN13]{aydin2013single}
Nur{\c{s}}en Ayd{\i}n, {\c{S}}~{\.I}lker Birbil, JBG Frenk, and Nilay Noyan.
\newblock Single-leg airline revenue management with overbooking.
\newblock {\em Transportation Science}, 47(4):560--583, 2013.

\bibitem[AD14]{agrawal2014fast}
Shipra Agrawal and Nikhil~R Devanur.
\newblock Fast algorithms for online stochastic convex programming.
\newblock In {\em Proceedings of the twenty-sixth annual ACM-SIAM symposium on
  Discrete algorithms}, pages 1405--1424. SIAM, 2014.

\bibitem[ADL16]{agrawal2016efficient}
Shipra Agrawal, Nikhil~R Devanur, and Lihong Li.
\newblock An efficient algorithm for contextual bandits with knapsacks, and an
  extension to concave objectives.
\newblock In {\em Conference on Learning Theory}, pages 4--18. PMLR, 2016.

\bibitem[AG19]{arlotto2019uniformly}
Alessandro Arlotto and Itai Gurvich.
\newblock Uniformly bounded regret in the multisecretary problem.
\newblock {\em Stochastic Systems}, 2019.

\bibitem[AHL12]{alaei2012online}
Saeed Alaei, MohammadTaghi Hajiaghayi, and Vahid Liaghat.
\newblock Online prophet-inequality matching with applications to ad
  allocation.
\newblock In {\em Proceedings of the 13th ACM Conference on Electronic
  Commerce}, pages 18--35, 2012.

\bibitem[AHL13]{alaei2013online}
Saeed Alaei, MohammadTaghi Hajiaghayi, and Vahid Liaghat.
\newblock The online stochastic generalized assignment problem.
\newblock In {\em Approximation, Randomization, and Combinatorial Optimization.
  Algorithms and Techniques}, pages 11--25. Springer, 2013.

\bibitem[AX20]{arlotto2018logarithmic}
Alessandro Arlotto and Xinchang Xie.
\newblock Logarithmic regret in the dynamic and stochastic knapsack problem
  with equal rewards.
\newblock {\em Stochastic Systems}, 2020.

\bibitem[Bec58]{beckmann1958decision}
Martin~J Beckmann.
\newblock Decision and team problems in airline reservations.
\newblock {\em Econometrica: Journal of the Econometric Society}, pages
  134--145, 1958.

\bibitem[Bel87]{belobaba1987air}
Peter Belobaba.
\newblock {\em Air travel demand and airline seat inventory management}.
\newblock PhD thesis, Massachusetts Institute of Technology, 1987.

\bibitem[BG96]{bitran1996managing}
Gabriel~R Bitran and Stephen~M Gilbert.
\newblock Managing hotel reservations with uncertain arrivals.
\newblock {\em Operations Research}, 44(1):35--49, 1996.

\bibitem[BM95]{bitran1995application}
Gabriel~R Bitran and Susana~V Mondschein.
\newblock An application of yield management to the hotel industry considering
  multiple day stays.
\newblock {\em Operations research}, 43(3):427--443, 1995.

\bibitem[Bra19]{bray2019does}
Robert~L Bray.
\newblock Does the multisecretary problem always have bounded regret?
\newblock {\em arXiv preprint arXiv:1912.08917}, 2019.

\bibitem[BW18]{bumpensanti2018re}
Pornpawee Bumpensanti and He~Wang.
\newblock A re-solving heuristic with uniformly bounded loss for network
  revenue management.
\newblock {\em arXiv preprint arXiv:1802.06192}, 2018.

\bibitem[CLY21]{chen2021linear}
Guanting Chen, Xiaocheng Li, and Yinyu Ye.
\newblock How linear reward helps in online resource allocation.
\newblock {\em arXiv preprint arXiv:2101.11092}, 2021.

\bibitem[Coo02]{cooper2002asymptotic}
William~L Cooper.
\newblock Asymptotic behavior of an allocation policy for revenue management.
\newblock {\em Operations Research}, 50(4):720--727, 2002.

\bibitem[DKX19]{dai2019network}
Jiangang Dai, Anton~J Kleywegt, and Yongbo Xiao.
\newblock Network revenue management with cancellations and no-shows.
\newblock {\em Production and Operations Management}, 28(2):292--318, 2019.

\bibitem[ET10]{erdelyi2010dynamic}
Alexander Erdelyi and Huseyin Topaloglu.
\newblock A dynamic programming decomposition method for making overbooking
  decisions over an airline network.
\newblock {\em INFORMS Journal on Computing}, 22(3):443--456, 2010.

\bibitem[FB19]{freund2019uniform}
Daniel Freund and Siddhartha Banerjee.
\newblock Uniform loss algorithms for online stochastic decision-making with
  applications to bin packing.
\newblock {\em Available at SSRN 3479189}, 2019.

\bibitem[Fel57]{feller1957introduction}
William Feller.
\newblock An introduction to probability theory and its applications.
\newblock 1957.

\bibitem[GJ97]{geraghty1997revenue}
M~Karen Geraghty and Ernest Johnson.
\newblock Revenue management saves national car rental.
\newblock {\em Interfaces}, 27(1):107--127, 1997.

\bibitem[GR20]{gupta2020interior}
Varun Gupta and Ana Radovanovi{\'c}.
\newblock Interior-point-based online stochastic bin packing.
\newblock {\em Operations Research}, 68(5):1474--1492, 2020.

\bibitem[GT19a]{gallego2019overbooking}
Guillermo Gallego and Huseyin Topaloglu.
\newblock Overbooking.
\newblock In {\em Revenue Management and Pricing Analytics}, pages 83--105.
  Springer, 2019.

\bibitem[GT19b]{gallego2019revenue}
Guillermo Gallego and Huseyin Topaloglu.
\newblock {\em Revenue management and pricing analytics}, volume 209.
\newblock Springer, 2019.

\bibitem[JK12]{jasin2012re}
Stefanus Jasin and Sunil Kumar.
\newblock A re-solving heuristic with bounded revenue loss for network revenue
  management with customer choice.
\newblock {\em Mathematics of Operations Research}, 37(2):313--345, 2012.

\bibitem[KTT12]{kunnumkal2012randomized}
Sumit Kunnumkal, Kalyan Talluri, and Huseyin Topaloglu.
\newblock A randomized linear programming method for network revenue management
  with product-specific no-shows.
\newblock {\em Transportation Science}, 46(1):90--108, 2012.

\bibitem[Lit72]{littlewood1972forecasting}
Kenneth Littlewood.
\newblock Forecasting and control of passenger bookings.
\newblock {\em Airline Group International Federation of Operational Research
  Societies Proceedings, 1972}, 12:95--117, 1972.

\bibitem[LSJ99]{lautenbacher1999underlying}
Conrad~J Lautenbacher and Shaler Stidham~Jr.
\newblock The underlying markov decision process in the single-leg airline
  yield-management problem.
\newblock {\em Transportation Science}, 33(2):136--146, 1999.

\bibitem[MS87]{mangasarian1987lipschitz}
Olvi~L Mangasarian and T-H Shiau.
\newblock Lipschitz continuity of solutions of linear inequalities, programs
  and complementarity problems.
\newblock {\em SIAM Journal on Control and Optimization}, 25(3):583--595, 1987.

\bibitem[MV99]{metters1999yield}
Richard Metters and Vicente Vargas.
\newblock Yield management for the nonprofit sector.
\newblock {\em Journal of Service Research}, 1(3):215--226, 1999.

\bibitem[MVR99]{mcgill1999revenue}
Jeffrey~I McGill and Garrett~J Van~Ryzin.
\newblock Revenue management: Research overview and prospects.
\newblock {\em Transportation science}, 33(2):233--256, 1999.

\bibitem[Ope21]{opentable}
OpenTable.
\newblock {\em OpenTable Support}, 2021.
\newblock Accessed 02-02-2021.

\bibitem[Pin20]{pinelis2020monotonicity}
Iosif Pinelis.
\newblock Monotonicity properties of the poisson approximation to the binomial
  distribution.
\newblock {\em Statistics \& Probability Letters}, 167:108901, 2020.

\bibitem[Rot71]{rothstein1971airline}
Marvin Rothstein.
\newblock An airline overbooking model.
\newblock {\em Transportation Science}, 5(2):180--192, 1971.

\bibitem[RW08]{reiman2008asymptotically}
Martin~I Reiman and Qiong Wang.
\newblock An asymptotically optimal policy for a quantity-based network revenue
  management problem.
\newblock {\em Mathematics of Operations Research}, 33(2):257--282, 2008.

\bibitem[SWZ20]{sun2020near}
Rui Sun, Xinshang Wang, and Zijie Zhou.
\newblock Near-optimal primal-dual algorithms for quantity-based network
  revenue management.
\newblock {\em arXiv preprint arXiv:2011.06327}, 2020.

\bibitem[Tho61]{thompson1961statistical}
HR~Thompson.
\newblock Statistical problems in airline reservation control.
\newblock {\em Journal of the Operational Research Society}, 12(3):167--185,
  1961.

\bibitem[TVR99]{talluri1999randomized}
Kalyan Talluri and Garrett Van~Ryzin.
\newblock A randomized linear programming method for computing network bid
  prices.
\newblock {\em Transportation science}, 33(2):207--216, 1999.

\bibitem[TVR06]{talluri2006theory}
Kalyan~T Talluri and Garrett~J Van~Ryzin.
\newblock {\em The theory and practice of revenue management}, volume~68.
\newblock Springer Science \&amp; Business Media, 2006.

\bibitem[VB20]{vera2019bayesian}
Alberto Vera and Siddhartha Banerjee.
\newblock The bayesian prophet: A low-regret framework for online decision
  making.
\newblock {\em Management Science}, 2020.

\bibitem[VBG19]{vera2019online}
Alberto Vera, Siddhartha Banerjee, and Itai Gurvich.
\newblock Online allocation and pricing: Constant regret via bellman
  inequalities.
\newblock {\em arXiv preprint arXiv:1906.06361}, 2019.

\end{thebibliography}
\newpage



%
%
%
%
%

\appendix
\section{Omitted Proofs}\label{app:proofs}
\subsection{Existence of locally optimal global optimum}\label{app:glocal}
\begin{claim}\label{cl:global_optimality}
There exists a globally optimal solution $\vec{x}^\star$ to \eqref{eq:offline} with (estimated) future arrivals $\tildeN\kzedit{[}t\kzedit{]}$ and past accepted arrivals ~$\vec{x}[t-1]$ that is locally optimal for every type $j$.
\end{claim}

\emph{Proof of Claim \ref{cl:global_optimality}. } Let $\vec{x}^\star$ be a globally optimal solution that maximizes $\sum_ix_i^\star$, i.e., there exists no solution $\vec{x}'$ with objective as high as that of $\vec{x}^\star$ and $\sum_ix_i'>\sum_ix_i^\star$. We prove that $\vec{x}^\star$ is locally optimal for every type $j$. Let $\vec{x} = \vec{x}[t-1] + \vec{x}^\star$. 

Suppose $\vec{x}^\star$ is not locally optimal at some $i \in [k]$. We shall show that this implies
$$x_i^\star < \tildeN_i\kzedit{[}t\kzedit{]} \text{ and } \mathbb{P}\left[\sum_{n} X_{n,x_n} \geq B\right] = q_i.$$
Since $\vec{x}^\star$ is not locally optimal at $i$, it violates either Lemma \ref{lem:local_optimality} (i) or Lemma \ref{lem:local_optimality} (ii). If~$\vec{x}^\star$ violates Lemma \ref{lem:local_optimality} (ii), we must have
$x_i^\star > 0$ and $\mathbb{P}\left[\sum_{n \neq i} X_{n,x_n} + X_{i,x_i -1} \geq B\right] > q_i$. This leads to the contradiction that
$$V_B\left(\vec{x}\right) - V_B\left(\vec{x}-e_i\right) = \mathbb{P}\left[\sum_{n \neq i} X_{n,x_n} + X_{i,x_i-1} < B\right] \cdot 0 + \mathbb{P}\left[\sum_{n \neq i} X_{n,x_n} + X_{i,x_i-1} \geq B\right] \cdot p_i > v_i,$$
which means $\vec{x} - e_i$ yields a better objective than $\vec{x}$ does, and $\vec{x}^\star$ is not a global optimum.

Thus, $\vec{x}^\star$ must violate Lemma \ref{lem:local_optimality} (i), and we have $x_i^\star < \tildeN_i\kzedit{[}t\kzedit{]}$ and $\mathbb{P}\left[\sum_{n} X_{n,x_n} \geq B\right] \leq q_i$. In particular, $x_i^\star < \tildeN_i\kzedit{[}t\kzedit{]}$ and $\mathbb{P}\left[\sum_{n} X_{n,x_n} \geq B\right] < q_i$ would lead to the contradiction that $$V_B\left(\vec{x} + e_i\right) - V_B\left(\vec{x}\right) = \mathbb{P}\left[\sum_{n} X_{n,x_n} < B\right] \cdot 0 + \mathbb{P}\left[\sum_{n} X_{n,x_n} \geq B\right] \cdot p_i < v_i,$$
which means $\vec{x} + e_i$ yields a better objective than $\vec{x}$ does, and $\vec{x}^\star$ is not a global optimum. Thus, we must have $$x_i^\star < \tildeN_i\kzedit{[}t\kzedit{]} \text{ and } \mathbb{P}\left[\sum_{n} X_{n,x_n} \geq B\right] = q_i$$ if the global optimum $\vec{x}^\star$ is not locally optimal at $i$.

From $\mathbb{P}\left[\sum_{n} X_{n,x_n} \geq B\right] = q_i$ we know
$$V_B\left(\vec{x} + e_i\right) - V_B\left(\vec{x}\right) = \mathbb{P}\left[\sum_{n} X_{n,x_n} < B\right] \cdot 0 + \mathbb{P}\left[\sum_{n} X_{n,x_n} \geq B\right] \cdot p_i = v_i,$$
which means $\vec{x} + e_i$ yields the same objective as $\vec{x}$ does. Thus, $\vec{x}^\star+ e_i$ is a global optimum, 
which contradicts the assumption that $\vec{x}^\star$ maximizes $\sum_ix_i^\star$ among globally optimal solutions. \halmos


\subsection{Proof of Claims \ref{cl:berryesseen} and \ref{cl:exchange}}\label{app:cl_berry}

\subsubsection*{Proof of Claim \ref{cl:berryesseen}} Recall that $Y = \sum_{n \neq j} X_{n,x_n}+X_{j,x_j+\delta^\star}$ and $Z = \sum_{n} X_{n,x_n}$ for constant $\delta^\star$. Then, we denote the mean of random variable $Y$ and $Z$ as $\mu_Y$ and $\mu_Z$, and denote the standard deviation of $Y$ and $Z$ as $\sigma_Y$ and $\sigma_Z$. We compute
$$\mu_Z = \mathbb{E}[Z] = \mathbb{E}\left[\sum_{n = 1}^k X_{n,x_n}\right] = \sum_{n = 1}^k x_n p_n,$$
and similarly $\mu_Y = \mathbb{E}[Y] = \delta^\star p_j + \sum_{n = 1}^k x_n p_n$.
Moreover, 
$$\sigma_Z^2 = \text{Var}[Z] = \text{Var}\left[\sum_{n = 1}^k X_{n,x_n}\right] = \sum_{n = 1}^k \text{Var}\left[X_{n,x_n}\right] = \sum_{n = 1}^k x_n p_n (1-p_n),$$
and similarly $\sigma_Y^2 = \text{Var}[Y] = \delta^\star p_j (1-p_j) + \sum_{n = 1}^k x_n p_n (1-p_n)$.

Let $m = \sum_{n = 1}^k x_n$ denote the total number of customers in solution $\vec{x}$ and label the customers from $1$ to $m$. That is, we label each individual customer with $r \in [m]$ rather than labelling each customer type with $n \in [k]$. Since the results in this claim depend only on $m$ (and not on how many customers of each type are in $\vec{x}$), with a slight abuse of notation we let $m\to\infty$ when really referring to $\sum_n x_n\to\infty$. For $r \in \{1,...,m\}$, pick the $r$th customer, and denote its type by $n\in\{1,...,k\}$; we define the random variable $Z_r := X_{n,1} - p_n = X_{n,1} - p_r$.
Then, $\mathbb{E}[Z_r] = 0$ for every $r$ and $$Z = \sum_{r = 1}^m Z_r + \sum_{n = 1}^k x_n p_n = \sum_{r = 1}^m Z_r + \mu_Z.$$
Similarly, for fixed $\delta^\star$, we define for $r \in \{1,...,m+\delta^\star\}$, with the $r$th customer of type $n$, the random variables $Y_r := X_{n,1} - p_n = X_{n,1} - p_r$, so that $\mathbb{E}[Y_r] = r$ for every $r$ and 
$$Y = \sum_{r = 1}^{m+\delta^\star} Y_r + \delta^\star p_j + \sum_{n = 1}^k x_n p_n = \sum_{r = 1}^{m+\delta^\star} Y_r  + \mu_Y.$$

We normalize $Z$ and $Y$ by defining $$Z^\prime = \frac{Z-\mu_Z}{\sigma_Z} = \frac{\sum_{r = 1}^m Z_r}{\sigma_Z}, Y^\prime = \frac{Y-\mu_Y}{\sigma_Y} = \frac{\sum_{r = 1}^{m+\delta^\star} Y_r}{\sigma_Y}.$$
Denote the CDF of standard normal random variable $X_0 \sim \mathcal{N}(0,1)$ by $\Phi(x), \forall x \in \mathbb{R}$. To identify~$\delta^\star,m_0$ that satisfy~\eqref{diff}, we resort to the following version of the Berry-Esseen theorem for independent random variables.

\begin{proposition}[Theorem 2 of Chapter XVI.5 in \cite{feller1957introduction}]\label{pr:berryesseen}
Let $J_r$ be independent random variables such that $\mathbb{E}\left[J_r\right] = 0$, $\mathbb{E}\left[J_r^2\right] := \sigma_r^2$ and $\mathbb{E}\left[\left|J_r^3\right|\right] := \kappa_r$. Put $s_m^2 = \sum_{r = 1}^m \sigma_r^2$, $g_m = \sum_{r = 1}^m \kappa_r$, and denote by $J^\prime$ the normalized sum $(J_1 + ... + J_m)/s_m$. Then for all $x$ and $m$,
\begin{equation}\label{esseen}
\left|F_{J^\prime}(x) - \Phi(x)\right| \leq 6 \frac{g_m}{s_m^3}.
\end{equation}
\end{proposition}

For zero-mean random variables $Z_r$ defined above, we compute
$$g_m = \sum_{r = 1}^m \mathbb{E}\left[\left|Z_r^3\right|\right] = \sum_{r = 1}^m (1-p_r)p_r\left[(1-p_r)^2 + p_r^2\right] \leq \sum_{r = 1}^m (1-p_r)p_r,$$
and
$$s_m^3 = \left(\mathbb{E}\left[Z_r^2\right]\right)^{1.5} = \left[\sum_{r = 1}^m (1-p_r)p_r\right]^{1.5}.$$

Since the $Z_r$ constructed above satisfy the conditions in Proposition \ref{pr:berryesseen}, from \eqref{esseen} we obtain
$$\left|F_{Z^\prime}(x) - \Phi(x)\right| \leq 6 \frac{g_m}{s_m^3} \leq \frac{6}{\sqrt{\sum_{r = 1}^m (1-p_r)p_r}} =  \frac{6}{\sigma_Z}.$$

Similarly, for zero-mean random variables $Y_r$ defined above,
$$\left|F_{Y^\prime}(x) - \Phi(x)\right| \leq \frac{6}{\sqrt{\delta^\star p_j(1-p_j) + \sum_{r = 1}^m (1-p_r)p_r}} = \frac{6}{\sigma_Y}.$$

We assume without loss of generality that $\bar{q}_j-\frac{6}{\sigma_Y} > 0$ and $\bar{q}_j+\frac{6}{\sigma_Y} < 1$, which is always true for sufficiently large $m$. With 
$$
F_{Y^\prime}(x) - \frac{6}{\sigma_Y}\leq \Phi(x),$$
we know, since $\Phi(\cdot)$ is non-decreasing, $$\Phi^{-1}\left(F_{Y^\prime}(x) - \frac{6}{\sigma_Y}\right) \leq \Phi^{-1}(\Phi(x))=x = F_{Y^\prime}^{-1}\left(F_{Y^\prime}(x)\right).$$ Thus, with $x = F_{Y^\prime}^{-1}(\bar{q}_j)$ we have
$$\Phi^{-1}\left(\bar{q}_j-\frac{6}{\sigma_Y}\right)\leq F_{Y^\prime}^{-1}(\bar{q}_j) $$ Similarly, we can show $$F_{Z^\prime}^{-1}(\bar{q}_j) \leq \Phi^{-1}\left(\bar{q}_j+\frac{6}{\sigma_Z}\right).$$

Further, the definition of $Y^\prime$ gives
$$F_{Y^\prime}^{-1}(\bar{q}_j) = \left[F_{Y}^{-1}(\bar{q}_j)-\mu_Y\right] \frac{1}{\sigma_Y},$$
so $F_{Y}^{-1}(\bar{q}_j) = \sigma_Y F_{Y^\prime}^{-1}(\bar{q}_j) + \mu_Y$. Similarly, $F_{Z}^{-1}(\bar{q}_j) = \sigma_Z F_{Z^\prime}^{-1}(\bar{q}_j) + \mu_Z$.

Therefore,
\begin{equation}\label{expansion}
\begin{split}
&F_{Y}^{-1}(\bar{q}_j) - F_{Z}^{-1}(\bar{q}_j) = \sigma_Y F_{Y^\prime}^{-1}(\bar{q}_j) + \mu_Y - \left[\sigma_Z F_{Z^\prime}^{-1}(\bar{q}_j) + \mu_Z\right]\\
\geq & \sigma_Y \cdot \Phi^{-1}\left(\bar{q}_j-\frac{6}{\sigma_Y}\right) + \mu_Y - \left[\sigma_Z \cdot \Phi^{-1}\left(\bar{q}_j+\frac{6}{\sigma_Z}\right) + \mu_Z\right]\\
= & \underbrace{(\mu_Y - \mu_Z)}_{(I)} + \underbrace{\Phi^{-1}\left(\bar{q}_j-\frac{6}{\sigma_Y}\right) (\sigma_Y-\sigma_Z)}_{(II)} -  \underbrace{\left[\Phi^{-1}\left(\bar{q}_j+\frac{6}{\sigma_Z}\right) - \Phi^{-1}\left(\bar{q}_j-\frac{6}{\sigma_Y}\right)\right] \sigma_Z}_{(III)}
\end{split}
\end{equation}

Take $$\delta^\star = \frac{3+12 \Phi^{-1 \prime}(\bar{q}_j)}{p_j}.$$ We aim to show that there exists $m_0 \in \mathbb{Z}^+$ such that \eqref{expansion} $\geq 1$ is guaranteed for $m \geq m_0$. To do so, we show that there exist $m_1,m_2 \in \mathbb{Z}^+$ such that the inequality holds true for any $m \geq m_0 = \max(m_1,m_2).$ Since $(I) = \delta^\star p_j$, we find $m_1$ and $m_2$ that bound (II) and (III), respectively.

{Let $w = \argmin_n p_n(1-p_n)$. 
We know $$\lim_{m \to +\infty} |\sigma_Y-\sigma_Z| = \lim_{m \to +\infty} \left|\sqrt{\delta^\star p_j(1-p_j) + \sum_{r = 1}^m (1-p_r)p_r} - \sqrt{\sum_{r = 1}^m (1-p_r)p_r}\right|$$
$$\leq \lim_{m \to +\infty} \left|\sqrt{\delta^\star p_j(1-p_j) + \sum_{r = 1}^m (1-p_{w})p_{w}} - \sqrt{\sum_{r = 1}^m (1-p_{w})p_{w}}\right| = 0,$$
and thus $\lim_{m \to +\infty} (\sigma_Y-\sigma_Z) = 0$. Similarly, we can show that $$\lim_{m \to +\infty} \frac{\sigma_Y}{\sigma_Z} = 1 \text{ and } \lim_{m \to +\infty} \frac{\sigma_Z}{\sigma_Y^2} = \frac{1}{\sigma_Z} = 0.$$}

To bound (II), we know that $\Phi^{-1}\left(\bar{q}_j-\frac{6}{\sigma_Y}\right)$ is bounded both above and below with respect to~$m$, and $\lim_{m \to +\infty} (\sigma_Y-\sigma_Z) = 0$. Thus, $\lim_{m \to +\infty} (II) = 0$. This implies in particular that there exists~$m_1 \in \mathbb{Z}^+$ such that $(II) \geq -1$ for all $m \geq m_1$.

To bound (III), we apply Taylor expansion. Pick any $\epsilon$ such that $0 < \epsilon < \max(\bar{q}_j,q_j)$.
\begin{enumerate}[label=(\roman*)]
    \item For any $\Delta \in [0,q_j-\epsilon]$,
    \begin{equation} \label{taylor1}
    \Phi^{-1}(\bar{q}_j+\Delta) = \Phi^{-1}(\bar{q}_j) + \Phi^{-1 \prime}(\bar{q}_j) \cdot \Delta + \Phi^{-1 \prime\prime}(a_1) \cdot \Delta^2
    \end{equation}
    for some $a_1 \in [\bar{q}_j,\bar{q}_j+\Delta]$
    \item For any $\Delta \in [\bar{q}_j-\epsilon,0]$,
    \begin{equation} \label{taylor2}
    \Phi^{-1}(\bar{q}_j-\Delta) = \Phi^{-1}(\bar{q}_j) - \Phi^{-1 \prime}(\bar{q}_j) \cdot \Delta + \Phi^{-1 \prime\prime}(a_2) \cdot \Delta^2
    \end{equation}
    for some $a_2 \in [\bar{q}_j-\Delta,\bar{q}_j]$.
\end{enumerate}

Plugging $\Delta = \frac{6}{\sigma_Z}$ into \eqref{taylor1}, we obtain
$$\Phi^{-1}(\bar{q}_j+\Delta) = \Phi^{-1}(\bar{q}_j) + \Phi^{-1 \prime}(\bar{q}_j) \cdot \frac{6}{\sigma_Z} + \Phi^{-1 \prime\prime}(a_1) \frac{36}{\sigma_Z^2},$$ where $a_1 \in [\bar{q}_j,\bar{q}_j+\frac{6}{\sigma_Z}].$ Similarly, for $\Delta = \frac{6}{\sigma_Y}$, we have $$\Phi^{-1}(\bar{q}_j-\Delta) = \Phi^{-1}(\bar{q}_j) - \Phi^{-1 \prime}(\bar{q}_j) \cdot \frac{6}{\sigma_Y} + \Phi^{-1 \prime\prime}(a_2) \frac{36}{\sigma_Y^2},$$
where $a_2 \in [\bar{q}_j-\frac{6}{\sigma_Y},\bar{q}_j].$ Thus, 
\begin{equation*}
    \begin{split}
        & (III) = \left[\Phi^{-1}\left(\bar{q}_j+\frac{6}{\sigma_Z}\right) - \Phi^{-1}(\bar{q}_j) + \Phi^{-1}(\bar{q}_j) - \Phi^{-1}\left(\bar{q}_j-\frac{6}{\sigma_Y}\right)\right] \sigma_Z\\
        = & \left[\Phi^{-1 \prime}(\bar{q}_j) \cdot \frac{6}{\sigma_Z} + \Phi^{-1 \prime\prime}(a_1) \frac{36}{\sigma_Z^2} + \Phi^{-1 \prime}(\bar{q}_j) \cdot \frac{6}{\sigma_Y} - \Phi^{-1 \prime\prime}(a_2) \frac{36}{\sigma_Y^2}\right] \sigma_Z\\
        = & 6 \Phi^{-1 \prime}(\bar{q}_j) \left[1 + \frac{\sigma_Z}{\sigma_Y}\right] + 36 \left[\Phi^{-1 \prime\prime}(a_1) \frac{1}{\sigma_Z} - \Phi^{-1 \prime\prime}(a_2) \frac{\sigma_Z}{\sigma_Y^2}\right]
    \end{split}
\end{equation*}

Since $\Phi^{-1 \prime\prime}(a_1)$ and $\Phi^{-1 \prime\prime}(a_2)$ are both bounded over compact sets, 
$$\lim_{m \to +\infty} 36 \left[\Phi^{-1 \prime\prime}(a_1) \frac{1}{\sigma_Z} - \Phi^{-1 \prime\prime}(a_2) \frac{\sigma_Z}{\sigma_Y^2}\right] = 0.$$ Moreover,
$$\lim_{m \to +\infty} 6 \Phi^{-1 \prime}(\bar{q}_j) \left[1 + \frac{\sigma_Y}{\sigma_Z}\right] = 12 \Phi^{-1 \prime}(\bar{q}_j).$$

Therefore, $\lim_{m \to +\infty} (III) = 12 \Phi^{-1 \prime}(\bar{q}_j),$ which is a constant. Thus, there exists~$m_2 \in \mathbb{Z}^+$ such that $(III) \leq 12 \Phi^{-1 \prime}(\bar{q}_j)+1$ for all $m \geq m_2$. Then, plugging the value of $\delta^\star$ into \eqref{expansion}, we have
$$F_{Y}^{-1}(\bar{q}_j) - F_{Z}^{-1}(\bar{q}_j) \geq \delta^\star p_j - 1 - 12 \Phi^{-1 \prime}(\bar{q}_j) - 1 \geq 1$$
for any $m \geq \max(m_1,m_2)$. \halmos

\subsubsection*{Proof of Claim \ref{cl:exchange}}
We denote the events
$$E_1 = \left\{ \sum_{n \neq i,j} X_{n,x_n}+X_{j,x_j-\delta_{j}-1} < B-1-X_{i,x_i} \right\},$$ $$E_2 = \left\{\sum_{n \neq i,j} X_{n,x_n}+X_{j,x_j-\delta_{j}-1} < B-X_{i,x_i+1} \right\},$$
which allows us to write the claim as $\mathbb{P}\left[E_1^c\right] \geq \mathbb{P}\left[E_2^c\right].$
One can equivalently show $\mathbb{P}\left[E_2\right] \geq \mathbb{P}\left[E_1\right]$, which follows immediately from $X_{i,x_i+1}=X_{i,x_i}+X_{i,1}\leq X_{i,x_i}+1$. \halmos


\subsection{Proof of Lemma \ref{lem:exchange} (ii)}\label{app:chernoff}
In Case (ii), we show that the objective of $\vec{x}^\star+\frac{R_{ij}}{p_i} \cdot e_i$ is always greater than that of $\vec{x}^\star+\frac{R_{ij}}{p_j} \cdot e_j$ for any $R_{ij} \in \mathbb{Z}^+$ such that $\frac{R_{ij}}{p_i}$ and $\frac{R_{ij}}{p_j}$ are both integers.
First, the revenue is the same with $\vec{x}^\star+\frac{R_{ij}}{p_i} \cdot e_i$ and~$\vec{x}^\star+\frac{R_{ij}}{p_j} \cdot e_j$, since $R_{ij}/p_i\cdot v_i=R_{ij}/p_j \cdot v_j$ follows from~$q_i=~q_j$. 
Second, we show that the difference in expected compensation, \kzedit{$\mathbb{E}[W_B\left(X+X_i\right)] - \mathbb{E}[W_B\left(X+X_j\right)]$}, is negative.

We apply the following probabilistic bound on binomial random variables.

\begin{proposition}[Corollary 4 in \cite{pinelis2020monotonicity}]\label{pr:binommean}
Take any $\lambda \in (0,\infty), n\in\mathbb{N}$, and denote $Y_n \sim Bin\left(n,\lambda/n\right)$.
\begin{enumerate}[label=(\roman*)]
    \item If $m \geq 1+\lambda$, then $\mathbb{P}\left[Y_n \geq m\right]$ is strictly increasing in $n$, i.e., 
    $$\mathbb{P}\left[Y_{n_1} \geq m\right] < \mathbb{P}\left[Y_{n_2} \geq m\right] \text{ if } n_1 < n_2.$$
    \item If $m \leq \lambda$, then $\mathbb{P}\left[Y_n \geq m\right]$ is strictly decreasing in $n$, i.e., 
    $$\mathbb{P}\left[Y_{n_1} \geq m\right] > \mathbb{P}\left[Y_{n_2} \geq m\right] \text{ if } n_1 < n_2.$$
\end{enumerate}
\end{proposition}
With $X_i \sim Bin\left(\frac{R}{p_i},p_i\right), X_j \sim Bin\left(\frac{R}{p_j},p_j\right)$, we can apply this result by setting $\lambda=R_{ij}, n_i=R_{ij}/p_i, n_j=R_{ij}/p_j$. Then, we find that for either $R'=R_{ij}$ or $R'=R_{ij}+1$
$$F_{X_i}(x) =1-\mathbb{P}[X_i \geq x+1] > 1-\mathbb{P}[X_j \geq x+1]= F_{X_j}(x) \text{ when } x > R^\prime,$$ $$\text{and } F_{X_i}(x)=1-\mathbb{P}[X_i \geq x+1] \leq 1-\mathbb{P}[X_j \geq x+1]=F_{X_j}(x) \text{ when } x \leq R^\prime.$$
To show that the difference in expected compensation is strictly negative, we will prove that
$$\forall a \geq 0:\kzedit{\mathbb{E}[W_B\left(X+X_i\right)|X = a]} = \mathbb{E}\left[\left(X_i - (B-a)\right)^+\right]$$ $$\leq \mathbb{E}\left[\left(X_j - (B-a)\right)^+\right] = \kzedit{\mathbb{E}[W_B\left(X+X_j\right)|X = a]}, $$
and the inequality is strict for some $a$ with $\mathbb{P}\left[X=a\right] > 0$. We distinguish between three cases: 
\begin{enumerate}[label=(\roman*)]
    \item $E_1 = \left\{a:B - a \leq 0 \right\}$;
    \item $E_2 = \left\{a:0 < B-a \leq R^\prime\right\}$;
    \item $E_3 = \left\{a:R^\prime < B-a\right\}$.
\end{enumerate}
\noindent Note that we may assume WLOG in this discussion that $R^\prime \leq B$ because we define $R^\prime$ as a constant that depends only on $p_i$, $p_j$, not on $B$.

In event $E_1$, there is no capacity left after $a$ customers show up, so 
$$\kzedit{\mathbb{E}[W_B\left(X+X_i\right)|X \in E_1] - \mathbb{E}[W_B\left(X+X_j\right)|X \in E_1]} = \mathbb{E}[X_i] - \mathbb{E}[X_j]{= R-R} = 0.$$

\noindent For events $E_2$ and $E_3$, we first establish a result for any $X_0 \sim Bin(n,p)$, where $n \in \mathbb{Z}^+$ and $0 < p < 1$. Given $B \geq 0$,
\begin{eqnarray*}
\kzedit{\mathbb{E}[W_B\left(X_0\right)]} = \mathbb{E}\left[(X_0-B)^+\right] = \sum_{x = B+1}^{n} \left(x-B\right) \mathbb{P}\left[X_0 = x\right] = \sum_{x = 1}^{n-B} x \mathbb{P}\left[X_0 = B+x\right] \\
= \sum_{x = 1}^{n-B} \sum_{y = x}^{n-B} \mathbb{P}\left[X_0 = B+y\right]  = \sum_{x = 0}^{n-B} \left(1-F_{X_0}(B+x)\right) = \sum_{x = B}^{n} \left(1-F_{X_0}(x)\right).
\end{eqnarray*}
For any $a \in E_2$, the remaining capacity is $B-a > 0$. We plug in $$n = \frac{R_{ij}}{p_i} \text{ for } \kzedit{\mathbb{E}[W_B\left(X+X_i\right)|X = a]} \text{ and } n=\frac{R_{ij}}{p_j} \text{ for } \kzedit{\mathbb{E}[W_B\left(X+X_j\right)|X = a]}.$$

Then,
\begin{align*}
    &\kzedit{\mathbb{E}[W_B\left(X+X_i\right)|X = a] - \mathbb{E}[W_B\left(X+X_j\right)|X = a]}\\ 
    = &\sum_{x = B-a}^{\frac{R_{ij}}{p_i}} \left(1-F_{X_i}(x)\right) - \sum_{x = B-a}^{\frac{R_{ij}}{p_j}} \left(1-F_{X_j}(x)\right)\\
    =& \mathbb{E}[X_i] - \sum_{x = 0}^{B-a} \left(1-F_{X_i}(x)\right) - \left[\mathbb{E}[X_j]-\sum_{x = 0}^{B-a} \left(1-F_{X_j}(x)\right)\right] \\
    =& \sum_{x = 0}^{B-a} \left(F_{X_i}(x) - F_{X_j}(x)\right) < 0.
\end{align*}

For any $a \in E_3$,
$$\kzedit{\mathbb{E}[W_B\left(X+X_i\right)|X = a] - \mathbb{E}[W_B\left(X+X_j\right)|X = a]}$$
$$= \sum_{x = B-a}^{\frac{R_{ij}}{p_j}} \left(1-F_{X_i}(x)\right) - \sum_{x = B-a}^{\frac{R_{ij}}{p_j}} \left(1-F_{X_j}(x)\right) = \sum_{x = B-a}^{\frac{R_{ij}}{p_j}} \left(F_{X_j}(x) - F_{X_i}(x)\right) < 0.$$

Thus, overall we have

\begin{equation*}
\begin{split}
& \kzedit{\mathbb{E}[W_B\left(X+X_i\right)] - \mathbb{E}[W_B\left(X+X_j\right)]}\\
=& \sum_{a = B}^{\infty} \mathbb{P}[X = a] \left[\kzedit{\mathbb{E}[W_B\left(X+X_i\right)|X = a] - \mathbb{E}[W_B\left(X+X_j\right)|X = a]}\right] \\
& + \sum_{a = B-R^\prime}^{B-1} \mathbb{P}[X = a] \left[\kzedit{\mathbb{E}[W_B\left(X+X_i\right)|X = a] - \mathbb{E}[W_B\left(X+X_j\right)|X = a]}\right] \\
& + \sum_{a = 0}^{B-R^\prime-1} \mathbb{P}[X = a] \left[\kzedit{\mathbb{E}[W_B\left(X+X_i\right)|X = a] - \mathbb{E}[W_B\left(X+X_j\right)|X = a]}\right]\\
<& \sum_{a = B}^{\infty} \mathbb{P}[X = a] \cdot 0 + \sum_{a = B-R^\prime}^{B-1} \mathbb{P}[X = a] \cdot 0 + \sum_{a = 0}^{B-R^\prime-1} \mathbb{P}[X = a] \cdot 0 = 0.
\end{split}
\end{equation*}

The strict inequality implies that the difference in expected compensation is strictly negative. This completes the proof of Lemma \ref{lem:exchange}. \halmos

\subsection{Proof of Lemma \ref{lem:global_sensitivity} (ii)}\label{app:global_sensitivity}
The proof is almost symmetric to the first part. Instead of Inequality \eqref{eq:delta_bound}, we show 
$$x_j^\star - x_j^\prime \leq \delta \sum_{n =1}^{\hat{j}} (x_n^\prime - x_n^\star), \forall j > \hat{j}.$$
Since $\vec{x}^\prime$ is locally optimal at every $i \geq \widetilde{j}$, we now show that if $x_i^\star - x_i^\prime > \delta \sum_{n = 1}^{\hat{j}} (x_n^\prime - x_n^\star)$ for some type $i > \hat{j}$, then $\vec{x}^\prime$ cannot be locally optimal at $i$, which is a contradiction.

Let $m = \sum_{j =1}^{\hat{j}} (x_j^\prime - x_j^\star)$. Similar to before, we construct a sequence of solutions $$\vec{x}^1,\vec{x}^2,...,\vec{x}^{m+1},$$ all of which are locally optimal at $i$ but vary in the number of customers of other types. Let $\vec{x}^1 = \vec{x}^\star$, which we know is locally optimal at $i$. For each $h \in \{1,...,m\}$, we add a customer of type $j \leq \hat{j}$: we start with a solution $\vec{x}^{h}$ that is locally optimal at $i$, and apply Lemma \ref{lem:local_sensitivity} to find $l_{h+1}\leq\delta$ such that solution $\vec{x}^{h+1} = \vec{x}^{h} + e_j - l_{h+1} e_i$ is locally optimal at $i$. 
After adding all $m=\sum_{n =1}^{\hat{j}} (x_n^\prime - x_n^\star)$ customers, we find solution $\vec{x}^{m+1}$ that 
satisfies local optimality at $i$,  $x_i^{m+1} > x_i^\prime$, $x_j^{m+1}=x_j'\forall j\leq \hat{j}$, and $x_j^{m+1}=x_j^\star\geq x_j'\forall j>\hat{j}$. 

Now we argue as in the first part: with $x_j^{m+1}\geq x_j'\forall j$ we have the first, and with $\vec{x}^{m+1}$ locally optimal at $i$ (from Lemma \ref{lem:local_optimality}) the second inequality in
$$\mathbb{P}\left[\sum_{n} X_{n,x_n'} \geq B \right]\leq\mathbb{P}\left[\sum_{n \neq i} X_{n,x_n^{m+1}} + X_{i,x_i^{m+1}-1} \geq B \right] \leq q_i.$$
Since $x_i'<N_i$ this contradicts, by Lemma \ref{lem:local_optimality} (i), that $\vec{x}'$ is locally optimal at $i$. \halmos


\subsection{Proof of Claim \ref{cl:final_lemma}}\label{app:cl_final_lemma}
Let $\widetilde{j}$ be the threshold index of index solution $\vec{x}'$. Since $x_j'=\tildeN_j\kzedit{[}t-1\kzedit{]}\frall j<i$, we know $\widetilde{j} \geq i$. From the proof of Lemma \ref{lem:global_sensitivity} we know

\begin{enumerate}[label=(\roman*)]
    \item If $\widetilde{j} = i$, $\vec{x}'$ is locally optimal at $i$
    \item If $\widetilde{j} > i$, $\vec{x}' - e_{\widetilde{j}}$ is locally optimal at $i$
\end{enumerate}

In particular, $\vec{x}' - e_{\widetilde{j}}$ is a feasible solution because $x_{\widetilde{j}}' \geq 1$ is guaranteed for threshold index $\widetilde{j}$ as defined. Here we assume without loss of generality that $\widetilde{j} > i$, and the proof of (i) follows from the same argument.

For index solution $\vec{x}^\star$, similarly, we have threshold index $\hat{j}$ with $x_{\hat{j}}^\star>0$. Since $x_j^\star=0$ for every~$j \geq i$, we know $\hat{j} < i$. Thus, recall from the proof of Lemma \ref{lem:global_sensitivity} that $\vec{x}^\star$ is locally optimal at $i$. 

To prove Claim \ref{cl:final_lemma}, 
suppose that $x_i^\prime - x_i^\star > \delta \sum_{n = 1}^{i-1} (x_n^\star - x_n^\prime)$. Then  $x_i^\star < N_i$ and we derive a contradiction by showing that
\begin{equation}\label{eq:loc_opt}
\mathbb{P}\left[\sum_{n} X_{n,x_n[t-2]} + \sum_{n} X_{n,x_n^\star} \geq B \right]  \leq q_i,
\end{equation}
which means that $\vec{x}^\star$ does not fulfill the local optimality conditions from  Lemma \ref{lem:local_optimality} (i) at $i$. 

Let $m = \sum_{n = 1}^{i-1} (x_n^\star - x_n^\prime)$. We construct a sequence of solutions $\vec{x}^1,\vec{x}^2,\vec{x}^3,...,\vec{x}^{m+1}$ where~$\vec{x}^{m+1}$ has the following properties:
\begin{itemize}
    \item $\vec{x}^{m+1}$ is locally optimal at $i$,
    \item $x_i^{m+1} > x_i^\star$,
    \item $x_j^{m+1} \geq x_j^\star$ for all $j > i$,
     and
    \item $x_j^{m+1} = x_j^\star$ for all $j < i$.
\end{itemize}

Thus, $$\mathbb{P}\left[\sum_{n} X_{n,x_n[t-2]}+\sum_{n} X_{n,x_n^\star} \geq B \right] \leq \mathbb{P}\left[\sum_{n} X_{n,x_n[t-2]}+\sum_{n \neq i} X_{n,x_n^{m+1}} + X_{i,x_i^{m+1}-1} \geq B \right] \leq q_i,$$ 

where the first inequality comes from the fact that $x_j^{m+1}\geq x_j^\star\forall j$, and the second inequality comes from the local optimality of $\vec{x}^{m+1}$ at $i$ (allowing us to apply Lemma \ref{lem:local_optimality} (ii)).

Therefore, such $\vec{x}^{m+1}$ implies Inequality \eqref{eq:loc_opt}, i.e., a contradiction to $\vec{x}^\star$ being locally optimal at~$i$. We derive $\vec{x}^{m+1}$ by inductively constructing $\vec{x}^h$ from $\vec{x}^{h-1}$, where $\vec{x}^h$  is locally optimal at $i$ for every $h \in \{1,2,...,m+1\}$, $x_i^{h} > x_i^\star + \delta (m+1-h)$, and $\vec{x}^{m+1}$ fulfills the last two properties above.

Recall that $\vec{x}' - e_{\widetilde{j}}$ is locally optimal at $i<\widetilde{j}$, so we start the construction from $\vec{x}^1 = \vec{x}'- e_{\widetilde{j}}$. Then, for each $h \in \{1,...,m\}$, we repeat the procedure of adding one customer of type $j < i$: given a solution $\vec{x}^{h}$ that is locally optimal at~$i$, Lemma \ref{lem:local_sensitivity} allows us to find the $l_{h+1}$ such that $$\vec{x}^{h+1} = \vec{x}^{h} + e_j - l_{h+1} e_i$$ is locally optimal at~$i$. Moreover, Lemma \ref{lem:local_sensitivity} guarantees that $l_{h+1} \leq \delta$, i.e., $x_i^{h} - x_i^{h+1} \leq \delta$. Thus, the first two properties are maintained throughout. Since we have $x^1_j \geq x_j^\star = 0$ for $j > i$ and we never remove a customer from such $j$, the third property holds. Since we add customers until $x_j^{m+1}=x^\star_j$ for $j<i$ the fourth property holds. This completes the proof of Claim \ref{cl:final_lemma}. \halmos

\section{Model extensions}\label{app:extensions}

In this appendix we extend our results in two separate directions: first, we show that type-dependent refunds for no-shows can be seamlessly incorporated into our results. Next, we argue that with type-dependent resource demand (as opposed to each arrival requiring just one unit of resources), our results continue to hold under an additional technical assumption.

\subsection{Model with product-dependent refunds for no-show} \label{app:refund}

Consider a model in which a no-show customer of type $j$ is given a refund of $0 \leq r_j < v_j$. Define $\bar{v}_j = v_j - r_j (1-p_j)$, which denotes the expected revenue from a customer of type $j$ in this new model. We show that all results in Section \ref{sec:main} naturally extend to this model by replacing $v_j$ with $\bar{v}_j$ in the proof.\footnote{This is a well-known technique, notably mentioned in Section 3.3.2 of \cite{gallego2019revenue}, for introducing product-dependent no-show refunds into the model.}

Since $v_j$'s are not used in the proofs of Lemma \ref{lem:local_sensitivity}, Lemma \ref{lem:small_loss_probability} and Theorem \ref{thm:coupling}, the only results in Section \ref{sec:main} that need to be updated are Lemma \ref{lem:local_optimality} and Theorem \ref{thm:index_loss}. 

For Lemma \ref{lem:local_optimality}, we now define $q_j = \frac{\bar{v}_j}{p_j}$. Then, in Lemma \ref{lem:local_optimality} (i), $\mathbb{P}\left[\sum_{n} X_{n,x_n} \geq B\right] > q_i$, by linearity of expectation, still implies that accepting $x_i^\star + 1$ customers of type $i$  yields less profit than accepting~$x_i^\star$. A similar argument can be made for Lemma \ref{lem:local_optimality} (ii). The formal proof is exactly the same as in Section \ref{sec:local_analysis} with this pair of new definitions of $q_j$ and $\bar{v}_j$.

For Theorem \ref{thm:index_loss}, similarly, by replacing $v_j$ with $\bar{v}_j$ all analyses on expected revenue still hold by linearity of expectation.

\subsection{Model with Type-dependent Resource Demand}\label{app:multi-unit}
We next consider a model in which a customer of type $j$ demands $d_j$ units of capacity, where $d_j \in \mathbb{Z}^+$. The $d_j$ units of resource demand must be served entirely to avoid rejecting the customer at departure. The compensation for rejecting a customer of type $j$ is proportional to the resource requirement $d_j$. In particular, since we have assumed that the compensation for a customer who demands a single unit of capacity to be $1$, here we assume that the compensation for rejecting a customer of type $j$ is $d_j.$

For our analysis in this part we re-define the critical ratio $q_j:=\frac{v_j}{d_j p_j}$, and~$\bar{q_j}:=1-\frac{v_j}{d_j p_j}$ as  its complement accordingly. Assume that the types are ordered such that $1 >q_1 > q_2 > \ldots > q_k.$ Intuitively, since $d_j p_j$ is the expected units of resource occupied by a customer of type $j$, $q_j$ captures the expected revenue from a customer of type $j$ per unit of capacity occupied by him or her. Notice that we now introduced the following technical condition which we rely on in the proof of Lemma~\ref{lem:exchange}~(Ext), and which our results in this appendix require.
\begin{assumption}
$q_i\neq q_j$ for every $i\neq j$.
\end{assumption}
As before, the index policy prioritizes customers with higher critical ratio over ones with lower critical ratios. With this in hand, we are ready to formally generalize the definitions and results in Section~\ref{sec:main}.

\begin{definition*}[Definition~\ref{defn:local_opt} (Ext)]
Consider a feasible solution $\vec{x}^\star$ to \eqref{eq:offline} with (estimated) future arrivals $\tildeN\kzedit{[}t\kzedit{]}$ and past accepted arrivals $\vec{x}[t-1]$. We call $\vec{x}^\star$ \emph{locally optimal at} $i$ if 
$$
x_i^\star=\max\argmax_{x_i:0\leq x_i \leq \tildeN_i\kzedit{[}t\kzedit{]}}\left\{
x_iv_i-
\mathbb{E}\left[\left(d_i X_{i,x_i[t-1]+x_i}+\sum_{j\neq i} d_j X_{j,x_j[t-1]+x_j^\star}-B\right)^+\right] \right\}.
$$
\end{definition*}

With this local optimality condition in hand, we extend Lemma \ref{lem:local_optimality} to the case with type-dependent resource demand.

\begin{lemma*}[Lemma~\ref{lem:local_optimality} (Ext)]
For a solution $\vec{x}^\star$ to \eqref{eq:offline} with (estimated) future arrivals $\tildeN\kzedit{[}t\kzedit{]}$ and past accepted arrivals ~$\vec{x}[t-1]$, let $\vec{x} = \vec{x}[t-1] + \vec{x}^\star$. Then, $\vec{x}^\star$ is locally optimal at $i$ if and only if both of the following conditions are satisfied:
\begin{enumerate}[label=(\roman*)]
    \item $x_i^\star = \tildeN_i\kzedit{[}t\kzedit{]}$ or $\mathbb{P}\left[\sum_{n} d_n X_{n,x_n} \geq B-d_i + 1\right] > q_i$;
    \item $x_i^\star = 0$ or $\mathbb{P}\left[\sum_{n \neq i} d_n X_{n,x_n} + d_i X_{i,x_i -1} \geq B-d_i+1\right] \leq q_i$.
\end{enumerate}
\end{lemma*}

The proof of Lemma \ref{lem:local_optimality} can be extended to show Lemma~\ref{lem:local_optimality} (Ext) with minor changes. Specifically, these changes include multiplying each binomial random variables $X$ by the corresponding resource requirement, changing $B$ to $N-d_n+1$ and plugging in the new definition of $q_n$ for all $n$. In a model with type-dependent resource demand, similar to the model without, there always exists some global optimum that is locally optimal for every type $j$. This follows from an extension of Appendix \ref{app:glocal} with the same changes as in Lemma~\ref{lem:local_optimality}. Thus, we may again compare our solution to a global optimum that fulfills the local optimality condition for every type. 

The next two results, Lemma \ref{lem:local_sensitivity} (Ext) and Theorem \ref{thm:index_loss} (Ext), continue to hold true as stated in the main body of the paper, though the proofs require additional work. For clarity of exposition, we restate the results here and provide proof sketches in Appendix \ref{app:local_sensitivity_ext} and Appendix \ref{app:index_loss_ext}, respectively.

\begin{lemma*}[Lemma~\ref{lem:local_sensitivity} (Ext)]
There exists a constant $\delta \geq 1$ such that for any solution $\vec{x}^\star$ to \eqref{eq:offline} with $\tildeN\kzedit{[}t\kzedit{]}$ and ~$\vec{x}[t-1]$, it is true that if $\vec{x}^\star$ is locally optimal at $j$, then $\exists l \in \{0,1,...,\delta\}$ such that $\vec{x}^\prime = \vec{x}^\star + e_i - l e_j$ is locally optimal at $j$.
\end{lemma*}

\begin{theorem*}[Theorem \ref{thm:index_loss} (Ext)]
There exists a constant $M_1$ such that
$$
\mathbb{E}_{\vec{A}}\left[\widehat{\opt}_{\vec{A}}-\calH_{\vec{A}}\kzedit{[}1\kzedit{]}\right]\leq M_1.
$$
\end{theorem*}

We next state the generalization of Lemma \ref{lem:small_loss_probability} to the model with type-dependent resource requirements; this has the same statement as in the original model. We omit the proof as it follows the same reasoning as in Section \ref{sec:lemma3} with the sole difference being due to constants changing in the generalizations of our claims and lemmas.

\begin{lemma*}[Lemma~\ref{lem:small_loss_probability} (Ext)]
There exists a constant $M_2$ such that
$$
\sum_{t=2}^T \mathbb{P}\left[\calH_{\vec{A}}\kzedit{[}t-1\kzedit{]}-\calH_{\vec{A}}\kzedit{[}t\kzedit{]}>0\right]\leq M_2.
$$
\end{lemma*}

From these bounds we derive, again, using compensated coupling \cite{vera2019bayesian,freund2019uniform}, the result in Theorem~\ref{thm:coupling}~(Ext).

\begin{theorem*}
[Theorem~\ref{thm:coupling} (Ext)]
There exists a constant $M$ such that $\mathbb{E}_{\vec{A}}\left[\widehat{\opt}_{\vec{A}}-\obj_{\vec{A}}\right] \leq M$.
\end{theorem*}

\emph{Proof Sketch of Theorem \ref{thm:coupling} (Ext)}.
We change the starred inequality in the proof of Theorem~\ref{thm:coupling} to $$M_1 + \mathbb{E}_{\vec{A}}\left[\sum_{t=2}^T\calH_{\vec{A}}\kzedit{[}t-1\kzedit{]}-\calH_{\vec{A}}\kzedit{[}t\kzedit{]}\right] \leq M_1 + d_{max} \sum_{t=2}^T \mathbb{P}\left[\calH_{\vec{A}}\kzedit{[}t-1\kzedit{]}-\calH_{\vec{A}}\kzedit{[}t\kzedit{]}>0\right]$$ because any action in one period affects the objective, through lost revenue or incurred compensation, by at most $d_{max} := \max_j d_j$. The other parts of the proof of Theorem \ref{thm:coupling} (Ext) remain exactly the same as before. \halmos

This concludes all main results in the model with type-dependent resource demand.

\subsubsection{Proof of Lemma~\ref{lem:local_sensitivity} (Ext)}\label{app:local_sensitivity_ext}

Since Lemma \ref{lem:local_sensitivity} is built upon Claim \ref{cl:berryesseen} and Claim \ref{cl:exchange}, we first provide generalizations of these claims. We denote for a given constant~$\delta^\star \geq 1$ and solution $\vec{x}$: $$Y = \sum_{n \neq j} d_n X_{n,x_n}+ d_j X_{j,x_j+\delta^\star}, Z = \sum_{n} d_n X_{n,x_n}, \text{ and } m = \sum_{n = 1}^k x_n.$$

\begin{claim*}[Claim~\ref{cl:berryesseen} (Ext)]
There exist constants $\delta_{j}^\star$ and $m_0 \in \kzedit{\mathbb{Z}^+}$ such that
\begin{equation}\label{diff_ext}
    F_{Y}^{-1} (\bar{q}_j)-F_{Z}^{-1} (\bar{q}_j) \geq d_i
\end{equation}
holds whenever $m \geq m_0$.
\end{claim*}

\begin{claim*}[Claim~\ref{cl:exchange} (Ext)]
If $\vec{x} - (\delta_j + 1)e_j + e_i$ is a feasible solution to the optimization problem \eqref{eq:offline}, then
\begin{equation*}
\begin{split}
    &\mathbb{P}\left[\sum_{n \neq i,j} d_n X_{n,x_n}+ d_j X_{j,x_j-\delta_{j}-1} \geq B-d_j + 1-d_i X_{i,x_i+1}\right] \\
    \leq& \mathbb{P}\left[\sum_{n \neq i,j} d_n X_{n,x_n}+d_j X_{j,x_j-\delta_{j}-1} \geq B-d_j+1 -d_i-d_i X_{i,x_i}\right]
\end{split}
\end{equation*}
\end{claim*}

\emph{Proof Sketch of Claim \ref{cl:berryesseen} (Ext)}.
To establish Claim \ref{cl:berryesseen} (Ext), we need to scale the equations accordingly since the random variables are scaled by $d_n$'s.

Specifically, in the proof of Claim \ref{cl:berryesseen} (Ext) we now have $$\mu_Z = \mathbb{E}[Z] = \sum_{n = 1}^k d_n x_n p_n \text{ and } \sigma_Z^2 = \text{Var}[Z] = \sum_{n = 1}^k d_n^2 x_n p_n (1-p_n).$$
Similarly, $\mu_Y = \mathbb{E}[Y] = \delta^\star d_j p_j + \sum_{n = 1}^k d_n x_n p_n$, and $$\sigma_Y^2 = \text{Var}[Y] =\delta^\star d_j^2 p_j (1-p_j) + \sum_{n = 1}^k d_n^2 x_n p_n (1-p_n).$$

Moreover, for fixed $\delta^\star$, we re-define for $r \in \{1,...,m+\delta^\star\}$, where the $r$th customer is of type $n$, the random variable $Z_r := d_n (X_{n,1} - p_n)=d_n(X_{n,1} - p_r)$. Then, $\mathbb{E}[Z_r] = 0$ for~$r=1,\ldots,m$ and $$Z = \sum_{r = 1}^m Z_r + \sum_{n = 1}^k d_n x_n p_n = \sum_{r = 1}^m Z_r + \mu_Z.$$
Similarly, we re-define the random variables $Y_r := d_n(X_{n,1} - p_n) = d_n(X_{n,1} - p_r)$, so that~$\mathbb{E}[Y_r] = 0$ for~$r=1,\ldots,m+\delta^\star$ and 
$$Y = \sum_{r = 1}^{m+\delta^\star} Y_r + \delta^\star d_j p_j + \sum_{n = 1}^k d_n x_n p_n = \sum_{r = 1}^{m+\delta^\star} Y_r  + \mu_Y.$$

Finally, we normalize $Z$ and $Y$ by defining $$Z^\prime = \frac{Z-\mu_Z}{\sigma_Z} = \frac{\sum_{r = 1}^m Z_r}{\sigma_Z}, Y^\prime = \frac{Y-\mu_Y}{\sigma_Y} = \frac{\sum_{r = 1}^{m+\delta^\star} Y_r}{\sigma_Y}.$$ Then we can apply Proposition \ref{pr:berryesseen} to the re-defined variables and the proof of Claim \ref{cl:berryesseen} (Ext) follows with $$\delta^\star = \frac{d_i+ 2+12 \Phi^{-1 \prime}(\bar{q}_j)}{d_j p_j}. \halmos$$ 

\emph{Proof Sketch of Claim \ref{cl:exchange} (Ext)}.
The proof of Claim \ref{cl:exchange} (Ext) follows the same reasoning  as that of Claim \ref{cl:exchange} once we observe that $$d_i X_{i,x_i+1} = d_i X_{i,x_i} + d_i X_{i,1} \leq d_i X_{i,x_i} + d_i. \halmos$$ 

\emph{Proof Sketch of Lemma \ref{lem:local_sensitivity} (Ext)}.
With the extended claims established, the proof of Lemma~\ref{lem:local_sensitivity} requires only minor changes. Specifically, we multiply all binomial random variable $X$'s by the corresponding resource requirement, update the local optimality conditions based on Lemma \ref{lem:local_optimality} and change $B-1$ to~$B-d_i$ throughout the analysis. \halmos 

\subsubsection{Proof of Theorem \ref{thm:index_loss} (Ext)}\label{app:index_loss_ext}

For the proof of Theorem \ref{thm:index_loss} (Ext) we let $d_{max} := \max_j d_j$ and $d_{min} := \min_j d_j$. To establish the result we first generalize Lemmas \ref{lem:exchange}---\ref{lem:global_sensitivity}.

\begin{lemma*}[Lemma~\ref{lem:exchange} (Ext)]
For every $i < j$ there exists a constant $R_{ij}$ such that, for any feasible solutions $\vec{x}^\star+\frac{R_{ij}}{d_i p_i} \cdot e_i$ and $\vec{x}^\star+\frac{R_{ij}}{d_j p_j} \cdot e_j$ to optimization problem \eqref{eq:offline} with $N\kzedit{[}t,T\kzedit{]}$ and $\vec{x}[t-1]$, the objective of $\vec{x}^\star+\frac{R_{ij}}{d_i p_i} \cdot e_i$ is always greater than that of $\vec{x}^\star+\frac{R_{ij}}{d_j p_j} \cdot e_j$.
\end{lemma*}

\emph{Proof Sketch of Lemma \ref{lem:exchange} (Ext)}.
For any constant $R \in \mathbb{Z}^+$ such that $\frac{R}{d_i p_i}$ and $\frac{R}{d_j p_j}$ are both integers, denote $$X_i \sim d_i Bin\left(\frac{R}{d_i p_i},p_i\right), X_j \sim d_j Bin\left(\frac{R}{d_j p_j},p_j\right).$$ Since we leave all customers among $\vec{x} = \vec{x}^\star + \vec{x}[t-1]$ unchanged, we denote the number of unchanged customers who consume a resource by a random variable $X \sim \sum_{n} d_n X_{n,x_n}.$

As before, the difference in revenue between any feasible solutions $\vec{x}^\star+\frac{R}{d_i p_i} \cdot e_i$ and $\vec{x}^\star+\frac{R}{d_j p_j} \cdot e_j$ is 
$$\frac{R}{d_i p_i} v_i  - \frac{R}{d_j p_j} v_j = \left(\frac{v_i}{d_i p_i}-\frac{v_j}{d_j p_j}\right)R \in \Omega(R).$$

Then, what is left to show is that $$\kzedit{\mathbb{E}[W_B\left(X+X_i\right)] - \mathbb{E}[W_B\left(X+X_j\right)]} \leq \mathbb{E}\left[(X_i-X_j)^+\right]\in O\left(\sqrt{R\log(R)}\right).$$ This follows with the same application of a Chernoff bound as before, setting $$\kappa=2 \max(d_i,d_j) \sqrt{R\log(R)}.$$ Putting everything together, we observe that there exists $R_{ij}$ such that the difference in objective between $\vec{x}^\star+\frac{R}{d_i p_i} \cdot e_i$ and $\vec{x}^\star+\frac{R}{d_j p_j} \cdot e_j$ is at least $$\left(\frac{v_i}{d_i p_i}-\frac{v_j}{d_j p_j}\right)R - O\left(\sqrt{R\log(R)}\right) > 0.\halmos$$

\begin{lemma*}[Lemma~\ref{lem:chernoff_difference} (Ext)]
Consider an optimal solution $\vec{x}^\star$ to the optimization problem \eqref{eq:offline} with future arrivals~$N\kzedit{[}t,T\kzedit{]}$, and past accepted arrivals $\vec{x}[t-1]$. Then for every $i$ and $j$ with $i<j$ there exists a constant $R_{ij}$ such that at least one of the following two is true:
\begin{enumerate}[label=(\roman*)]
    \item $x_i^\star > N\kzedit{[}t,T\kzedit{]}-\frac{R_{ij}}{d_i p_i}$ \item $x_j^\star < \frac{R_{ij}}{d_j p_j}$.
\end{enumerate}
\end{lemma*}

\emph{Proof Sketch of Lemma \ref{lem:chernoff_difference} (Ext)}.
Lemma \ref{lem:chernoff_difference} (Ext) is immediate from Lemma~\ref{lem:exchange} (Ext) since otherwise we could replace $\frac{R_{ij}}{d_j p_j}$ type $j$ customers with $\frac{R_{ij}}{d_i p_i}$ type $i$ customers to derive a contradiction that $\vec{x}^\star$ is not an optimal solution. \halmos

\begin{lemma*}[Lemma~\ref{lem:global_sensitivity} (Ext)]\label{lem:global_sensitivity_ext}
Consider an optimal solution $\vec{x}^\star$ that is locally optimal for every type~$j$, and an optimal index solution $\vec{x}^\prime$ that is locally optimal at its threshold index $\widetilde{j}$, both for the optimization problem \eqref{eq:offline} with arrivals $\vec{N}$. With $\delta$ as constructed in Lemma \ref{lem:local_sensitivity}, we have
\begin{enumerate}[label=(\roman*)]
    \item $\sum_j \left(x_j^\prime - x_j^\star\right)^+ \leq k \delta \left[2 +\sum_j \left(x_j^\star - x_j^\prime\right)^+\right]$,
    \item $\sum_j \left(x_j^\star - x_j^\prime\right)^+ \leq k \delta \left[1 +\sum_j \left(x_j^\prime - x_j^\star\right)^+\right]$. 
\end{enumerate}
\end{lemma*}
\emph{Proof Sketch of Lemma \ref{lem:global_sensitivity} (Ext)}.
We sketch the proof of Lemma~\ref{lem:global_sensitivity} (Ext) (ii), and the proof of Lemma~\ref{lem:global_sensitivity} (Ext) (i) follows the same reasoning. We shall show that for every type $j> \widetilde{j}$, we can find an integer $l \leq \delta$ such that $\vec{x}^\prime + l e_j$ is locally optimal at $j$. Then, we follow the rest of the proof of Lemma~\ref{lem:global_sensitivity} (ii) to obtain an upper bound that is a constant $k\delta$ larger than in the original Lemma~\ref{lem:global_sensitivity}~(ii).

From Lemma~\ref{lem:local_optimality}~(i), since $\vec{x}'$ is locally optimal at $\widetilde{j}$, we find that either~$x_{\widetilde{j}}^\prime = N_{\widetilde{j}},$ in which case $x_{\widetilde{j}+1}^\prime = 0$, and
$$
\mathbb{P}\left[\sum_{n} d_n X_{n,x_n^\prime} \geq B - d_{\widetilde{j}+1}+1 \right] > q_{\widetilde{j}+1},\;
\text{ or }\;
\mathbb{P}\left[\sum_{n} d_n X_{n,x_n^\prime} \geq B - d_{\widetilde{j}}+1 \right] > q_{\widetilde{j}}.$$
Now, for any $j > \widetilde{j}$, we will check that
\begin{equation}\label{eq:every_index}
    \mathbb{P}\left[\sum_{n} d_n X_{n,x_n^\prime} + d_j X_{j,\delta} \geq B - d_j + 1\right] > q_j.
\end{equation} 

When \eqref{eq:every_index} is satisfied, we know from Lemma \ref{lem:local_optimality} that for $\vec{x}^\prime+ l e_j$ to be locally optimal at~$j$ we must have $l \leq \delta$. Since $q_{\widetilde{j}} \geq q_{\widetilde{j}+1} \geq q_j$, to show \eqref{eq:every_index}, it suffices to show that
\begin{equation}\label{eq:greater_than_both}
    \begin{split}
        &\mathbb{P}\left[\sum_{n} d_n X_{n,x_n^\prime} + d_j X_{j,\delta} \geq B - d_j + 1\right]\\
        \geq& \max\left(\mathbb{P}\left[\sum_{n} d_n X_{n,x_n^\prime} \geq B - d_{\widetilde{j}+1}+1 \right],\mathbb{P}\left[\sum_{n} d_n X_{n,x_n^\prime} \geq B - d_{\widetilde{j}}+1 \right]\right) 
    \end{split}
\end{equation}

We prove \eqref{eq:greater_than_both} by discussing two cases:
\begin{enumerate}[label=(\alph*)]
    \item $d_j \geq \max(d_{\widetilde{j}},d_{\widetilde{j}+1})$
    \item $d_j < d_{\widetilde{j}}$ or $d_j < d_{\widetilde{j}+1}$
\end{enumerate}

In case (a), we get \eqref{eq:greater_than_both} from 
\begin{equation*}
    \begin{split}
        &\mathbb{P}\left[\sum_{n} d_n X_{n,x_n^\prime} + d_j X_{j,\delta} \geq B - d_j + 1\right]\\
        \geq& \mathbb{P}\left[\sum_{n} d_n X_{n,x_n^\prime} \geq B - d_j + 1\right]\\
        \geq& \max\left(\mathbb{P}\left[\sum_{n} d_n X_{n,x_n^\prime} \geq B - d_{\widetilde{j}+1}+1 \right],\mathbb{P}\left[\sum_{n} d_n X_{n,x_n^\prime} \geq B - d_{\widetilde{j}}+1 \right]\right)
    \end{split}
\end{equation*}

In case (b), let $Y = \sum_{n} d_n X_{n,x_n^\prime}+ d_j X_{j,\delta}$ and $Z = \sum_{n} d_n X_{n,x_n^\prime}.$ From Lemma \ref{lem:local_sensitivity} and Claim~\ref{cl:berryesseen} we know $$F_{Y}^{-1} (\bar{q}_j)-F_{Z}^{-1} (\bar{q}_j) \geq d_{max} > d_{max} - d_{min}.$$ That is, $$\mathbb{P}\left[Y \geq B - d_{min} + 1\right] \geq \mathbb{P}\left[Z \geq B - d_{max} + 1\right].$$

Thus,
\begin{equation*}
    \begin{split}
        &\mathbb{P}\left[\sum_{n} d_n X_{n,x_n^\prime} + d_j X_{j,\delta} \geq B - d_j + 1\right]\\
        \geq& \mathbb{P}\left[\sum_{n} d_n X_{n,x_n^\prime} + d_j X_{j,\delta} \geq B - d_{min} + 1\right]\\
        \geq& \mathbb{P}\left[\sum_{n} d_n X_{n,x_n^\prime} \geq B - d_{max} + 1\right]\\
        \geq& \max\left(\mathbb{P}\left[\sum_{n} d_n X_{n,x_n^\prime} \geq B - d_{\widetilde{j}+1}+1 \right],\mathbb{P}\left[\sum_{n} d_n X_{n,x_n^\prime} \geq B - d_{\widetilde{j}}+1 \right]\right),
    \end{split}
\end{equation*}
which completes the proof of Lemma \ref{lem:global_sensitivity} (Ext). \halmos

Note that a similar construction is needed in the proof of Claim \ref{cl:final_lemma}, affecting only the size of a constant, and is omitted for brevity. 

\emph{Proof Sketch of Theorem \ref{thm:index_loss} (Ext)}.
Equipped with Lemma~\ref{lem:exchange} (Ext), Lemma~\ref{lem:chernoff_difference} (Ext) and Lemma~\ref{lem:global_sensitivity} (Ext), we can then extend the proof of Theorem \ref{thm:index_loss} to show Theorem \ref{thm:index_loss} (Ext) by taking $$M_1 := \delta k (k+3) \cdot \frac{R_{max}}{p_{min} d_{min}} \cdot d_{max}.\halmos$$

\section{Instance-dependence}\label{app:instance_loss}


In this section we analyze the loss of policies in instances where $v_i,p_i$ are allowed to change with $T$. In Appendix \ref{app:eg_general} we provide such an instance for which  any online policy incurs a loss of $\Omega(\sqrt{T})$  due to the inherent uncertainty in arrivals. Then, in Appendix \ref{app:eg_index} we provide a different instance  where even the clairvoyant index policy, unaffected by the uncertainty of arrivals, incurs a loss of $\Omega(\sqrt{T})$; the later example highlights a limitation of index policies when the critical ratio of different types are ``too close'' to each other.

\subsection{Lower bound of loss with instance-dependent parameters}\label{app:eg_general}

Consider an example with the following three types of customers
\begin{align*}
    & \lambda_1 = \frac{1}{6}, v_1 = \frac{1}{2}, p_1 = 1\\
    & \lambda_2 = \frac{1}{3}, v_2 = \frac{1}{\sqrt{T}}, p_2 = \frac{3}{\sqrt{T}}\\
    & \lambda_3 = \frac{1}{2}, v_3 = 0, p_3 = 1.
\end{align*}

Moreover, suppose that $B = \frac{T}{6}.$ The intuition here is that we prefer customers of type $1$ over those of type $2$ and should never accept a customer of type $3$. However, due to the stochastic nature of the arrivals, we do not know how many type $1$ customers arrive and are thus likely to make mistakes in deciding on how many type $2$ customers to accept. For example, when at least~$\frac{T}{6}$ arrivals of type~1 occur, we would not want to accept any type 2 customers; but when at most~$\frac{T}{6}-\sqrt{T}$ arrivals of type 1 occur, we would want to accept $\Omega(T)$ type 2 customers. More specifically we show that even with full knowledge of the first $\frac{T}{2}$ arrivals, the expected loss of any online policy is~$\Omega(\sqrt{T}).$ We first establish bounds on the probability of the following events:
\begin{itemize}
    \item $E_1 = \left\{N_2\left[1,\frac{T}{2}\right] \geq \frac{T}{10}\right\} \cap \left\{N_2 \leq \frac{2}{5}T \right\}:$ We know that $$\mathbb{E}\left[N_2\left[1,\frac{T}{2}\right]\right] = \frac{1}{6}T \text{ and } \mathbb{E}[N_2] = \frac{1}{3}T.$$ With $\epsilon = \frac{6}{15}$ a Chernoff bound gives $$\mathbb{P}\left[N_2\left[1,\frac{T}{2}\right] \geq \frac{T}{10}\right] \geq 1-e^{-b_1 T}$$ for some constant $b_1 > 0$. Similarly we have $\mathbb{P}\left[N_2 \leq \frac{2}{5}T\right] \geq 1-e^{-b_2 T}$ for some constant $b_2 > 0$. Thus, $\mathbb{P}[E_1] \geq 1-e^{-c_1 T}$ for some constant $c_1 > 0.$ 
    \item $E_2 = \left\{ N_1\left[1,\frac{T}{2}\right] \in \left[\frac{T}{12} - \sqrt{T},\frac{T}{12}\right]\right\}:$ Since $$\mathbb{P}\left[N_1\left[1,\frac{T}{2}\right] \geq \frac{T}{12} - \sqrt{T}\right] =  \mathbb{P}\left[\frac{N_1\left[1,\frac{T}{2}\right] - \frac{T}{12}}{\sqrt{\frac{T}{2}(1-\lambda_1)\lambda_1}} \geq -\frac{\sqrt{T}}{\sqrt{\frac{T}{2}(1-\lambda_1)\lambda_1}}\right] = \mathbb{P}\left[\frac{N_1\left[1,\frac{T}{2}\right] - \frac{T}{12}}{\sqrt{\frac{T}{2}(1-\lambda_1)\lambda_1}} \geq -\frac{1}{\sqrt{\frac{(1-\lambda_1)\lambda_1}{2}}}\right]$$ and similarly $$\mathbb{P}\left[N_1\left[1,\frac{T}{2}\right] \leq \frac{T}{12}\right] = \mathbb{P}\left[\frac{N_1\left[1,\frac{T}{2}\right] - \frac{T}{12}}{\sqrt{\frac{T}{2}(1-\lambda_1)\lambda_1}} \leq 0\right],$$ 
    from the Berry{-}Esseen Theorem (\cite{feller1957introduction}, Chapter XVI.5, Theorem 2), we know
    that $$\mathbb{P}\left[N_1\left[1,\frac{T}{2}\right] \geq \frac{T}{12} - \sqrt{T}\right] \geq 1- \Phi\left(-\frac{1}{\sqrt{\frac{(1-\lambda_1)\lambda_1}{2}}}\right) - \frac{b}{\sqrt{T}}$$ and $$\mathbb{P}\left[N_1\left[1,\frac{T}{2}\right] \leq \frac{T}{12}\right] \geq \Phi(0) + \frac{b}{\sqrt{T}}$$ for some constant $b > 0.$ Thus, $\mathbb{P}[E_2] \geq c_2$ for some constant $c_2 > 0.$
    \item $E_3 = \left\{N_1\left[\frac{T}{2}+1,T\right] \leq \frac{T}{12} - 2 \sqrt{T}\right\}:$ Again, from the Berry{-}Esseen Theorem we know that $\mathbb{P}[E_3] \geq c_3$ for some constant $c_3 > 0$.
    \item $E_4 = \left\{N_1\left[\frac{T}{2}+1,T\right] \geq \frac{T}{12} + 2 \sqrt{T}\right\}:\mathbb{P}[E_4] \geq c_4$ for some constant $c_4 > 0$.
\end{itemize}

We now condition our analyses on the event $E_1 \cap E_2$. Since
\begin{align*}
    \mathbb{P}[E_1 \cap E_2] =& \mathbb{P}[E_1] + \mathbb{P}[E_2] - \mathbb{P}[E_1 \cup E_2]\\
    \geq & \mathbb{P}[E_1] + \mathbb{P}[E_2] -1\\
    \geq & c_2 - e^{-c_1 T},
\end{align*}
to show that the expected loss is $\Omega(\sqrt{T})$ it suffices to show that that is the case conditioned on~$E_1 \cap E_2$. Under $E_1 \cap E_2$, with full knowledge of the first $\frac{T}{2}$ arrivals, we show that the expected loss is $\Omega(\sqrt{T})$ no matter the decisions made for the first $\frac{T}{2}$ arrivals. Specifically, we distinguish between the following three cases that partition the set of  all possible decisions a policy can make for the first~$\frac{T}{2}$ arrivals: 
\begin{enumerate}[label = (\roman*)]
    \item $x_1\left[\frac{T}{2}\right] \leq \frac{T}{12} - 2 \sqrt{T}$
    \item $x_1\left[\frac{T}{2}\right] > \frac{T}{12} - 2 \sqrt{T}$ and $x_2\left[\frac{T}{2}\right] \leq \frac{T}{20}$
    \item $x_1\left[\frac{T}{2}\right] > \frac{T}{12} - 2 \sqrt{T}$ and $x_2\left[\frac{T}{2}\right] > \frac{T}{20}$
\end{enumerate}

In case (i), we know from event $E_3$ that with constant probability $c_3$ the number of type $1$ arrivals in the last~$\frac{T}{2}$ periods is less than $\frac{T}{12} - 2 \sqrt{T}.$ Thus, conditioning on event $E_3,$ we have $$N_1 = N_1\left[1,\frac{T}{2}\right] + N_1\left[\frac{T}{2}+1,T\right] \leq \frac{T}{12} + \frac{T}{12} - 2 \sqrt{T} = B - 2 \sqrt{T},$$ so the optimal offline solution in this case should accept all type $1$ arrivals and some type $2$ arrivals. However, since $$x_1\left[\frac{T}{2}\right] \leq \frac{T}{12} - 2 \sqrt{T} \leq N_1\left[1,\frac{T}{2}\right] - \sqrt{T},$$ in case (i) we accept at least $\sqrt{T}$ fewer type $1$ customers than would be feasible without having to pay any compensation (when not accepting any type 2 customers). Thus, even if we accept all type $2$ customers (upper bounded by $\frac{2}{5}T$ under~$E_1$) without incurring compensation for any of the type $2$ customers that we overbook, the incurred loss from not accepting $\sqrt{T}$ additional type 1 customers is at least $$v_1 \sqrt{T} - v_2 \frac{2}{5} T = \frac{1}{10} \sqrt{T} \in \Omega(\sqrt{T}).$$

In case (ii), we again condition on $E_3$ and thus have $N_1 \leq B - 2 \sqrt{T}.$ Since $$x_2\left[\frac{T}{2}\right] \leq \frac{T}{20} \leq N_2\left[1,\frac{T}{2}\right] - \frac{T}{20},$$ it is feasible to accept at least $\frac{T}{20}$ more type $2$ customers. 
The increase in profit from accepting $\frac{T}{20}$ additional type $2$ customers will be at least
\begin{align*}
    & v_2 \frac{T}{20} - \mathbb{E}\left[\left(Bin\left(x_2 + \frac{T}{20},\frac{3}{\sqrt{T}}\right) + x_1 - B\right)^+\right],
\end{align*}
which is what it would be without including the compensation for the original solution $\vec{x}.$ Then, we know from $x_1 \leq N_1 \leq B - 2 \sqrt{T}$ and $N_2 \leq \frac{2}{5}T$ that the increase in profit is lower bounded by 
$$v_2 \frac{T}{20} - \mathbb{E}\left[\left(Bin\left(\frac{2}{5}T,\frac{3}{\sqrt{T}}\right)-2 \sqrt{T}\right)^+\right].$$

Further, we derive that
\begin{equation}\label{eq:last_cher}
\begin{split}
&\mathbb{E}\left[\left(Bin\left(\frac{2}{5}T,\frac{3}{\sqrt{T}}\right)-2 \sqrt{T}\right)^+\right]\\
\leq & 3\sqrt{\sqrt{T} \log(\sqrt{T})} \cdot \mathbb{P}\left[Bin\left(\frac{2}{5}T,\frac{3}{\sqrt{T}}\right) \leq 2 \sqrt{T} + 6\sqrt{\sqrt{T} \log(\sqrt{T})}\right] \\
&+ \frac{2}{5}T \cdot \mathbb{P}\left[Bin\left(\frac{2}{5}T,\frac{3}{\sqrt{T}}\right) > 2 \sqrt{T} + 6 \sqrt{\sqrt{T} \log(\sqrt{T})}\right]\\
\leq & 3 \sqrt{\sqrt{T} \log(\sqrt{T})} + T \cdot \mathbb{P}\left[Bin\left(\frac{2}{5}T,\frac{3}{\sqrt{T}}\right) > 6 \sqrt{\sqrt{T} \log(\sqrt{T})}\right]\\
\leq & 3 \sqrt{\sqrt{T} \log(\sqrt{T})} + T \cdot \sqrt{T}^{-10}\\
\in & O\left(\sqrt{\sqrt{T} \log(\sqrt{T})}\right),
\end{split}
\end{equation}

where the third inequality follows from a Chernoff bound with $\epsilon = 5 \sqrt{\frac{\log{\sqrt{T}}}{\sqrt{T}}}$. Thus, the increase in profit from accepting $\frac{T}{20}$ additional type $2$ customers is lower bounded by
$$v_2 \frac{T}{20} - \mathbb{E}\left[\left(Bin\left(\frac{2}{5}T,\frac{3}{\sqrt{T}}\right)-2 \sqrt{T}\right)^+\right] \in \Omega(\sqrt{T}).$$

In case (iii) we condition on $E_4$, which happens with constant probability, to bound the loss as~$\Omega(\sqrt{T})$. 
Conditioning on $E_4,$ we have $$N_1 = N_1\left[1,\frac{T}{2}\right] + N_1\left[\frac{T}{2}+1,T\right] \geq \frac{T}{12} - \sqrt{T} + \frac{T}{12} + 2 \sqrt{T} = B + \sqrt{T}.$$
Observe first that with $x_1 > N_1 - \frac{3}{20} \sqrt{T} \geq B+\frac{17}{20} \sqrt{T}$ the type 1 customers by themselves already cause $\Omega(\sqrt{T})$ loss due to the fact that every type $1$ customer accepted, beyond the first~$B$, is guaranteed to incur compensation. Thus, we aim to show that when $x_1 \leq N_1 - \frac{3}{20} \sqrt{T}$ then replacing~$\frac{1}{20} T$ type $2$ customers with $\frac{3}{20} \sqrt{T}$ type $1$ customers improves profit by $\Omega(\sqrt{T})$.

Given any solution $\vec{x}$, replacing $\frac{1}{20} T$ type $2$ customers with $\frac{3}{20} \sqrt{T}$ type $1$ customers increases the expected compensation by at most
\begin{align*}
    & \mathbb{E}\left[\left(x_1 + \frac{3}{20} \sqrt{T} - Bin\left(x_2 - \frac{1}{20} T,\frac{3}{\sqrt{T}}\right) - B\right)^+\right] - \mathbb{E}\left[\left(x_1 + Bin\left(x_2,\frac{3}{\sqrt{T}}\right) - B\right)^+\right]\\
    \leq & \mathbb{E}\left[\left(\frac{3}{20} \sqrt{T} - Bin\left(x_2 - \frac{1}{20} T,\frac{3}{\sqrt{T}}\right) - Bin\left(x_2,\frac{3}{\sqrt{T}}\right) \right)^+\right],
\end{align*}
where the inequality follows from $\mathbb{E}[(a+b)^+ - (a+c)^+] \leq \mathbb{E}[(b-c)^+]$ with $a,b,c$ set respectively as~$a = x_1 - B$, $b = \frac{3}{20} \sqrt{T} - Bin\left(x_2 - \frac{1}{20} T,\frac{3}{\sqrt{T}}\right)$, and $c = Bin\left(x_2,\frac{3}{\sqrt{T}}\right).$
Thus, replacing~$\frac{1}{20} T$ type $2$ customers with $\frac{3}{20} \sqrt{T}$ type $1$ customers improves profit by at least 
\begin{align*}
    &v_1 \frac{3}{20} \sqrt{T} - v_2 \frac{1}{20} T - \mathbb{E}\left[\left(\frac{3}{20} \sqrt{T} - Bin\left(x_2 - \frac{1}{20} T,\frac{3}{\sqrt{T}}\right) - Bin\left(x_2,\frac{3}{\sqrt{T}}\right) \right)^+\right]
    \\=& \frac{3}{40} \sqrt{T} - \frac{1}{20} \sqrt{T} - \mathbb{E}\left[\left(\frac{3}{20} \sqrt{T} - Bin\left(x_2 - \frac{1}{20} T,\frac{3}{\sqrt{T}}\right) - Bin\left(x_2,\frac{3}{\sqrt{T}}\right) \right)^+\right]
    \\
    \geq& \frac{1}{40} \sqrt{T} - O(\sqrt{\sqrt{T} \log(\sqrt{T})}) \in \Omega(\sqrt{T}),
\end{align*}
where the inequality comes from the same argument as in \eqref{eq:last_cher}. 

This completes all cases of the proof. \halmos

\subsection{Limitation of the clairvoyant index policy}\label{app:eg_index}

In this section we provide an instance where even the \emph{clairvoyant} index policy, with no uncertainty about arrivals, incurs a loss of $\Omega(\sqrt{T})$. At a high level, this happens because the index policy is unable to differentiate among types that are ``too close" to each other, e.g., when the difference in their critical ratio is on the order of $\frac{1}{T}$. Interestingly, the problematic case is when they are extremely close without being equal; when they are exactly equal, the second part of Lemma \ref{lem:exchange} showed that index solutions successfully distinguish between them.

Consider an example with the following two types of customers
\begin{align*}
    & \lambda_1 = \frac{1}{2}, v_1 = \frac{1}{4}, p_1 = \frac{1}{2}\\
    & \lambda_2 = \frac{1}{2}, v_1 = \frac{1}{2}-\frac{1}{T}, p_1 = 1.
\end{align*}

Moreover, let $B = \frac{T}{6}$ be an integer. By Chernoff bound, we know that with high probability we would find $$N_1 \geq \frac{2}{5}T \text{ and } N_2 \geq \frac{T}{6},$$ so it suffices to show that the clairvoyant index policy incurs $\Omega(\sqrt{T})$ loss under this event. 


Since $\frac{v_1}{p_1} > \frac{v_2}{p_2}, \forall T,$ we know that the offline index policy accepts the type $1$ customers before accepting any type $2$ customers. 
We claim that the clairvoyant index solution is given by $x_1 = \frac{T}{3}$ and $x_2 = 0.$ To see this, observe that by symmetry we have $$\mathbb{P}\left[Bin\left(\frac{T}{3},\frac{1}{2}\right) \geq \frac{T}{6}\right]=1/2 = \frac{v_1}{p_1}.$$ As a result, $$\mathbb{P}\left[Bin\left(\frac{T}{3}+1,\frac{1}{2}\right) \geq \frac{T}{6}\right]>\frac{v_1}{p_1}$$ and $$\mathbb{P}\left[Bin\left(\frac{T}{3}-1,\frac{1}{2}\right) \geq \frac{T}{6}\right]<\frac{v_1}{p_1}.$$ Then, by the local optimality condition in Lemma \ref{lem:local_optimality} we conclude that $x_1 = \frac{T}{3}$ and $x_2 = 0$ is locally optimal at type $1$. Since $x_1 < N_1,$ no type $2$ customer should be accepted and this clairvoyant index solution is unique and globally optimal. 

Now, a solution that accepts $\frac{T}{6}$ type $2$ customers and no type $1$ customers, achieves $\frac{T}{12}-\frac{1}{6}$ revenue with no compensation. Thus, we aim to show that the clairvoyant index solution (which obtains~$\frac{T}{12}$ revenue) incurs an expected compensation of $\Omega(\sqrt{T})$. Indeed,

\begin{align*}
    &\mathbb{E}\left[\left(Bin\left(\frac{T}{3},\frac{1}{2}\right) - B\right)^+\right]\\
    \geq &\sqrt{T} \cdot \mathbb{P}\left[Bin\left(\frac{T}{3} - \sqrt{T},\frac{1}{2}\right) - B \geq \sqrt{T}\right]\\
    = & \sqrt{T} \cdot \mathbb{P}\left[Bin\left(\frac{T}{3} - \sqrt{T},\frac{1}{2}\right) \geq \frac{T}{6} - \frac{\sqrt{T}}{2}+ \frac{3}{2}\sqrt{T}\right]\\
    \geq & \sqrt{T} \cdot c \in \Omega(\sqrt{T})
\end{align*}
for some constant $c > 0$ based on the Berry-Esseen Theorem. This demonstrates that index solutions may incur $\Omega(\sqrt{T})$ loss, even in the absence of uncertainty, when different types are \emph{too close} in their critical ratios.

\section{Proof of $\Omega(T)$ loss using policy in \cite{erdelyi2010dynamic}}\label{app:dpd_loss}

We first formally define the policy proposed in \cite{erdelyi2010dynamic} by providing the DLP (Deterministic Linear Program) and DPD (Dynamic Programming Decomposition) formulations in \cite{erdelyi2010dynamic} adapted to our notation.  The DLP policy in \cite{erdelyi2010dynamic} relies on the following deterministic linear relaxation, where $z_j$ stands for the expected number of type $j$ customers to accept over the course of the planning horizon, and $w_j$ for the number of accepted type $j$ customers to whom the policy plans to deny service. 

\begin{equation*}
\begin{array}{ll@{}ll}
\text{maximize}  & \displaystyle\sum\limits_{j=1}^{k} v_{j}z_{j} - \displaystyle\sum\limits_{j=1}^{k} w_j \\
\text{subject to}& \displaystyle\sum\limits_{j=1}^{k} (p_j z_j - w_j) \leq B, \\
&z_{j} \leq \displaystyle\sum\limits_{t=1}^{T} \lambda_j, \forall j \in [k]\\
&w_{j} \leq p_j z_j, \forall j \in [k]\\
&z_j,w_j \geq 0, \forall j \in [k]
\end{array}
\end{equation*}

Let $\Pi^\star$ be the optimal dual value of the dual variable associated with the first constraint in the linear program. Then, the DLP policy of \cite{erdelyi2010dynamic} can be expressed as follows.\footnote{The sole difference to the policy described by \cite{erdelyi2010dynamic} is that we replace condition $v_j \geq \min\{p_j \Pi^\star,p_j\}$ by $v_j \geq p_j \Pi^\star$ using the fact that we made the assumption throughout that $\frac{v_j}{p_j} < 1$.}

\begin{algorithm}[ht]
\caption{DLP Policy in \cite{erdelyi2010dynamic}}\label{alg:dlp}
\begin{algorithmic}[1]
\State Initialize $x_{j}[0]=0\;\frall j$\;
 \For{$t=1,\ldots,T$}
  \State Observe type $j$ of arrival in period $t$\;
  \State \textbf{if}\;{$v_j \geq p_j \Pi^\star$:}{
  Accept arrival of type $j$ in period $t$ and set $x_j[t]=x_j[t-1]+1$
  }
  \State \textbf{else}:\;{Reject arrival of type $j$ in period $t$ and set $x_j[t]=x_j[t-1]$}
  \State Set $x_{j'}[t]=x_{j'}[t-1],\frall j'\neq j$
 \EndFor
\end{algorithmic}
\end{algorithm}

\cite{erdelyi2010dynamic} use the DLP policy to derive estimates of the fraction of type $j$ customers among all accepted customers, which is denoted $\alpha_j$. Specifically, DPD simulates the trajectory of the DLP policy under~$M$ realizations of arrivals. With $x_j^m$ denoting the number of accepted type $j$ customers over the entire time horizon in the $m$th realization, they define
$$\alpha_j = \frac{\sum_{m = 1}^M x_j^m}{\sum_{m = 1}^M \sum_{j'} x_{j'}^m}.$$

The DPD policy uses the values of $\alpha_j$ for each $j$ to approximate the expected penalty cost. Specifically, it defines $\Gamma_{\vec{\alpha}}(x)$ as the number of customers who show up, among $x$ accepted customers, when a fraction $\alpha_j$ of the accepted customers is of type $j$. Mathematically, $\Gamma_{\vec{\alpha}}(x)$ is the sum of $k$ random variables, the $j$th of which is $Bin(\left \lfloor{x \alpha_j}\right \rfloor, p_j)$ with probability $\left \lfloor{x \alpha_j}\right \rfloor +1-x \alpha_j$ and $Bin(\left \lceil{x \alpha_j}\right \rceil, p_j)$ otherwise.\footnote{Notice that $x\alpha_j$ need not be integer, and consequently, for the binomial random variables to be well-defined, we decompose $x\alpha_j$ into $\lfloor x\alpha_j\rfloor$ and $\lceil x\alpha_j\rceil$ with appropriately chosen probability.} 
With this state aggregation technique, the DPD policy uses the following DP to decide dynamically whether to accept a customer at time $t$.

\begin{equation} \label{eq:dpd_dp}
\begin{split}
& U_{t}\left(x\right)\\
=& \sum_{j = 1}^{k} \lambda_j \max\{v_j + U_{t+1}(x+1),U_{t+1}(x)\}, \\
&\text{ and } U_T(x) = -\mathbb{E}[(\Gamma_{\vec{\alpha}}(x)-B)^+]],
\end{split}
\end{equation}
In this way, it suffices for the program to know the total number of reservations to obtain the expected penalty cost and the DP becomes computationally tractable. Then, relying on a solution to \eqref{eq:dpd_dp}, \cite{erdelyi2010dynamic} formulate the DPD policy (see Algorithm \ref{alg:dpd}).

Our main result in this section, Proposition \ref{prop:lower_bound_dpd}, shows that this state aggregation technique can lead to performance losses of $\Omega(T)$ when the no-show probabilities are heterogeneous across types. Crucially, this result is not driven by the estimation technique used to derive the $\alpha_j$ values, i.e., the simulation of the DLP policy under different trajectories. Instead, it is based on the aggregation of $\alpha_j$ diluting the potentially very different value/demand requirement of customers with different no-show probabilities. For example, and this is what the proposition is based on, a customer class with small value and very high no-show probability can be preferable to one with high value and very small no-show probability. While our class of index policies appropriately trades off values and no-show probabilities, the DPD policy always prefers the higher-value customers. Finally, it is worth noting that the numerical results of \cite{erdelyi2010dynamic} do not display this problem as they are based on homogeneous no-show probabilities across classes (in which case the state aggregation is not required in the first place).

\begin{algorithm}[ht]
\caption{DPD Policy in \cite{erdelyi2010dynamic}}\label{alg:dpd}
\begin{algorithmic}[1]
\State Initialize $x_{ j}[0]=0\;\frall j$\;
\For{$t=1,\ldots,T$}
  \State Observe type $j$ of arrival in period $t$\;
  If{$v_j \geq U_{t+1}(\sum_j x_j[t-1]) - U_{t+1}(\sum_j x_j[t-1]+1)$}:{
  Accept arrival of type $j$ in period $t$ and set $x_j[t]=x_j[t-1]+1$}
  \State \textbf{else}: {Reject arrival of type $j$ in period $t$ and set $x_j[t]=x_j[t-1]$}
  \State Set $x_{j'}[t]=x_{j'}[t-1],\frall j'\neq j$
 \EndFor
\end{algorithmic}
\end{algorithm}


\begin{proposition}\label{prop:lower_bound_dpd}
When the no-show probabilities $p_1,\ldots,p_k$ are not all equal, then the loss of Algorithm \ref{alg:dpd} can be of order $\Omega(T)$.
\end{proposition}


\emph{Proof.}
Consider two fare classes, type $1$ and type $2$. Each type 1 customer generates a revenue of $r$ per customer and shows up with probability $p_1 = 1$; each type 2 customer generates a revenue of $\epsilon$ per customer and shows up with probability $p_2 = 0$. There is a penalty cost of~$1$ for each denied customer, and $1 > \frac{1}{2} = r > \epsilon > 0.$  To analyze the asymptotic loss in this set-up, we assume $B = \frac{T}{2}.$ At each period $t$, an incoming customer is of type $1$ with probability $\frac{4}{5}$ and of type $2$ with probability~$\frac{1}{5}$. Recall that we use $N_j$ for the number of arrivals of type $j$ over the entire time horizon.

From $p_2 \Pi^\star = 0 < \epsilon$ we know the DLP Policy would always accept types $2$ customers. Depending on the relative values of $r, p_1, \text{ and } \Pi^\star$ we are in one of the following two cases:
\begin{enumerate}[label=(\alph*)]
    \item $r < p_1 \Pi^\star$
    \item $r \geq p_1 \Pi^\star$
\end{enumerate}

In case (a), the DLP policy accepts only type $2$ customers, which yields $\alpha_1 = 0$ and $\alpha_2 = 1$. The DPD policy then solves for \eqref{eq:dpd_dp} assuming that no accepted customers would show up, which leads itself to accept all customers. Since $\mathbb{E}[N_1] = \frac{4}{5}T > \frac{3}{4}T$, it is trivial to show by Chernoff bound that $\mathbb{P}\left[N_1 \geq \frac{3}{4}T\right] \geq 1-e^{-\Theta(T)}.$ Conditioned on the high-probability event that $N_1 \geq \frac{3}{4}T$, both the DPD policy and the offline optimal solution accept all type $2$ customers. However, the offline optimal solution only accepts $B$ type $1$ customers, while the DPD policy accepts at least $\frac{3}{4}T - B = \frac{1}{4}T$ more type $1$ customers than the offline optimal solution does, all of whom show up and receive a compensation, which leads to an expected loss of $\Omega(T).$

We next show that in case (b) the loss of Algorithm \ref{alg:dpd} is also of order $\Omega(T)$.\footnote{While we restrict ourselves here to the $\alpha_j$'s arising from the DLP policy, the arguments extend to prove that the loss of Algorithm \ref{alg:dpd} is of order $\Omega(T)$ for any estimates of $\alpha_j$'s.} First, we derive a high probability bound on the value of $\alpha_j$'s arising from the DLP policy. Denote the number of type $j$ customers accepted by DLP over the entire time horizon by $\bar{x}_j$. Since $v_j \geq p_j \Pi^\star, \forall j$ in case (b), the DLP policy simply accepts all customers. Then, by Chernoff bound we know with probability at least $1-e^{-\Theta(T)}$ we have, with $\epsilon_1 = \frac{1}{38},\epsilon_2 = \frac{1}{19}$, that $$\bar{x}_1 \in \left[\left(\frac{4}{5}-\epsilon_1\right)T,\left(\frac{4}{5}+\epsilon_2\right)T\right] \text{ and } \bar{x}_2 \in \left[\left(\frac{1}{5}-\epsilon_1\right)T,\left(\frac{1}{5}+\epsilon_2\right)T\right].$$ Thus, we obtain $\alpha_1 = \frac{\bar{x}_1}{\bar{x}_1+\bar{x}_2} \in \left[\frac{7}{10},\frac{9}{10}\right]$ with probability greater than $1-e^{-\Theta(T)}$. DPD then assumes that a fixed portion, $\alpha_j$, of all accepted customers are of type $j$ (and thus shows up with probability~$p_j$), and solves \eqref{eq:dpd_dp}. Specifically, DPD solves the following fictional problem:
\begin{itemize}
    \item Each arrival of type $1$ has value $r$ and each arrival of type $2$ has value $\epsilon$ 
    \item At departure, a fraction $\alpha_1$ of customers show up
\end{itemize}
Note that the objective of the fictional problem given accepted customers $\vec{x}$ is either
$$rx_1+\epsilon x_2 - \left(\left \lceil(x_1+x_2)\alpha_1\right \rceil-B\right)^+ \text{ or }rx_1+\epsilon x_2 - \left(\left \lfloor(x_1+x_2)\alpha_1\right \rfloor-B\right)^+$$
with probability depending on the fractional part of $(x_1+x_2)\alpha_1$. To simplify the arithmetic we assume instead, without loss of generality, that the objective is deterministically $$I(\vec{x}):=rx_1+\epsilon x_2 - \left((x_1+x_2)\alpha_1-B\right)^+.$$ 
This changes the objective by at most $1$ and thus does not affect our asymptotic analyses.



In contrast, in the real problem, with $p_1=1$ and $p_2=0$, the objective with accepted customers~$\vec{x}$ is $$R(\vec{x}):=rx_1+\epsilon x_2 - (x_1-B)^+.$$ 
Given an arrival vector $\vec{A}$, we denote the solutions, on $\vec{A}$, to the fictional problem based on the DPD policy by $\vec{x}^{DPD}$, the offline optimal solution to the fictional problem by $\vec{x}^{fic}$ and the offline optimal solution to the real problem by $\vec{x}^{real}$ --- notice that to avoid lengthy case distinctions we assume that $\vec{x}^{fic}$ can be fractional. Our goal is to show $$\mathbb{E}_{\vec{A}}[R(\vec{x}^{real})-R(\vec{x}^{DPD})] \in \Omega(T).$$

In the rest of the proof we show
\begin{enumerate}[label=(\roman*)]
    \item $\mathbb{E}_{\vec{A}}[R(\vec{x}^{real})-R(\vec{x}^{fic})] \in \Omega(T)$,
    \item $\mathbb{E}_{\vec{A}}[|R(\vec{x}^{fic})-R(\vec{x}^{DPD})|] \in o(T)$.
\end{enumerate}

Assuming both of the statements above are true, we know $$\mathbb{E}_{\vec{A}}[R(\vec{x}^{real})-R(\vec{x}^{DPD})] \in \Omega(T) \pm o(T) = \Omega(T).$$

We shall make use of the event that
$$E := \left\{N_1 \geq \frac{3}{4}T, N_2 \geq \frac{T}{24} \text{ and } \alpha_1 \in \left[\frac{7}{10},\frac{9}{10}\right]\right\}.$$ In particular, since $\mathbb{E}[N_1] = \frac{4}{5}T > \frac{3}{4}T$, it is trivial to show by Chernoff bound that $\mathbb{P}\left[N_1 \geq \frac{3}{4}T\right] \geq 1-e^{-\Theta(T)}.$ Similarly, $\mathbb{P}\left[N_2 \geq \frac{T}{24}\right] \geq 1-e^{-\Theta(T)}.$ That is, event $E$ happens with probability greater than $1-e^{-\Theta(T)}$. All arguments from now on are conditioned on event $E$.

\subsubsection*{Step 1: Proof of (i)}
We first show that $\mathbb{E}_{\vec{A}}[R(\vec{x}^{real})-R(\vec{x}^{fic})|E] \in \Omega(T)$. Conditioning on event $E$, we know that $$R(\vec{x}^{real}) = N_2 \epsilon + \min(N_1,B) r = N_2 \epsilon + B r$$ given the realized arrival vector $\vec{A}$, where $x^{real}_1 = B$ and $x^{real}_2 = N_2 \geq \frac{T}{24}$. 

For $\vec{x}^{fic}$, note that customers of type $1$ are strictly preferred over customers of type $2$ in the fictional problem since $r > \epsilon$ and the two types are otherwise equivalent. Thus, it is optimal to accept all type $1$ customers before accepting any type $2$ customer in the fictional problem. Next, observe that, under event $E$ with $\alpha_1\geq\frac{7}{10}>\frac{1}{2}=r$, $\vec{x}^{fic}$ accepts at most~$\frac{B}{\alpha_1}$ customers in total. This is because additional customers increase the objective by, at most,~$r-\alpha_1<0$.
We can similarly argue that $\vec{x}^{fic}$ accepts at least  $\frac{B}{\alpha_1}$ customers in total, i.e., it accepts exactly $\frac{B}{\alpha_1}$ customers. Since we also know that $N_1 > \frac{B}{\alpha_1}$ and that type~$1$ customers are preferred over type~$2$ customers, we obtain $x^{fic}_1=\frac{B}{\alpha_1} < N_1$ and $x^{fic}_2 = 0.$

Therefore, the number of rejected type $2$ customers in $\vec{x}^{fic}$ is $N_2 \geq \frac{T}{24}$, and the loss of~$R(\vec{x}^{fic})$ from rejecting type $2$ customers is at least $\frac{T}{24} \epsilon$. Moreover, $\vec{x}^{real}$ accepts exactly~$B$ type 1 customers (under $E$), whereas $\vec{x}^{fic}$ accepts at least $B$ of them. With $p_1=1$, each accepted type 1 customer, beyond the first $B$, only  decreases the real objective function  $R(\cdot)$; as a result, we find that
$$\mathbb{E}_{\vec{A}}[R(\vec{x}^{real})-R(\vec{x}^{fic})|E] \geq \frac{T}{24} \epsilon = \Omega(T).$$

Finally, since $\mathbb{P}[E^c] \leq e^{-\Theta(T)}$, we obtain
$$\mathbb{E}_{\vec{A}}[R(\vec{x}^{real})-R(\vec{x}^{fic})] \geq \Omega(T) \cdot \mathbb{P}[E] + 0 \cdot (1-\mathbb{P}[E]) = \Omega(T).$$

\subsubsection*{Step 2: Proof of (ii)}
Since we know that each accepted customer leads to at most an additional revenue of 1 or compensation of 1, we have $$\mathbb{E}_{\vec{A}}[|R(\vec{x}^{fic})-R(\vec{x}^{DPD})|] \leq \mathbb{E}_{\vec{A}}\left[k\Vert \vec{x}^{fic}-\vec{x}^{DPD}\Vert_{\infty}\right].$$ Thus, to prove (ii), it suffices to show that $\mathbb{E}_{\vec{A}}\left[\Vert\vec{x}^{fic}-\vec{x}^{DPD}\Vert_{\infty}\right] \in o(T)$. We will prove this as follows: we first argue that the DPD policy obtains $O(1)$ expected loss relative to~$\vec{x}^{fic}$ to the fictional problem with objective $I(\cdot)$. We then use Markov's inequality to obtain a coarse probabilistic bound on the DPD policy having $\Omega(\sqrt{T\log(T)})$ loss. Next, we assume (see Equation \eqref{eq:omega_t}) that whenever $\vec{x}^{DPD}$ has small loss, under $I(\cdot)$, relative to~$\vec{x}^{fic}$, 
the two solutions must be close together. Worded differently: $\vec{x}^{DPD}$ cannot have small loss with respect to~$I(\cdot)$ without being close to $\vec{x}^{fic}$. We conclude the proof by first using our probabilistic bound on $\vec{x}^{DPD}$ having large loss, and then proving the statement in Equation~\eqref{eq:omega_t}.

Now, from Theorem 2.1 of \cite{freund2019uniform} we know there exists\footnote{Since $I(\vec{x})$ is deterministically a linear objective, we apply Theorem 2.1 of \cite{freund2019uniform} instead of Theorem \ref{thm:coupling} of our paper.} an online algorithm with solution $\vec{x}^{online}$ to the fictional problem such that $\mathbb{E}_{\vec{A}}\left[I(\vec{x}^{fic})-I(\vec{x}^{online})\right] \in O(1).$ Since the DPD policy provides the optimal online solution to the fictional problem, we then obtain $$\mathbb{E}_{\vec{A}}\left[I(\vec{x}^{fic})-I(\vec{x}^{DPD})\right] \leq \mathbb{E}_{\vec{A}}\left[I(\vec{x}^{fic})-I(\vec{x}^{online})\right] \in O(1),$$ i.e., there exists a constant $M > 0$ such that $$\mathbb{E}_{\vec{A}}\left[I(\vec{x}^{fic})-I(\vec{x}^{DPD})\right] \leq M.$$

We then apply Markov's inequality to find that 
$$
\mathbb{P}\left[I(\vec{x}^{fic})-I(\vec{x}^{DPD}) \geq \sqrt{T\log(T)}\Big|E\right] \cdot \mathbb{P}[E]\leq \frac{ \mathbb{E}_{\vec{A}}\left[I(\vec{x}^{fic})-I(\vec{x}^{DPD})|E\right]\mathbb{P}[E]}{\sqrt{T\log(T)}},
$$
which in turn can be bounded as
\begin{equation*}
    \begin{split}
    \frac{\mathbb{E}_{\vec{A}}\left[I(\vec{x}^{fic})-I(\vec{x}^{DPD})|E\right] \cdot \mathbb{P}[E]}{\sqrt{T\log(T)}}
    \leq & \frac{\mathbb{E}_{\vec{A}}\left[I(\vec{x}^{fic})-I(\vec{x}^{DPD})\right]}{\sqrt{T\log(T)}}
    \leq  \frac{M}{\sqrt{T\log(T)}}.
    \end{split}
\end{equation*}

That is, $$\mathbb{P}[I(\vec{x}^{fic})-I(\vec{x}^{DPD}) \geq \sqrt{T\log(T)}|E] \leq \frac{M}{\sqrt{T\log(T)}\mathbb{P}[E]}.$$

We shall show that there exists a constant $c_1 > 0$ such that
\begin{equation}\label{eq:omega_t}
    \Vert \vec{x}^{fic}-\vec{x}^{DPD}\Vert_{\infty} \geq c_1 \sqrt{T\log(T)} \text{ implies }I(\vec{x}^{fic})-I(\vec{x}^{DPD}) \geq \sqrt{T\log(T)}.
\end{equation}

From that follows the first inequality in
\begin{equation*}
    \begin{split}
    &\mathbb{P}\left[\Vert\vec{x}^{fic}-\vec{x}^{DPD}\Vert_{\infty} \geq c_1 \sqrt{T\log(T)}|E\right]\\
    \leq& \mathbb{P}\left[I(\vec{x}^{fic})-I(\vec{x}^{DPD}) \geq \sqrt{T\log(T)}|E\right]\\
    \leq& \frac{M}{\sqrt{T\log(T)}\mathbb{P}[E]}.
    \end{split}
\end{equation*}

Then, assuming \eqref{eq:omega_t}, we obtain
\begin{equation*}
    \begin{split}
    &\mathbb{E}_{\vec{A}}\left[\Vert\vec{x}^{fic}-\vec{x}^{DPD}\Vert_{\infty}\right] \\
    \leq& \mathbb{E}_{\vec{A}}\left[\Vert\vec{x}^{fic}-\vec{x}^{DPD}\Vert_{\infty}\right|E] \cdot \mathbb{P}[E] + T(1-\mathbb{P}[E])\\
    \leq & (c_1 \sqrt{T\log(T)} \cdot \mathbb{P}[\Vert\vec{x}^{fic}-\vec{x}^{DPD}\Vert_{\infty} < c_1 \sqrt{T\log(T)}|E] \\
    &+ T \cdot \mathbb{P}[\Vert\vec{x}^{fic}-\vec{x}^{DPD}\Vert_{\infty} \geq c_1 \sqrt{T\log(T)}|E] ) \mathbb{P}[E] + o(1)\\
    \leq & c_1 \sqrt{T\log(T)} + M \sqrt{\frac{T}{\log(T)}}+ o(1)\\
    \in & o(T)
    \end{split}
\end{equation*}
as required.

Now, to show \eqref{eq:omega_t}, recall from step 1 that we have $x^{fic}_1=\frac{B}{\alpha_1}<N_1,  x_2^{fic}=0$. 
We assume for simplicity that $\frac{B}{\alpha_1} \in \mathbb{Z}$ to avoid a length case distinction between $x^{DPD}_1 + x^{DPD}_2=\left \lfloor{\frac{B}{\alpha_1}}\right \rfloor$ and $x^{DPD}_1 + x^{DPD}_2=\left \lceil{\frac{B}{\alpha_1}}\right \rceil$. We now show that $x^{DPD}_1 + x^{DPD}_2 = \frac{B}{\alpha_1}$, i.e., DPD also accepts~$\frac{B}{\alpha_1}$ customers in total.   
Under event $E$, since $r < \alpha_1$, an additional customer would not be accepted by the DPD policy when $x^{DPD}_1 + x^{DPD}_2 \geq \frac{B}{\alpha_1}.$ Similarly, under event $E$ we know the optimal online policy would always accept all of the last $\left\lfloor\frac{a}{\alpha_1}\right\rfloor$ arrivals when $B-(x^{DPD}_1 + x^{DPD}_2) \alpha_1 = a > 0$. 
Thus, $x^{DPD}_1 + x^{DPD}_2 \geq \frac{B}{\alpha_1}$, i.e., $x^{DPD}_1 + x^{DPD}_2 = \frac{B}{\alpha_1} = x^{fic}_1 + x^{fic}_2.$

Now, take $c_1 = \frac{1}{r-\epsilon}$. We know from $$x^{DPD}_1 + x^{DPD}_2 = \frac{B}{\alpha_1} = x^{fic}_1 + x^{fic}_2 \text{ and }\Vert \vec{x}^{fic}-\vec{x}^{DPD}\Vert_{\infty} \geq c_1 \sqrt{T\log(T)}$$ that $$x^{fic}_1 - x^{DPD}_1 \geq c_1 \sqrt{T\log(T)} \text{ and } x^{DPD}_2 - x^{fic}_2 \geq c_1 \sqrt{T\log(T)}.$$ Then, $$I(\vec{x}^{fic})-I(\vec{x}^{DPD}) \geq c_1 \sqrt{T\log(T)}(r-\epsilon) \geq \sqrt{T\log(T)},$$ which completes the proof of (ii). \halmos

\end{document}